\input amstex\documentstyle {amsppt}  
\pagewidth{12.5 cm}\pageheight{19 cm}\magnification\magstep1
\topmatter
\title Character sheaves on disconnected groups, II\endtitle
\author G. Lusztig\endauthor
\address Department of Mathematics, M.I.T., Cambridge, MA 02139\endaddress
\thanks Supported in part by the National Science Foundation\endthanks
\endtopmatter   
\document
\define\da{\dagger}
\define\ucf{\un{\cf}}
\define\uS{\un S}

\define\Lie{\text{\rm Lie }}

\define\si{\sim}

\define\sqc{\sqcup}

\define\qua{\quad}

\define\hX{\hat X}
\define\hY{\hat Y}
\define\hf{\hat f}
\define\hboc{\hat\boc}

\define\tct{\ti\ct}
\define\tcm{\ti\cm}
\define\tcf{\ti\cf}

\define\tis{\ti\s}
\define\tbc{\ti\boc}

\define\bj{\bar j}

\define\bY{\bar Y}

\define\bpi{\bar\p}
\define\bsi{\bar\s}
\define\bboc{\bar\boc}
\define\lb{\linebreak}

\define\op{\oplus}

\define\em{\emptyset}
\define\imp{\implies}
\define\ra{\rangle}

\define\m{\mapsto}
\define\do{\dots}
\define\la{\langle}
\define\bsl{\backslash}

\define\lra{\leftrightarrow}

\define\sub{\subset}
\define\bxt{\boxtimes}
\define\T{\times}
\define\ti{\tilde}
\define\nl{\newline}
\redefine\i{^{-1}}
\define\fra{\frac}
\define\un{\underline}

\define\ot{\otimes}
\define\bbq{\bar{\QQ}_l}

\define\Ad{\text{\rm Ad}}
\define\Hom{\text{\rm Hom}}
\define\End{\text{\rm End}}

\define\Ker{\text{\rm Ker}}

\define\supp{\text{\rm supp}}

\define\a{\alpha}
\redefine\b{\beta}
\redefine\c{\chi}
\define\g{\gamma}
\redefine\d{\delta}
\define\e{\epsilon}

\redefine\o{\omega}
\define\p{\pi}
\define\ph{\phi}

\define\r{\rho}
\define\s{\sigma}
\redefine\t{\tau}

\redefine\l{\lambda}

\redefine\D{\Delta}

\define\Si{\Sigma}

\define\Ph{\Phi}

\define\boc{\bold c}
\define\dd{\bold d}

\define\kk{\bold k}

\redefine\AA{\bold A}

\define\EE{\bold E}
\define\FF{\bold F}

\define\NN{\bold N}

\define\QQ{\bold Q}

\define\ZZ{\bold Z}

\define\ca{\Cal A}

\define\cc{\Cal C}
\define\cd{\Cal D}
\define\ce{\Cal E}
\define\cf{\Cal F}

\define\ch{\Cal H}
\define\ci{\Cal I}

\define\cl{\Cal L}
\define\cm{\Cal M}
\define\cn{\Cal N}

\define\cs{\Cal S}
\define\ct{\Cal T}

\define\cv{\Cal V}
\define\cw{\Cal W}
\define\cz{\Cal Z}
\define\cx{\Cal X}
\define\cy{\Cal Y}

\define\fg{\frak g}

\define\fD{\frak D}

\define\fZ{\frak Z}

\define\tg{\ti g}

\define\tp{\ti p}

\define\tA{\ti A}

\define\tD{\ti D}

\define\tG{\ti G}

\define\tL{\ti L}
\define\tM{\ti M}

\define\tX{\ti X}
\define\tY{\ti Y}

\define\Irr{\text{\rm Irr }}

\define\bg{\bar g}

\define\bS{\bar S}

\define\bce{\bar\ce}
\define\tce{\ti\ce}
\define\BBD{BBD}
\define\BO{B}
\define\HS{HS}

\define\IC{L2}
\define\COX{L6}
\define\CLA{L7}
\define\CSII{L8}
\define\AD{L9}
\define\LSS{LS1}
\define\LS{LS2}
\define\MA{M1}
\define\MAA{M2}
\define\MS{MS}
\define\SP{Sp}
\define\SPR{Spr}
 \head Introduction \endhead
This paper is a part of a series (beginning with \cite{\AD}) which attempts to
develop a theory of character sheaves on a not necessarily connected reductive
algebraic group $G$. The main theme of this paper is establishing the "generalized
Springer correspondence" in complete generality. 

More precisely, we consider the set $\cn$ consisting of all pairs $(\boc,\cf)$
where $\boc$ is a unipotent $G^0$-conjugacy class in $G$ and $\cf$ is an irreducible
local system on $\boc$, equivariant for the conjugation action of $G^0$. We define a 
finite sequence of finite Coxeter groups $W^1,W^2,\do,W^k$ and we establish a canonical
bijection $\cn\lra\sqc_i\Irr W^i$ where $\Irr W^i$ is the set of isomorphism classes of 
irreducible representations of $W^i$. (This is done in 11.10 after preliminary work in 
\S7-\S11.) This extends a bijection established in \cite{\IC} in which only pairs 
$(\boc,\cf)$ with $\boc\sub G^0$ were considered. (If the characteristic $p$ of the
ground field is $0$, the condition $\boc\sub G^0$ is automatically satisfied, but if 
$p>1$ this is not so.) This also extends the original Springer 
correspondence \cite{\SPR} in which not only $\boc$ is assumed to be contained in $G^0$ but $\cf$ 
is subject to certain restrictions. The methods we use are generalizations of those in
\cite{\IC}, although a number of new technical difficulties must be overcome. On the 
other hand, the utilization of Deligne's theory of weights makes some of our proofs 
simpler than the corresponding proofs in \cite{\IC}. Now $\cn$ has a natural partition 
into blocks corresponding to the obvious partition of $\sqc_i\Irr W^i$ into pieces 
indexed by $i$. Of particular interest are the pairs $(\boc,\cf)$ which form a block 
by themselves (the corresponding $W^i$ is trivial). Such pairs (which could be called
cuspidal) are classified in \S12. In \S13 we describe combinatorially the generalized 
Springer correspondence for disconnected classical groups in characteristic $2$ in 
analogy with the method of \cite{\LS} (for some partial results in this direction, 
see \cite{\MS}). \S14 is a complement to the discussion in \cite{\LS, \S4} of the 
generalized Springer correspondence for Spin groups in characteristic $\ne 2$; namely 
it expresses that correspondence by a closed combinatorial formula instead of the 
inductive procedure of \cite{\LS}

We adhere to the notation of \cite{\AD}. In particular, $G$ denotes a fixed, not 
necessarily connected, reductive algebraic group over an algebraically closed field 
$\kk$. 

Here is some additional notation. The cardinal of a finite set $F$ is denoted by $|F|$.
Let $\nu$ or $\nu_G$ be the number of positive roots of $G^0$. The closure of a subset 
$S$ of an algebraic variety $V$ is denoted by $\bS$. The dual of a local system $\ce$ 
is denoted by $\ce^*$. We generally ignore Tate twists. For an algebraic variety $X$,
let $\fD:\cd(X)@>>>\cd(X)$ be Verdier duality. 

References to \cite{\AD} are given without specifying \cite{\AD}; for example, we 
write 1.7(a) instead on \cite{\AD, 1.7(a)}.

\head Contents\endhead
7. Sheaves on the variety of quasi-semisimple classes.

8. Study of local systems on unipotent conjugacy classes.

9. A restriction theorem.

10. Preparatory results.

11. The structure of the algebra $\EE$.

12. Classification of objects in $\cn_D^0$.

13. Symbols.

14. Spin groups.

\head 7. Sheaves on the variety of quasi-semisimple classes\endhead
\subhead 7.1\endsubhead
Let $D$ be a connected component of $G$.  Let $D//G^0$ be the set of closed 
$G^0$-conjugacy classes in $D$, that is, the set of quasi-semisimple $G^0$-conjugacy 
classes in $D$, see 1.4(c). By geometric invariant theory,  $D//G^0$ has a natural 
structure of affine variety and there is a well defined morphism $\s:D@>>>D//G^0$ 
such that, for $g\in D$, $\s(g)$ is the unique closed $G^0$-conjugacy class in $D$ 
contained in the closure of the $G^0$-conjugacy class of $g$. Let $g\in D$ be 
quasi-semisimple and let $T_1$ be a maximal torus of $Z_G(g)^0$. Then $\s$ induces a 
morphism 

(a) $gT_1@>>>D//G^0$ 
\nl
which is a finite (ramified) covering inducing a bijection between the set of orbits 
of the action of the finite group $\{n\in G^0;ngT_1n\i=gT_1\}/T_1$ on $gT_1$ and 
$D//G^0$, see 1.14(a)-(d).

In our case, $\s$ can be described explicitly at follows. Let $g\in D$. By 1.9 
(applied to $Z_G(g_s)$ instead of $G$), we can find $u\in g_uZ_G(g_s)^0$ such that 
$u$ is unipotent, quasi-semisimple in $Z_G(g_s)$; then $g_su\in D$ is quasi-semisimple 
in $G$ (see 1.4(c)) and $\s(g)$ is the $G^0$-conjugacy class of $g_su$ (this is 
independent of the choice of $u$ since $u$ is unique up to $Z_G(g_s)^0$-conjugacy.) To
verify that $\s(g)$ is as stated above, we must show that $g_su$ is in the closure of 
the $G^0$-conjugacy class of $g=g_sg_u$. It is enough to show that $u$ is in the 
closure of the $Z_G(g_s)^0$-conjugacy class of $g_u$; this follows from 
\cite{\SP, II, 2.21} applied to $Z_G(g_s)$ instead of $G$. Clearly $\s$ is constant on 
$G^0$-conjugacy classes in $D$. 

\medpagebreak

In the case where $D$ contains unipotent elements, there is a unique point
$\o\in D//G^0$ such that $\s\i(\o)=\{g\in D;g\text{ unipotent}\}$.

\medpagebreak

The following result will be used several times in this paper. Let $Z=Z'\cup Z''$ be 
a partition of a variety $Z$ with $Z'$ open and $Z''$ closed. Let $\cf$ be a local 
system on $Z$. Assume that $\dim Z\le n$. Then

(b) {\it the natural sequence 
$0@>>>H_c^{2n}(Z',\cf)@>>>H_c^{2n}(Z,\cf)@>>>H_c^{2n}(Z'',\cf)@>>>0$ is exact.}
\nl
We can find an open subset $U$ of $Z$ such that $\dim(Z-U)<n$ and $U\cap Z''$ is 
open in $U$. Then 
$$\align&H_c^{2n}(Z',\cf)=H_c^{2n}(Z'\cap U,\cf),H_c^{2n}(Z,\cf)=H_c^{2n}(Z\cap U,\cf),
\\&H_c^{2n}(Z'',\cf)=H_c^{2n}(Z''\cap U,\cf)\endalign$$
and it is enough to show that the natural sequence 

$0@>>>H_c^{2n}(Z'\cap U,\cf)@>>>H_c^{2n}(Z\cap U,\cf)@>>>H_c^{2n}(Z''\cap U,\cf)
@>>>0$ 
\nl
is exact. This holds since $Z'\cap U,Z''\cap U$ form a partition of $Z\cap U$ into
subsets that are both open and closed.

\subhead 7.2\endsubhead
Let $(L,S)\in\AA$ with $S\sub D$ and let $P$ be a parabolic of $G$ with Levi $L$ 
such that $S\sub N_GP$. Let $Y'_{L,S}$ be as in 3.14. Let $\d$ be the connected
component of $N_GL\cap N_GP$ that contains $S$. Let ${}^\d\cz_L,S^*$ be as in 3.11.
Let $A_{L,S}=\s(Y'_{L,S})$. We show:

(a) {\it if $g\in S^*,g'\in\bS U_P$ are such that $\s(g')=\s(g)$ then there exists
$u\in U_P$ such that $ug'u\i\in\bS,L(ug'u\i)=L$.}
\nl
Since $g'_s\in N_GP$ is semisimple, it normalizes some Levi of $P$ that is, some 
$U_P$-conjugate of $L$. Hence replacing $g',x$ by $u'{}\i g'u',xu'$ for some 
$u'\in U_P$ we may assume in addition that $g'_s\in N_GL\cap N_GP$. We have $g'=hv$ 
where $h\in\bS,v\in U_P$. As in the proof of Lemma 5.5 we see that $g'_s=h_s$ and 
$T(g')=T(h)$. Hence $L(g')=L(h)$. Since $h$ is isolated in $N_GL\cap N_GP$ (Lemma 
2.8) we have $L\sub L(h)$ (see 3.8(b)) and $L\sub L(g')$. 

Let $v\in g_uZ_G(g_s)^0$ be such that $v$ is unipotent, quasi-semisimple in 
$Z_G(g_s)$; let $v'\in g'_uZ_G(g'_s)^0$ be such that $v'$ is unipotent, 
quasi-semisimple in $Z_G(g'_s)$. The assumption $\s(g')=\s(g)$ implies that there 
exists $x\in G^0$ with $xg'_sv'x\i=g_sv$, hence $xg'_sx\i=g_s,xv'x\i=v$. Now 
$xv'x\i\in xg'_ux\i Z_G(g_s)^0$ hence $v\in xg'_ux\i Z_G(g_s)^0$. It follows that 
$xg'_ux\i Z_G(g_s)^0=g_uZ_G(g_s)^0$ that is, $xg'_ux\i=g_u\mod Z_G(g_s)^0$. Thus,
$$\align&T(xg'x\i)=(\cz_{Z_G(xg'_sx\i)^0}\cap Z_G(xg'_ux\i))^0=(\cz_{Z_G(g_s)^0}
\cap Z_G(xg'_ux\i))^0\\&=(\cz_{Z_G(g_s)^0}\cap Z_G(g_u))^0=T(g)\endalign$$
hence $L(xg'x\i)=L(g)$. Since $g\in S^*$, we have $L(g)=L$ hence $L(xg'x\i)=L$ and
$L(g')=x\i Lx$, $\dim L(g')=\dim(x\i Lx)=\dim L$. This, combined with $L\sub L(g')$ 
implies $L=L(g')$. Since $g'\in N_G(L(g'))$, we have $g'\in N_GL$. We have also 
$g'=hv$ where $h\in S\sub N_GL,v\in U_P$. Since $hv\in N_GL$ we have 
$v\in N_GL\cap U_P$ hence $v=1$ and $g'=h\in\bS$. This proves (a).

Exactly the same proof gives the following variant of (a):

(b) {\it Let $g\in S^*,g'\in SU_P$ be such that $\s(g')=\s(g)$. Then there exists 
$u\in U_P$ such that $ug'u\i\in S,L(ug'u\i)=L$, hence $ug'u\i\in S^*$.}

For future reference we note the following result.

(c) {\it Let $(L',S')\in\AA$. The following two conditions are equivalent:

(i) $Y_{L',S'}$ is contained in the closure of $Y_{L,S}$;

(ii) there exists $x\in G^0$ such that $xLx\i\sub L'$ and $S'$ is contained in the 
closed subset $T_x=\cup_{y\in L'}y\bS U_{Q_x}y\i$ of $N_GL'$; here $Q_x=xPx\i\cap L'$, 
a parabolic of $L'$ with Levi $xLx\i$.}

{\it Proof of (i)$\imp$(ii).} Let $g\in S'{}^*$. We will show that 

{\it there exists $x\in G^0$ such that $xLx\i\sub L'$ and $g\in T_x$.}
\nl
(This would imply that $S'\sub T_x$ since $T_x$ is stable under $L'$-conjugacy and 
under left translation by $\cz_{L'}^0$.) From (i) we deduce $g\in Y'_{L,S}$ (see Lemma 
3.14). Replacing $(L,S,P)$ by a $G^0$-conjugate, we may assume that $g=hu$ with 
$h\in\bS,u\in U_P$. Replacing $(L,S)$ by a $U_P$-conjugate, we may further assume that 
$g_s=h_s\in N_GL\cap N_GP$, $T_G(g)\sub(\cz_L\cap Z_G(h))^0$ (see the proof of
Lemma 3.15). It follows that $Z_{G^0}((\cz_L\cap Z_L(h))^0)\sub Z_{G^0}(T(g))$ hence 
$L\sub L'$. (We use 1.10, 3.9 and that $g\in S'{}^*$.) It is enough to show that 
$g\in T_1$. 

Let $Q'=P\cap L'$. This is a parabolic of $L'$ with Levi $L$ (since $L\sub L'$) and
$U_{Q'}=U_P\cap L'$. There is a unique parabolic $P'$ of $G^0$ with Levi $L'$ such that
$P\sub P'$. Since $g$ normalizes $L'$ (recall that $g\in S'{}^*$) and $P$ (since 
$g\in\bS U_P$) we have $g\in N_GP'$. Since $h\in\bS$ we have $h\in N_GL\cap N_GP$. Now 
$huP'u\i h\i=P'$ hence $h\i P'h=uP'u\i=P'$ (recall that $u\in U_P\sub P\sub P'$). Also 
$h\in N_GP$. Hence $h\i P'h$ and $P'$ both contain $P$ and are $G^0$-conjugate (we have
$u\i(h\i P'h)u=P',u\in G^0$) hence $h\i P'h=P'$. Now $h\in N_GP'\cap N_GL$ hence 
$h\in N_GL'$ (since $L'$ is the unique Levi of $P'$ that contains $L$). We see that 
$h\in N_GL'\cap N_GP'$. Since $hu\in N_GL'\cap N_GP'$ it follows that
$u\in N_GL'\cap P'=L'$. Thus $u\in L'\cap U_P=U_{Q'}$. We see that $g\in\bS U_Q$ hence 
$g\in T_1$. Thus (ii) holds.

{\it Proof of (ii)$\imp$(i).} We may assume that $L\sub L'$ and 
$S'\sub\cup_{y\in L'}y\bS U_Qy\i$ where $Q=P\cap L'$ (a parabolic of $L'$ with Levi 
$L$). Since $U_Q\sub U_P$, we have $S'\sub Y'_{L,S}$. Since $Y'_{L,S}$ is stable under 
$G^0$-conjugacy, it follows that $Y_{L',S'}\sub Y'_{L,S}$.

\proclaim{Lemma 7.3}(a) $A_{L,S}=\s(S)$.

(b) $A_{L,S}$ is an irreducible, closed subvariety of $D//G^0$ of dimension
$\dim{}^\d\cz_L^0$.

(c) Let $\cy=\{g\in Y'_{L,S};\dim L(g)>\dim L\}$. Then $\cy$ is closed in $Y'_{L,S}$
and $\{g\in Y'_{L,S};\dim L(g)=\dim L\}$ is open in $Y'_{L,S}$.

(d) $\{g\in Y'_{L,S};\dim L(g)=\dim L\}$ is exactly the inverse image of $\s(S^*)$ 
under $\s:Y'_{L,S}@>>>A_{L,S}$.

(e) The inverse image of $\s(S^*)$ under $\s:Y'_{L,S}@>>>A_{L,S}$ is open in 
$Y'_{L,S}$.

(f) $\s(S^*)$ is open in $A_{L,S}$.
\endproclaim
We prove (a). Let $g\in Y'_{L,S}$. We must show that $\s(g)\in\s(S)$. As in the 
proof of 3.15 we may assume that $g=hv$ where $h\in\bS,v\in U_P$ and $g_s=h_s$; 
moreover, by 1.22 we have $h_s=h'_s$ with $h'\in S,h'{}\i h\in Z_G(h_s)^0$ hence 
$h'_u{}\i h_u\in Z_G(g_s)^0$. We have $g_u=h_uv$ hence 
$v\in U_P\cap Z_G(g_s)=U_P\cap Z_G(g_s)^0$. Thus, 
$h'_u{}\i g_u=h'_u{}\i h_uv\in Z_G(g_s)^0$. Choose 
$u\in g_uZ_G(g_s)^0=h_uZ_G(h_s)^0=h'_uZ_G(h_s)^0$ such that $u$ is unipotent, 
quasi-semisimple in $Z_G(g_s)$. Then $\s(g)=\s(h)=\s(h')$ is the $G^0$-conjugacy 
class of $g_su=h_su=h'_su$ and $\s(g)\in\s(S)$, as required.

We prove (b). Let $h\in S$. We have $S=\cup_{x\in L}x{}^\d\cz_L^0hx\i$ hence 
$\s(S)=\{\s(zh);z\in{}^\d\cz_L^0\}$. Let $u\in h_uZ_L(h_s)^0$ be such that $u$ is 
unipotent quasi-semisimple in $Z_{N_GL\cap N_GP}(h_s)$. For any $z\in{}^\d\cz_L^0$ 
we have $Z_L(h_s)=Z_L(zh_s)\sub Z_G(zh_s)$ hence $u\in h_uZ_G(zh_s)^0$. Moreover, 
$u$ is unipotent quasi-semisimple in $Z_{N_GL\cap N_GP}(zh_s)$ hence $u$ is 
unipotent quasi-semisimple in $Z_G(zh_s)$. It follows that 

$\s(zh)=\s(zh_su)=\s(zg)$
\nl
where $g=h_su\in D$ is quasi-semisimple in $G$. We see that 

$\s(S)=\{\s(zg);z\in{}^\d\cz_L^0\}$. 
\nl
Let $T_1$ be a maximal torus of $Z_G(g)^0$. Now
${}^\d\cz_L^0g$ is an irreducible closed subvariety of $gT_1$ of dimension 
$\dim{}^\d\cz_L^0$; since $\s:gT_1@>>>D//G^0$ is a finite (ramified) covering (see 
7.1), it follows that $A_{L,S}=\s(S)=\s({}^\d\cz_L^0g)$ is an irreducible closed 
subvariety of $D//G^0$ of dimension $\dim{}^\d\cz_L^0$. This proves (b). 

We prove (c). By 3.15 there exist $(L_i,S_i),i\in[1,n]$, such that the closure of
$Y_{L,S}$ in $G$ equals $\cup_{i=1}^nY_{L_i,S_i}$ and $Y_{L_i,S_i}$ are disjoint. If
$g\in Y_{L_i,S_i}$ then $\dim L(g)=\dim L_i$. Let 
$J=\{i\in[1,n];\dim L_i>\dim L\}$. Then $\cy=\cup_{i\in J}Y_{L_i,S_i}$. Again by 
3.15, the closure of $Y_{L_i,S_i}$ in $G$ is $\cup_{j\in J_i}Y_{L_j,S_j}$ for some 
subset $J_i\sub[1,n]$. By 7.2(c), for $j\in J_i$ we have $\dim L_j\ge\dim L_i$. Hence 
if $i\in J,j\in J_i$ we have $\dim L_j\ge\dim L_i>\dim L$ so that $J_i\sub J$. We 
see that for $i\in J$, the closure of $Y_{L_i,S_i}$ in $G$ is contained in $\cy$. Thus
$\cy$ is closed. This implies the last assertion of (c) since 
$\dim L(g)\ge\dim L$ for $g\in Y'_{L,S}$, see 7.2(c). This proves (c).

We prove (d). We must show that, if $g\in Y'_{L,S}$, then the following two 
conditions are equivalent:

(i) $\dim L(g)=\dim L$;

(ii) $\s(g)=\s(g_1)$ for some $g_1\in S^*$.
\nl
Assume that (ii) holds. We may assume that $g\in\bS U_P$. By 7.2(a) we have 
$\dim L(g)=\dim L$. Thus (i) holds. Assume now that (i) holds. As in the proof of 
3.15 we may assume that $g=hv$ where $h\in\bS,v\in U_P$ and $g_s=h_s$; moreover, by
1.22 we have $h_s=h'_s$ with $h'\in S,h'{}\i h\in Z_G(h_s)^0$ hence 
$h'_u{}\i h_u\in Z_G(g_s)^0$. As in the proof of (a) we have $\s(g)=\s(h')$. As in 
the proof of 3.15 we have $T(g)=T(h)$. Hence $L(h)=L(g)$. Since 
$h_s=h'_s,h'{}\i h\in Z_G(h_s)^0$ we have $T(h)=T(h')$ (see 2.1(e)) hence 
$L(h)=L(h')$. Thus $L(h')=L(g)$ and $\dim L(h')=\dim L$. Since $h'$ is isolated in
$N_GL\cap N_GP$ we have $L\sub L(h')$ (see 3.8). It follows that $L=L(h')$. Since 
$h'\in S$ we have $h'\in S^*$. Thus (ii) holds. This proves (d).

(e) follows by combining (c),(d). Now a subset of $A_{L,S}$ is open in $A_{L,S}$ if
and only if its inverse image under $\s:Y'_{L,S}@>>>A_{L,S}$ is open in $Y'_{L,S}$.
(We may identify $A_{L,S}$ with the variety $Y'_{L,S}//G^0$ of closed 
$G^0$-conjugacy classes in $Y'_{L,S}$.) Hence (f) follows from (e). The lemma is 
proved.

\subhead 7.4\endsubhead
In the remainder of this section we fix parabolics $P',P''$ of $G^0$ with Levi 
$L',L''$ respectively and an isolated stratum $S'$ (resp. $S''$) of 
$\tL'=N_GL'\cap N_GP'$ (resp. $\tL''=N_GL''\cap N_GP''$) contained in the connected 
component $\d'$ (resp, $\d''$) of $\tL'$ (resp. $\tL''$) with 
$\d'\sub D,\d''\sub D$. We also fix $\ce'\in\cs(S')$, $\ce''\in\cs(S'')$ such that 
$(S',\ce'),(S'',\ce'')$ are cuspidal pairs of $\tL',\tL''$. Let
$$d_0=2\nu-\nu_{L'}-\nu_{L''}+\fra{1}{2}(\dim({}^{\d'}\cz_{L'}^0\bsl S')+
\dim({}^{\d''}\cz_{L''}^0\bsl S'')).$$
Let $\bY'$ (resp. $\bY''$) be the closure of $Y_{L',S'}$ (resp. 
$Y_{L'',S''}$) in $D$. As in 3.14, let 
$$\align&X'=\{(g,x'P')\in G\T G^0/P';x'{}\i gx'\in\bS'U_{P'}\},\\&
X''=\{(g,x''P'')\in G\T G^0/P'';x''{}\i gx''\in\bS''U_{P''}\}.\endalign$$
Define $\ph':X'@>>>\bY',\ph'':X''@>>>\bY''$ by $\ph'(g,x'P')=g,\ph''(g,x''P'')=g$.
For any stratum $S'_1$ of $\tL'$ contained in $\bS'$ let 

$X'_{S'_1}=\{(g,x'P')\in G\T G^0/P';x'{}\i gx'\in S'_1U_{P'}\}\sub X'$.
\nl
For any stratum $S''_1$ of $\tL''$ contained in $\bS''$ let 

$X''_{S''_1}=\{(g,x''P'')\in G\T G^0/P'';x''{}\i gx''\in S''_1U_{P''}\}\sub X''$. 
\nl
Let $\bce'$ (resp. $\bce''$) be the local system on $X'_{S'}$ (resp. $X''_{S''}$) 
defined in terms of $\ce'$ (resp. $\ce''$) in the same way as $\bce$ was defined in
terms of $\ce$ in 5.6. Let 

$K'=IC(X',\bce')\in\cd(X')$, $K''=IC(X'',\bce'')\in\cd(X'')$.
\nl
Let 
$$\fZ=\{(g,x'P',x''P'')\in G\T G^0/P'\T G^0/P'';x'{}\i gx'\in\bS'U_P,
x''{}\i gx''\in\bS''U_{P''}\}.$$
Let $A'=A_{L',S'},A''=A_{L'',S''},A=A'\cap A''\sub D//G^0$. Define 
$\tis:\fZ@>>>A$ by $\tis(g,x'P',x''P'')=\s(g)$. For $a\in A$ let 
$\fZ^a=\tis\i(a)\sub\fZ$. For $S'_1,S''_1$ as above we set
$$\align&\fZ_{S'_1,S''_1}=\{(g,x'P',x''P'')\in G\T G^0/P'\T G^0/P'';
x'{}\i gx'\in S'_1U_{P'},\\&x''{}\i gx''\in S''_1U_{P''}\},\endalign$$
$$\align\fZ^a_{S'_1,S''_1}=&\{(g,x'P',x''P'')\in G\T G^0/P'\T G^0/P'';\\&\s(g)=a,
x'{}\i gx'\in S'_1U_{P'},x''{}\i gx''\in S''_1U_{P''}\}.\endalign$$
The $\fZ^a_{S'_1,S''_1}$ form a partition of $\fZ^a$ into locally closed 
subvarieties with $\fZ^a_{S',S''}$ open.

\subhead 7.5\endsubhead
We show:
$$\align\dim\fZ^a_{S'_1,S''_1}\le d_0&-\fra{1}{2}(\dim({}^{\d'}\cz_{L'}^0\bsl S')
-\dim({}^{\d'}\cz_{L'}^0\bsl S'_1))\\&-\fra{1}{2}(\dim({}^{\d''}\cz_{L''}^0\bsl S'')
-\dim({}^{\d''}\cz_{L''}^0\bsl S''_1)).\tag a\endalign$$
Define $m:\fZ^a_{S'_1,S''_1}@>>>\s\i(a)$ by $m(g,x'P',x''P'')=g$. Now $\s\i(a)$ is 
a union of finitely many $G^0$-conjugacy classes (since the semisimple parts of 
elements of $\s\i(a)$ lie in a single $G^0$-conjugacy class). Hence it is enough to
estimate $\dim(m\i(\boc))$ for any $G^0$-conjugacy class $\boc$ contained in 
$\s\i(a)$. Any fibre of $m:m\i(\boc)@>>>\boc$ is isomorphic to a product of two 
varieties of the kind appearing in 4.4(b). Hence by 4.4(b), this fibre has dimension 
at most
$$(\nu-\fra{1}{2}\dim\boc)-(\nu_{L'}-\fra{1}{2}\dim({}^{\d'}\cz_{L'}^0\bsl S'_1))
+(\nu-\fra{1}{2}\dim\boc)-(\nu_{L''}
-\fra{1}{2}\dim({}^{\d''}\cz_{L''}^0\bsl S''_1))$$
and (a) follows.

From (a) we deduce immediately:
$$\dim\fZ^a\le d_0.\tag b$$ 
The inverse image of $K'\bxt K''\in\cd(X'\T X'')$ under 
$\fZ@>>>X'\T X'',(g,x'P',x''P'')\m((g,x'P'),(g,x''P''))$ is denoted again by 
$K'\bxt K''$. Its restrictions to various subvarieties of $\fZ$ are denoted again 
by $K'\bxt K''$. Similarly, the inverse image of the local system $\bce'\bxt\bce''$
under $\fZ_{S',S''}@>>>X'_{S'}\T X''_{S''},(g,x'P',x''P'')\m((g,x'P'),(g,x''P''))$ 
is denoted again by $\bce'\bxt\bce''$. Its restrictions to various subvarieties of 
$\fZ_{S',S''}$ are denoted again by $\bce'\bxt\bce''$. Let 

$\tct=\ch^{2d_0}\tis_!(K'\bxt K'')$, $\ct=\ch^{2d_0}\s^1_!(\bce'\bxt\bce'')$ 
\nl
where $\s^1:\fZ_{S',S''}@>>>A$ is the restriction of $\tis$; these are constructible 
sheaves on $A$. 

\proclaim{Lemma 7.6} The natural homomorphism $\ct@>>>\tct$ is an isomorphism.
\endproclaim
It is enough to show that, for any $a\in A$, the natural homomorphism
$$H^{2d_0}_c(\fZ^a_{S',S''},\bce'\bxt\bce'')=H^{2d_0}_c(\fZ^a_{S',S''},K'\bxt K'')
@>>>H^{2d_0}_c(\fZ^a,K'\bxt K'')$$
is an isomorphism. First we show:
$$\text{ if $(S'_1,S''_1)\ne(S',S'')$ and $H^t_c(\fZ^a_{S'_1,S''_1},K'\bxt K'')\ne 
0$ then }t<2d_0.\tag a$$
We may assume that $S'_1\ne S'$. (The case where $S''_1\ne S''$ is similar.) From 
the hypercohomology spectral sequence of $K'\bxt K''$ on $\fZ^a_{S'_1,S''_1}$ we see
that there exist $i,j',j''$ such that $t=i+j'+j''$ and
$H^i_c(\fZ^a_{S'_1,S''_1},\ch^{j'}K'\bxt\ch^{j''}K'')\ne 0$. It follows that
$i\le 2\dim\fZ^a_{S'_1,S''_1}$ and $X'_{S'_1}\sub\supp\ch^{j'}K'$,
$X''_{S''_1}\sub\supp\ch^{j''}K''$. (Note that $\ch^{j'}K'$ is a local system of
constant rank on $X'_{S'_1}$ and $\ch^{j''}K''$ is a local system of constant rank
on $X''_{S''_1}$.) Using 7.5(a) we have
$$i\le 2d_0-\dim({}^{\d'}\cz_{L'}^0\bsl S')+\dim({}^{\d'}\cz_{L'}^0\bsl S'_1)
-\dim({}^{\d''}\cz_{L''}^0\bsl S'')+\dim({}^{\d''}\cz_{L''}^0\bsl S''_1).$$
Moreover, by the definition of $K',K''$ we have
$$j'<\dim X'-\dim X'_{S'_1}=\dim({}^{\d'}\cz_{L'}^0\bsl S')-
\dim({}^{\d'}\cz_{L'}^0\bsl S'_1),$$
$$j''\le\dim X''-\dim X''_{S''_1}=\dim({}^{\d''}\cz_{L''}^0\bsl S'')-
\dim({}^{\d''}\cz_{L''}^0\bsl S''_1)$$
(since $S'_1\ne S'$). It follows that $t=i+j'+j''<2d_0$. This proves (a).

From (a) we see that $H^t_c(\fZ^a_{S'_1,S''_1},K'\bxt K'')=0$ for any
$(S'_1,S''_1)\ne(S',S'')$ and any $t\ge 2d_0$; hence
$$H^t_c(\fZ^a-\fZ^a_{S',S''},K'\bxt K'')=0 \text{ for any } t\ge 2d_0.\tag b$$
As a part of the long cohomology exact sequence of the partition
$\fZ^a=\fZ^a_{S',S''}\cup(\fZ^a-\fZ^a_{S',S''})$ we have the exact sequence
$$H^{2d_0-1}_c(\fZ^a-\fZ^a_{S',S''},K'\bxt K'')@>f>>
H^{2d_0}_c(\fZ^a_{S',S''},K'\bxt K'')@>>>H^{2d_0}_c(\fZ^a,K'\bxt K'')@>>>0,$$
where the last $0$ comes from (b). It is enough to prove that $f=0$. We may assume 
that $\kk$ is an algebraic closure of a finite field $\FF_q$, that $G$ has a fixed 
$\FF_q$-structure with Frobenius map $F:G@>>>G$, that $L',P',L'',P''$ are 
$F$-stable, any $S'_1\sub\bS',S''_1\sub\bS''$ as above is $F$-stable, that $F(a)=a$
and that we have isomorphisms $F^*\ce'@>\si>>\ce',F^*\ce''@>\si>>\ce''$ which make 
$\ce',\ce''$ into local systems of pure weight $0$. Then we have natural 
(Frobenius) endomorphisms of $H^{2d_0-1}_c(\fZ^a-\fZ^a_{S',S''},K'\bxt K'')$,
$H^{2d_0}_c(\fZ^a_{S',S''},K'\bxt K'')$ compatible with $f$. To show that $f=0$, it
is enough to show that 

(c) $H^{2d_0}_c(\fZ^a_{S',S''},K'\bxt K'')$ is pure of weight $2d_0$;

(d) $H^{2d_0-1}_c(\fZ^a-\fZ^a_{S',S''},K'\bxt K'')$ is mixed of weight $\le 2d_0-1$.
\nl
Now $H^{2d_0}_c(\fZ^a_{S',S''},K'\bxt K'')=H^{2d_0}_c(\fZ^a_{S',S''},
\bce'\bxt\bce'')$ and $\dim\fZ^a_{S',S''}\le d_0$; (c) follows. Using the partition
$\fZ^a-\fZ^a_{S',S''}=\cup_{(S'_1,S''_1)\ne(S',S'')}\fZ^a_{S'_1,S''_1}$, we see 
that to prove (d), it is enough to prove:

if $(S'_1,S''_1)\ne(S',S'')$ then $H^{2d_0-1}_c(\fZ^a_{S'_1,S''_1},K'\bxt K'')$ is 
mixed of weight $\le 2d_0-1$.
\nl
Using again the hypercohomology spectral sequence of $K'\bxt K''$ on 
$\fZ^a_{S'_1,S''_1}$ we see that it is enough to prove:

(e) if $i,j',j''$ are such that $2d_0-1=i+j'+j''$ then 
$H^i_c(\fZ^a_{S'_1,S''_1},\ch^{j'}K'\bxt\ch^{j''}K'')$ is mixed of weight 
$\le 2d_0-1$.
\nl
By Gabber's theorem \cite{\BBD, 5.3.2}, the local system $\ch^{j'}K'$ on 
$X'_{S'_1}$ is mixed of weight $\le j'$ and the local system $\ch^{j''}K''$ on 
$X''_{S''_1}$ is mixed of weight $\le j''$. Hence the local system 
$\ch^{j'}K'\bxt\ch^{j''}K''$ on $\fZ^a_{S'_1,S''_1}$ is mixed of weight 
$\le j'+j''$. Using Deligne's theorem \cite{\BBD, 5.1.14(i)}, we deduce that 

$H^i_c(\fZ^a_{S'_1,S''_1},\ch^{j'}K'\bxt\ch^{j''}K'')\ne 0$
\nl
is mixed of weight $\le i+j'+j''$. This proves (e). The lemma is proved.

\subhead 7.7\endsubhead
Let $E$ be a locally closed subset of $G^0/P'\T G^0/P''$ which is a union of 
$G^0$-orbits (for the diagonal $G^0$-action). We set
$$\align{}^E\fZ_{S',S''}&=\{(g,x'P',x''P'')\in G\T G^0/P'\T G^0/P'';
x'{}\i gx'\in S'U_{P'},\\&x''{}\i gx''\in S''U_{P''},(x'P',x''P'')\in E\}.\endalign$$
Let ${}^E\ct=\ch^{2d_0}\s^1_!(\bce'\bxt\bce'')$ (a constructible sheaf on $A$) where 
$\s^1:{}^E\fZ_{S',S''}@>>>A$ is $(g,x'P',x''P'')\m\s(g)$.

If $E$ is a $G^0$-orbit, we say that $E$ is good if for some/any $(x'P',x''P'')\in E$, 
$x'P'x'{}\i,x''P''x''{}\i$ have a common Levi; we say that $E$ is bad if it is not 
good. 

\proclaim{Lemma 7.8}(a) If $E$ is a bad orbit then ${}^E\ct=0$;

(b) if $E$ is a good orbit and there is no $n\in G^0$ such that $(P',nP'')\in E$ and
$n\i L'n=L'',n\i S'n=S''$ then ${}^E\ct=0$;

(c) if $E$ is a good orbit and there exists $n\in G^0$ such that $(P',nP'')\in E$ 
and $n\i L'n=L'',n\i S'n=S''$ (so that $A'=A''=A$) then $\fD({}^E\ct)[-2\dim A]$ is
a constructible sheaf on $A$.
\endproclaim
The fibre of the fibration $pr_{23}:{}^E\fZ_{S',S''}@>>>E$ at 
$(x'P',x''P'')\in E$ where $x',x''\in G^0$ may be identified with

$V=\{g\in G;x'{}\i gx'\in S'U_{P'},x''{}\i gx''\in S''U_{P''}\}$.
\nl
Define $j:V@>>>S'\T S''$ by 
$$j(g)=(\text{$S'$-component of $x'{}\i gx'$, $S''$-component of $x''{}\i gx''$}).$$
Let ${}^E\ct'=\ch^{2d_0-2\dim E}\s'_!(j^*(\ce'\bxt\ce''))$ where $\s':V@>>>A$ is 
$g\m\s(g)$. Let 
${}^E\ct''=\ch^{2d_0+2\dim H}\s''_!(\bbq\bxt j^*(\ce'\bxt\ce''))$ where 
$H=x'P'x'{}\i\cap x''P''x''{}\i$ and $\s'':G^0\T V@>>>A$ is $(x,g)\m\s(g)$. Using 
the spectral sequence attached to the composition $G^0\T V@>pr_2>>V@>\s'>>A$ (equal
to $\s''$) and that attached to the composition 
$G^0\T V@>>>H\bsl(G^0\T V)={}^E\fZ_{S',S''}@>\s^1>>A$ (equal to $\s)$ we obtain 
${}^E\ct''={}^E\ct$, ${}^E\ct''={}^E\ct'$. (We use that 
$\ch^i\s^1_!(\bce'\bxt\bce'')=0$ for $i>2d_0$ and 
$\ch^i\s'_!(j^*(\ce'\bxt\ce''))=0$ for $i>2d_0-2\dim E$.) It follows that 
${}^E\ct={}^E\ct'$.

Let $Q'=x'P'x'{}\i,Q''=x''P''x''{}\i$. Choose Levi subgroups $M'$ of $Q'$ and $M''$
of $Q''$ such that $M',M''$ contain a common maximal torus. Let 
$\tM'=N_GM'\cap N_GQ'$, $\tM''=N_GM''\cap N_GQ''$. Let $\Si'$ (resp. $\Si''$) be the 
unique stratum of $\tM'$ (resp. $\tM''$) such that $\Si'U_{Q'}=x'S'U_{P'}x'{}\i$ 
(resp. $\Si''U_{Q''}=x''S''U_{P''}x''{}\i$). Now $\Ad(x')$ (resp. $\Ad(x'')$) carries
the inverse image of $\ce'$ (resp. $\ce''$) under $pr_1:S'U_{P'}@>>>S'$ 
(resp. $pr_1:S''U_{P''}@>>>S''$) to a local system on $\Si'U_{Q'}$ (resp. 
$\Si''U_{Q''}$) which is the inverse image under $pr_1:\Si'U_{Q'}@>>>\Si'$ (resp. 
$pr_1:\Si''U_{Q'}@>>>\Si''$) of a local system $\cf'$ (resp. $\cf''$) on $\Si'$ (resp. 
$\Si''$). If $g\in V$ then $g\in N_GQ'\cap N_GQ''$. By 1.25(a),(b) we can write 
uniquely $g=zu''u=zu'v$ where 
$$z\in\tM'\cap\tM'',u''\in M'\cap U_{Q''},u\in U_{Q'}\cap Q'',u'\in M''\cap U_{Q'},
v\in U_{Q''}\cap Q'.$$
Thus $V$ can be identified with
$$\align&\{(u,v,u'',u',z)\\&\in(U_{Q'}\cap Q'')\T(U_{Q''}\cap Q')\T(M'\cap U_{Q''})
\T(M''\cap U_{Q'})\T(\tM'\cap\tM'');\\&u''u=u'v,zu''\in\Si',zu'\in\Si''\}.\endalign$$
The map $\s':V@>>>A$ becomes $(u,v,u'',u',z)\m\s(z)$. (We use that 
$h\in\tM',v\in U_{Q'}\imp\s(hv)=\s(v)$ and the analogous fact for $\tM'',U_{Q''}$ 
instead of $\tM',U_{Q'}$; see the argument in the proof of 7.3(a).) In this 
description, $V$ is fibred over
$$V_1=\{(u'',u',z)\in(M'\cap U_{Q''})\T(M''\cap U_{Q'})\T(\tM'\cap\tM'');
zu''\in\Si',zu'\in\Si''\}$$
with all fibres isomorphic to $U_{Q'}\cap U_{Q''}$ (see the argument in 4.2). The 
map $\s':V@>>>A$ factors through a map $\bsi:V_1@>>>A,(u'',u',z)\m\s(z)$. Since 
$U_{Q'}\cap U_{Q''}$ is an affine space of dimension
$2\nu-\nu_{L'}-\nu_{L''}-\dim E$, we see that 
$\ct'{}^E=\ch^r\bsi_!(\bj^*(\cf'\bxt\cf''))$ where 
$r=\dim({}^{\d'}\cz_{L'}^0\bsl S')+\dim({}^{\d''}\cz_{L''}^0\bsl S'')$ and 
$\bj:V_1@>>>\Si'\T\Si''$ is defined by $(u'',u',z)\m(zu'',zu')$. For any $a\in A$ let
$V_1^a=\bsi\i(a)\sub V_1$.

Assume that we are in the setup of (a). It is enough to show that for any $a\in A$,
$H^r_c(V_1^a,\bj^*(\cf'\bxt\cf''))=0$. Let $\tp_3:V_1^a@>>>\tM'\cap\tM''$ be the 
third projection. By an argument in the proof of Lemma 4.2, the image of $\tp_3$ is
a union of finitely many $(M'\cap M'')$-conjugacy classes $\e_1,\e_2,\do,\e_t$ in 
$\tM'\cap\tM''$. Since $\dim V_1\le\fra{1}{2}r$, it is enough to show that for any 
$z\in\e_i$ we have
$$H^{r-2\dim\e_i}_c(\tp_3\i(z),\bj^*(\cf'\bxt\cf''))=0.$$
Now $\tp_3\i(z)$ is a product $R'\T R''$ where $R'$ is the set of all elements in 
$(\tM'\cap N_GQ'')\cap\Si'$ whose image under $\tM'\cap N_GQ''@>>>\tM'\cap\tM''$ is 
equal to $z$ and $R''$ is the set of all elements in $(\tM''\cap N_GQ')\cap\Si''$ 
whose image under $\tM''\cap N_GQ'@>>>\tM'\cap\tM''$ is equal to $z$. Since 
$$\align&2\dim(R')\le d'=\dim({}^{\d'}\cz_{L'}^0\bsl S')-\dim\e_i,\\&
2\dim(R'')\le d''=\dim({}^{\d''}\cz_{L''}^0\bsl S'')-\dim\e_i\endalign$$
(see 4.2(a)) and $d'+d''=r-2\dim\e_i$, we are reduced to showing that
$$H^{d'}_c(R',\cf')\ot H^{d''}_c(R'',\cf'')=0.$$
Since $Q',Q''$ have no common Levi, we see that either $M'\cap Q''$ is a proper 
parabolic of $M'$ or $M''\cap Q'$ is a proper parabolic of $M''$. In the first case
we have $H^{d'}_c(R',\cf')=0$ since $(\Si',\cf')$ is a cuspidal pair for $\tM'$. In 
the second case we have $H^{d''}_c(R'',\cf'')=0$ since $(\Si'',\cf'')$ is a cuspidal
pair for $\tM''$. Thus (a) is proved.

Assume that we are in the setup of (b). It is enough to show that for any $a\in A$,
$H^r_c(V_1^a,\bj^*(\cf'\bxt\cf''))=0$. If this is not so, we can find $a\in A$ such
that $V_1^a\ne\em$. Now $Q',Q''$ have a common Levi hence 
$M'=M'',M'\cap U_{Q''}=\{1\},M''\cap U_{Q'}=\{1\}$ and we may identify $V_1^a$ with 
$\{z\in\Si'\cap\Si'';\s(z)=a\}$. Thus, $\Si'\cap\Si''\ne\em$. Since $\Si',\Si''$ are 
strata of $N_GM'=N_GM''$ with non-empty intersection, we must have $\Si'=\Si''$. 

Since $Q'=x'P'x'{}\i$, we see that $x'{}\i M'x'$ is a Levi of $P'$. We can find 
$v'\in U_{P'}$ such that $L'=v'{}\i x'{}\i M'x'v'$. We have 
$S'U_{P'}=x'{}\i\Si'x'U_{P'}=v'{}\i x'{}\i\Si'x'v'U_{P'}$. Since 
$S'\sub L',v'{}\i x'{}\i\Si'x'v'\sub L'$ it follows that 
$S'=v'{}\i x'{}\i\Si'x'v'$.
\nl
Similarly we can find $v''\in U_{P''}$ such that 

$L''=v''{}\i x''{}\i M''x''v''$, $S''=v''{}\i x''{}\i\Si''x''v''$.
\nl
Since $M'=M'',\Si'=\Si''$, we have

$x'v'L'v'{}\i x'{}\i=x''v''L''v''{}\i x''{}\i$,
$x'v'S'v'{}\i x'{}\i=x''v''S''v''{}\i x''{}\i$.
\nl
Thus $n\i L'n=L'',n\i S'n=S''$, 
where $n=v'{}\i x'{}\i x''v''\in G^0$, $(P',nP'')\in E$. This contradicts the 
assumption of (b). Thus, (b) is proved.

Assume that we are in the setup of (c). Let $n\in G^0$ be such that 
$(P',nP'')\in E$, $n\i L'n=L'',n\i S'n=S''$. Taking $x'=1,x''=n$ in the arguments 
above, we have $Q'=P',Q''=nP''n\i$ and we can take $M'=M''=L',\Si'=\Si''=S'$; we have 
${}^E\ct=\ch^r\bsi_!(\cf)$ where $\bsi:S'@>>>A$ is the restriction of $\s$, 
$\cf=\ce'\ot\Ad(n\i)^*\ce''$, $r=2\dim({}^{\d'}\cz_{L'}^0\bsl S')$.

Let $\tL'=N_GL'\cap N_GP'$. Replacing $G,D$ in the definition of $D//G^0,\s$ in 7.1
by $\tL',\d'$ we obtain an affine variety $\d'//L'$ whose points are the 
quasi-semisimple $L$-conjugacy classes in $\d'$ and a morphism 
$\s_0:\d'@>>>\d'//L'$. Now any quasi-semisimple $L$-conjugacy class in $\d'$ is 
contained in a unique quasi-semisimple $G^0$-conjugacy class in $D$. Thus we obtain
a map $\p:\d'//L'@>>>D//G^0$ which is in fact a morphism of affine varieties. 

We show that $\p$ is a finite (ramified) covering. Let $g\in\d'$ be 
quasi-semisimple in $\tL'$ and let $T_1$ be a maximal torus of $Z_{L'}(g)^0$. By an
argument similar to that in the proof of 1.12(b) we see that $T_1$ is also a 
maximal torus of $Z_G(g)^0$. By 7.1, $\s$ defines a finite (ramified) covering 
$gT_1@>>>D//G^0$; similarly, $\s_0$ defines a finite (ramified) covering 
$gT_1@>>>\d'//L'$. Clearly, $gT_1@>\s>>D//G^0$ is the composition 
$gT_1@>\s_0>>\d'//L'@>\p>>D//G^0$ and our assertion follows.

Let $\uS'=\s_0(S')$. Now ${}^{\d'}\cz_{L'}^0$ acts on $\d'//L'$ by left 
multiplication (with finite isotropy groups) and this action is compatible under 
$\s_0$ with the action of ${}^{\d'}\cz_{L'}^0$ on $\d'$ by left multiplication. 
Since the ${}^{\d'}\cz_{L'}^0$ action permutes transitively the $L'$-conjugacy 
classes in $S'$, we see that $\uS'$ is a single ${}^{\d'}\cz_{L'}^0$-orbit in 
$\d'//L'$; in particular, it is a smooth, locally closed subvariety of $\d'//L'$. 
Now $\bsi:S'@>>>A$ is a composition $S'@>\bsi_0>>\uS'@>\bpi>>A$ where $\bsi_0,\bpi$
are the restrictions of $\s_0,\p$.

We can find $t\ge 1$ such that $\cf$ is equivariant for the action $z:h\m z^th$ of
${}^{\d'}\cz_{L'}^0$ on $S'$ and this action induces an action of 
${}^{\d'}\cz_{L'}^0$ on $\uS'$ such that $\bsi_0:S'@>>>\uS'$ is 
${}^{\d'}\cz_{L'}^0$-equivariant. Hence the constructible sheaf 
$\ucf=\ch^r(\bsi_0)_!(\cf)$ on $\uS'$ is ${}^{\d'}\cz_{L'}^0$-equivariant. Since 
this action (which depends on $t$) is transitive on $\uS'$, we see that $\ucf$ is a
local system on $\uS'$. Since $S'$ is smooth, we have 
$\fD(\ucf)=\ucf^*[2\dim\uS']$. Now $\bpi$ is a finite (ramified) covering (since 
$\p$ is). Hence $\bpi_!$ commutes with Verdier duality and we have 
$\ch^r\bsi_!(\cf)=\bpi_!\ucf$. Hence 
$\fD(\bpi_!\ucf)=\bpi_!\ucf^*[2\dim\uS']=\bpi_!\ucf^*[2\dim A]$. Here $\bpi_!\ucf^*$ 
is a constructible sheaf on $A$ since $\bpi$ is finite. This proves (c).

\proclaim{Lemma 7.9} (a) $\fD(\ct)[-2\dim A]$ is a constructible sheaf on $A$.

(b) If $(L',S'),(L'',S'')$ are not conjugate under some element of $G^0$ then 
$\ct=0$.
\endproclaim
More generally, let $E$ be a locally closed subset of $G^0/P'\T G^0/P''$ which is a 
union of $G^0$-orbits (for the diagonal $G^0$-action). We show that 
$\fD({}^E\ct)[-2\dim A]$ is a constructible sheaf on $A$ and that, in the setup of 
(b), ${}^E\ct=0$. We argue by induction on the number $N$ of $G^0$-orbits contained
in $E$. If $N=0$ the result is trivial. If $N=1$ the result follows from Lemma 7.8.
Assume now that $N\ge 2$. We can find a partition $E=E'\cup E''$ so that $E',E''$ 
are $G^0$-invariant subsets of $E$ with $E'$ closed in $E$, $E''$ open in $E$ and 
$E''$ is a single $G^0$-orbit. By the induction hypothesis the result holds for $E'$
and $E''$. We have a natural exact sequence of constructible sheaves 
$$0@>>>\e''@>>>\e@>>>\e'@>>>0\tag c$$
where $\e''={}^{E''}\ct,\e={}^E\ct,\e'={}^{E'}\ct$. (To prove exactness it is enough
to show that for any $a\in A$, the natural sequence
$$\align&0@>>>H^{2d_0}_c({}^{E''}\fZ_{S',S''}\cap\fZ^a_{S',S''},\bce'\bxt\bce'')@>>>\\&
H^{2d_0}_c({}^{E}\fZ_{S',S''}\cap\fZ^a_{S',S''},\bce'\bxt\bce'')@>>>
H^{2d_0}_c({}^{E'}\fZ_{S',S''}\cap\fZ^a_{S',S''},\bce'\bxt\bce'')@>>>0\endalign$$
is exact. This follows from 7.1(b).)

Now (c) can be regarded as a distinguished triangle $(\e'',\e,\e')$ in $\cd(A)$. 
Hence the Verdier duals $(\fD\e',\fD\e,\fD\e'')$ form a distinguished triangle in 
$\cd(A)$ and we have a long cohomology exact sequence of constructible sheaves 

$\do@>>>\ch^i(\fD\e')@>>>\ch^i(\fD\e)@>>>\ch^i(\fD\e'')@>>>\do$
\nl
on $A$. Hence for each $i$ we have an exact sequence 

$\ch^i(\fD\e'[-2\dim A])@>>>\ch^i(\fD\e[-2\dim A])@>>>\ch^i(\fD\e''[-2\dim A])$. 
\nl
By the induction hypothesis, $\ch^i(\fD\e'[-2\dim A])=0,\ch^i(\fD\e''[-2\dim A])=0$
for $i\ne 0$. Hence $\ch^i(\fD\e[-2\dim A])=0$ for $i\ne 0$. This shows that 
$\fD\e[-2\dim A]$ is a constructible sheaf. In the setup of (b), we see from (c) and
the induction hypothesis that ${}^E\ct=0$. The lemma is proved.

\subhead 7.10\endsubhead
In the remainder of this section we assume that $L'=L''=L,S'=S''=S,P'=P''=P$ and 
that $\ce'=\ce''{}^*=\ce$ is irreducible. Let $\p:\tY_{L,S}@>>>Y_{L,S}$ be as in 
3.13. Recall that $\p$ is a principal $\cw_S$-bundle where 
$\cw_S=\{n\in N_{G^0}L;nSn\i=S\}/S$. (The $\cw_S$-action on $\tY_{L,S}$ is 
$w\m f_w,(g,xL)\m(g,xn_w\i L)$.) For $w\in\cw_S$ let $n_w\in N_{G^0}L$ be a 
representative of $w$ and let $\EE_w=\Hom(\Ad(n_w)^*\ce,\ce)$. (This vector space is
canonically defined, independent of the choice of $n_w$, by the $L$-equivariance of 
$\ce$.) Let $\cw_{\ce}=\{w\in\cw_S;\Ad(n_w)^*\ce\cong\ce\}$ (a subgroup of $\cw_S$).
We have 
$$\dim\EE_w=1 \text{ if $w\in\cw_{\ce}$; $\dim\EE_w=0$ if } w\in\cw_S-\cw_{\ce}.$$
Let $\tce$ be the local system on $\tY_{L,S}$ defined in 5.6. Recall that 
$a^*\tce=b^*\ce$ where $a(g,x)=(g,xL),b(g,x)=x\i gx$ in the diagram 
$\tY_{L,S}@<a<<\hY@>b>>S$ with $\hY=\{(g,x)\in G\T G^0;x\i gx\in S^*\}$. Define 
$\hf_w:\hY@>>>\hY$ by $(g,xL)\m(g,xn_w\i L)$. Then 
$a\hf_{n_w}=f_wa,b\hf_w=\Ad(n_w)b$. 

Let $\EE$ be the (semisimple) algebra of all endomorphisms of the local system 
$\p_!\tce$ on $Y_{L,S}$. We have canonically
$$\EE=\op_{w\in\cw_{\ce}}\EE_w.\tag a$$
Indeed, since $a,b$ are fibrations with smooth connected fibres, we have 
$$\align&\EE=\Hom(\p_!\tce,\p_!\tce)=\Hom(\p^*\p_!\tce,\tce)
=\op_{w\in\cw_S}\Hom(f_w^*\tce,\tce)\\&=\op_{w\in\cw_S}\Hom(a^*f_w^*\tce,a^*\tce)=
\op_{w\in\cw_S}\Hom(\hf_w^*a^*\tce,a^*\tce)\\&
=\op_{w\in\cw_S}\Hom(\hf_w^*b^*\ce,b^*\ce)
=\op_{w\in\cw_S}\Hom(b^*\Ad(n_w)^*\ce,b^*\ce)\\&=\op_{w\in\cw_S}
\Hom(\Ad(n_w)^*\ce,\ce)=\op_{w\in\cw_S}\EE_w=\op_{w\in\cw_{\ce}}\EE_w.\endalign$$
In particular, $\EE_w$ may be identified with a subspace of $\EE$. We see that
$\dim\EE=|\cw_\ce|$. 

\proclaim{Proposition 7.11} For $w\in\cw_S$ let $E_w$ be the $G^0$-orbit on 
$G^0/P\T G^0/P$ which contains $(P,n_wPn_w\i)$. There is a canonical isomorphism 
$\ct\cong\op_{w\in\cw_S}{}^{E_w}\ct$ of sheaves over $A=A_{L,S}$.
\endproclaim
Recall that $\s^1:\fZ_{S,S}@>>>A$ is given by $(g,x'P,x''P)\m\s(g)$. Let 
$\fZ'_{S,S}$ be the inverse image under $\s^1$ of the open subset $\s(S^*)$ of $A$ 
(see Lemma 7.3(f)). We have
$$\align\fZ'_{S,S}=&\{(g,x'P,x''P)\in G\T G^0/P\T G^0/P;\\&x'{}\i gx'\in SU_P,
x''{}\i gx''\in SU_P,\s(g)\in\s(S^*)\}.\endalign$$
Using 7.2(b) we see that the condition $\s(g)\in\s(S^*)$ is equivalent to the 
condition $g\in Y_{L,S}$, so that
$$\align\fZ'_{S,S}=&\{(g,x'P,x''P)\in G\T G^0/P\T G^0/P;\\&g\in Y_{L,S},
x'{}\i gx'\in SU_P,x''{}\i gx''\in SU_P\}.\endalign$$
Using Lemma 5.5 we see that $\fZ'_{S,S}$ may be identified with
$$\align\cy&
=\tY_{L,S}\T_{Y_{L,S}}\tY_{L,S}\\&=\{(g,x'L,x''L)\in G\T G^0/L\T G^0/L;
x'{}\i gx'\in S^*,x''{}\i gx''\in S^*\}.\endalign$$
Hence we have $\ct|_{\s(S^*)}=\ch^{2d_0}\s^2_!(\tce\bxt\tce^*)$ where 
$\s^2:\cy@>>>\s(S^*)$ is $(g,x'L,x''L)\m\s(g)$ and $\tce$ is the local systems on
$\tY_{L,S}$ defined in 5.6. Similarly, for any $G^0$-orbit $E$ on $G^0/P\T G^0/P$, 
we have ${}^E\ct|_{\s(S^*)}=\ch^{2d_0}\s^{2,E}_!(\tce\bxt\tce^*)$ where 
$$\align{}^E\cy=&\{(g,x'L,x''L)\in G\T G^0/L\T G^0/L;\\&x'{}\i gx'\in S^*,
x''{}\i gx''\in S^*,(x'P,x''P)\in E\}\endalign$$
and $\s^{2,E}:{}^E\cy@>>>\s(S^*)$ is $(g,x'L,x''L)\m\s(g)$.

Since $\tY_{L,S}@>>>Y_{L,S}$ is a principal bundle with group $\cw_S$, we have a 
finite partition $\cy=\sqc_{w\in\cw_S}{}^w\cy$ into open and closed subsets, where
${}^w\cy=\{(g,x'L,x''L)\in\cy;x''L=x'n_wL\}$. 

It is clear that ${}^w\cy={}^{E_w}\cy$. It follows that ${}^E\cy=\em$ for any 
$G^0$-orbit $E$ on $G^0/P\T G^0/P$ not of the form $E_w$. We see that
$\ch^{2d_0}\s^2_!(\tce\bxt\tce^*)$ is canonically isomorphic to
$\op_{w\in\cw_S}\ch^{2d_0}\s^{2,E}_!(\tce\bxt\tce^*)$. Thus we obtain a canonical 
isomorphism 
$$\ct|_{\s(S^*)}@>\si>>\op_{w\in\cw_S}{}^{E_w}\ct|_{\s(S^*)}.\tag c$$ 
By Lemma 7.9 and its proof, $\ct$ and $\op_{w\in\cw_S}{}^{E_w}\ct$ are intersection 
cohomology complexes on the irreducible variety $A$. Since $\s(S^*)$ is open dense 
in $A$, the isomorphism (c) extends uniquely to an isomorphism
$\ct@>\si>>\op_{w\in\cw_S}{}^{E_w}\ct$. The proposition is proved.
 
\subhead 7.12\endsubhead
We set $X=X'=X'',\bY=\bY'=\bY'',\ph=\ph'=\ph''$. We set $K=K',K^*=K''$. Since 
$\ph_!K=IC(\bY,\p_!\tce)$ (see 5.7) we have $\EE=\End(\p_!\tce)=\End(\ph_!K)$. In 
particular $\ph_!K$ is naturally an $\EE$-module and $\ph_!K\ot\ph_!K^*$ is 
naturally an $\EE$-module (with $\EE$ acting on the first factor). This induces an 
$\EE$-module structure on $\ch^{2d_0}\s_!(\ph_!K\ot\ph_!K^*)=\tct$ (here 
$\s:\bY@>>>A$ is the restriction of $\s:D@>>>D//G^0$). Hence for any $a\in A$ we 
obtain an $\EE$-module structure on the stalk $\tct_a$.

\proclaim{Lemma 7.13}Let $w\in\cw_{\ce}$ and let $b_w$ be a basis element of 
$\EE_w$. Multiplication by $b_w$ in the $\EE$-module structure of 
$\tct=\ct=\op_{w'\in\cw_S}{}^{E_w}\ct$ (see Lemma 7.6, Proposition 7.11) defines for
any $w'\in\cw_S$ an isomorphism ${}^{E_{w'}}\ct@>\si>>{}^{E_{ww'}}\ct$.
\endproclaim
Since ${}^{E_{w'}}\ct,{}^{E_{ww'}}\ct$ are intersection cohomology complexes on $A$ 
(see Lemma 7.9 and its proof) and $\s(S^*)$ is open dense in $A$ (see Lemma 7.3), it
is enough to prove that, for any $a\in\s(S^*)$, multiplication by $b_w$ defines an 
isomorphism of stalks ${}^{E_w'}\ct_a@>\si>>{}^{E_{ww'}}\ct_a$. Let $\cy$ be as in 
7.11 and let $\cy^a=\{(g,x'L,x''L)\in\cy;\s(g)=a\}$. As in 7.11 we have a partition
$\cy=\sqc_{w'\in\cw_S}{}^{w'}\cy$ hence a partition 
$\cy^a=\sqc_{w'\in\cw_S}{}^{w'}\cy^a$ where ${}^{w'}\cy^a={}^{w'}\cy\cap\cy^a$ are 
both open and closed in $\cy^a$. For any $w'\in\cw_S$ we have
${}^{E_w'}\ct_a=H^{2d_0}_c({}^{w'}\cy^a,\tce\bxt\tce^*)$ and we must prove that 
multiplication by $b_w$ in the $\EE$-module structure of 
$$H^{2d_0}_c(\cy^a,\tce\bxt\tce^*)=
\op_{w'\in\cw_S}H^{2d_0}_c({}^{w'}\cy^a,\tce\bxt\tce^*)$$
defines for any $w'\in\cw_S$ an isomorphism
$$H^{2d_0}_c({}^{w'}\cy^a,\tce\bxt\tce^*)@>\si>> 
H^{2d_0}_c({}^{ww'}\cy^a,\tce\bxt\tce^*).\tag a$$
(The $\EE$-module structure on
$$H^{2d_0}_c(\cy^a,\tce\bxt\tce^*)
=H^{2d_0}_c(Y_{L,S}\cap\s\i(a),\p_!\tce\ot\p_!\tce^*)$$
is induced by the $\EE$-module structure on $\p_!\tce$; recall that
$\EE=\End(p_!\tce)$.) We have an isomorphism 
$$f:{}^{w'}\cy^a@>\si>>{}^{ww'}\cy^a,(g,x'L,x''L)\m(g,x'n_w\i L,x''L).$$
Since $w\in\cw_{\ce}$, $b_w$ defines an isomorphism 
$f^*(\tce\ot\tce^*)\cong\tce\ot\tce^*$ and this induces an isomorphism (a) which is
just multiplication by $b_w$. The lemma is proved.

\subhead 7.14\endsubhead
We preserve the setup of 7.12.
Let $\Irr\EE$ be a set of representatives for the isomorphism classes of simple
$\EE$-modules. Given a semisimple object $M$ of some abelian category such that $M$
is an $\EE$-module we shall write $M_\r=\Hom_\EE(\r,M)$ for $\r\in\Irr\EE$ and we 
have $M=\op_{\r\in\Irr\EE}(\r\ot M_\r)$ with $\EE$ acting only on the $\r$-factor 
and where $M_\r$ is in our abelian category. In particular we have
$\p_!\tce=\op_{\r\in\Irr\EE}\r\ot(\p_!\tce)_\r$ where $(\p_!\tce)_\r$ is an
irreducible local system on $Y_{L,S}$. We have 
$\ph_!K=\op_{\r\in\Irr\EE}\r\ot(\ph_!K)_\r$ where 
$(\ph_!K)_\r=IC(\bY,(\p_!\tce)_\r)$; moreover, for $a\in A$ we 
have $\tct_a=\op_{\r\in\Irr\EE}\r\ot(\tct_a)_\r$.

\subhead 7.15\endsubhead
Let $a\in A$. We set $\bY^a=\{g\in\bY;\s(g)=a\}$, $X^a=\ph\i(\bY^a)\sub X$. 
Let $\ph^a:X^a@>>>\bY^a$ be the restriction of $\ph:X@>>>\bY$. Let 
$S^a=S\cap\s\i(a),\bS^a=\bS\cap\s\i(a)$. For a stratum $S_1$ of $N_GL\cap N_GP$
contained in $\bS$ let $S_1^a=S_1\cap\s\i(a)$. Let $\ce^a=\ce|_{S^a}$. Let 
$\bce^a$ be the restriction of $\bce$ to $X_S^a=X^a\cap X_S$.

\proclaim{Lemma 7.16} (a) $S^a$ is an open dense smooth (non-empty) subset of 
$\bS^a$ of pure dimension $\dim S-\dim A$ and $IC(\bS^a,\ce^a)$ is the restriction 
of $IC(\bS,\ce)$ to $\bS^a$.

(b) Both $X^a,\bY^a$ have pure dimension $d_0=2\nu-2\nu_L+\dim S-\dim A$.

(c) $X_S^a$ is an open dense smooth subset of $X^a$ hence $K^a:=IC(X^a,\bce^a)$ is 
well defined. We have $K^a=K|_{X^a}$.

(d) We have $(\ph_!K)|_{\bY^a}=\ph^a_!K^a$. Moreover, $\ph^a_!K^a[d_0]$ is a 
semisimple perverse sheaf on $\bY^a$.

(e) We have $(\ph_!K)_\r|_{\bY^a}\ne 0$ for any $\r\in\Irr\EE$.
\endproclaim
We prove (a). $S^a$ is non-empty by Lemma 7.3(a). As in the proof of Lemma 7.8(c), 
the (surjective) map $\s:S@>>>A$ is a composition $S@>\bsi_0>>\uS@>\bpi>>A$. 
Moreover, by Lemma 7.3(a) (for $N_GL\cap N_GP$ instead of $G$), the morphism 
$\bsi_0:S@>>>\uS$ extends uniquely to a morphism $\bsi_1:\bS@>>>\uS$. Since $\bpi$ 
is finite, the set $F:=\bpi\i(a)$ is finite. We have 
$\bS^a=\bsi_1\i(F),S^a=\bsi_1\i(F)\cap S$. Since $\bsi_1$ is equivariant for an 
action of a torus which acts transitively on $\uS$ (see the proof of Lemma 7.8(c)) 
we see that (a) holds.

We prove (b). We have $X^a=\{(g,xP)\in G\T G^0/P;x\i gx\in\bS^aU_P\}$. The second 
projection makes $X^a$ into a fibration over $G^0/P$ with all fibres isomorphic to 
$\bS^aU_P$ hence $X^a$ has pure dimension as indicated. (We use (a).) From 
$\ph(X)=\bY$ we deduce $\ph(X^a)=\bY^a$ hence $\dim\bY^a\le\dim X^a$. Assume that 
$\bY^a$ has some irreducible component of dimension $e<\dim X^a$. Then that 
component contains an open dense subset $U$ such that $\dim\ph\i(g)=\dim X^a-e$ for
all $g\in U$. The fibre product $X^a\T_{\bY^a}X^a$ contains the fibre product 
$\ph\i(U)\T_U\ph\i(U)$  hence it has dimension 
$\ge\dim U+2(\dim X^a-e)=2\dim X^a-e$. By 7.5(b) we have 
$\dim(X^a\T_{\bY^a}X^a)\le 2\nu-2\nu_L+\dim S-\dim A=\dim X^a$. It follows that 
$2\dim X^a-e\le\dim X^a$ so that $e\ge\dim X^a$ contradicting $e<\dim X^a$. This 
proves (b).

To prove the first assertion of (c) it is enough to show that
$\{(g,x)\in G\T G^0;x\i gx\in S^aU_P\}$ is an open dense smooth subset of 
$X'=\{(g,x)\in G\T G^0;x\i gx\in\bS^aU_P\}$. This follows from (a). Consider the 
diagram 
$X^a@<a''_0<<X'@>b''_0>>\bS^a$ where $a''_0,b''_0$ are restrictions of $a'',b''$ in 
5.6(a). We have
$$\align a''_0{}^*IC(X^a,\bce^a)&=b''_0{}^*IC(\bS^a,\ce^a)=(b''{}^*IC(\bS,\ce))|_{X'}
\\&=(a''{}^*IC(X,\bce))|_{X'}=a''_0{}^*(IC(X,\bce)|_{X^a}.\endalign$$
(The second equality follows from (b); the third equality follows from 5.6(b).) 
Since $a''$ is a principal $P$-bundle it follows that 
$IC(X^a,\bce^a)=IC(X,\bce)|_{X^a}$ and (c) is proved. 

The first assertion of (d) follows immediately from (c). Since $\ph^a$ is proper, to
show that $(\ph_!K)|_{\bY^a}[d_0]$ is a perverse sheaf, it suffices to prove that:

(f) for any $i\ge 0$ we have $\dim\supp\ch^i(\ph^a_!K^a)\le\dim\bY^a-i$
\nl
and also the analogous assertion in which $K^a$ is replaced by $IC(X^a,(\bce^a)^*)$,
which is entirely similar to (f). As in the proof of Proposition 5.7 it is enough to 
prove that for any stratum $S_1$ of $N_GL\cap N_GP$ contained in $\bS$ we have
$$\dim\{g\in\bY^a;\dim(\ph\i(g)\cap X_{S_1})\ge\fra{i}{2}
-\fra{1}{2}(\dim S-\dim S_1)\}\le\dim\bY^a-i.$$
If this is violated for some $i\ge 0$, it would follow that the space of triples
$$\{(g,xP,x'P)\in\bY^a\T G^0/P\T G^0/P;x\i gx\in S_1U_P,x'{}\i gx'\in S_1U_P\}$$
has dimension $>\dim\bY^a-(\dim S-\dim S_1)$. Using 7.5(a) we see that this space of
triples has dimension $\le 2\nu-2\nu_L+\dim S_1-\dim A$. It follows that
$2\nu-2\nu_L+\dim S'-\dim A>\dim\bY^a-(\dim S-\dim S')$ hence
$\dim\bY^a<2\nu-2\nu_L+\dim S-\dim A$ contradicting (a). This proves that
$\ph^a_!K^a[d_0]$ is a perverse sheaf. It is semisimple by the decomposition theorem 
\cite{\BBD}. This proves (d).

We prove (e). We first show that ${}^{E_1}\ct_a\ne 0$. By the argument in the proof 
of Lemma 7.8(c) we see that it is enough to show that 
$H^{2\dim S-2\dim A}_c(S^a,\ce\ot\ce^*)\ne 0$. Since $\bbq$ is a direct summand of
$\ce\ot\ce^*$ it is enough to show that $H^{2\dim S-2\dim A}_c(S^a,\bbq)\ne 0$. This
follows from (a).

From Lemma 7.13 we see that the $\EE$-module structure defines an injective map
$\EE\ot{}^{E_1}\ct_a@>>>\tct_a$. Since ${}^{E_1}\ct_a\ne 0$, we have 
$(\EE\ot{}^{E_1}\ct_a)_\r\ne 0$ for any $\r\in\Irr\EE$, hence $(\tct_a)_\r\ne 0$. We 
have $\tct_a=H^{2d_0}_c(\bY^a,\ph_!K\ot\ph_!K^*)$ hence
$$\op_{\r\in\Irr\EE}\r\ot(\tct_a)_\r=
\op_{\r\in\Irr\EE}\r\ot H^{2d_0}_c(\bY^a,(\ph_!K)_\r\ot\ph_!K^*).$$
Hence for any $\r\in\Irr\EE$ we have
$(\tct_a)_\r=H^{2d_0}_c(\bY^a,(\ph_!K)_\r\ot\ph_!K^*)$. Since $(\tct_a)_\r\ne 0$, it
follows that $(\ph_!K)_\r|_{\bY^a}\ne 0$. The lemma is proved.

\head 8. Study of local systems on unipotent conjugacy classes\endhead  
\subhead 8.1\endsubhead
In this section we fix a connected component $D$ of $G$ such that $D$ contains some 
unipotent elements. Let $P$ be a parabolic of $G^0$ with Levi $L$ and let $S$ be an 
isolated stratum of $\tL=N_GL\cap N_GP$ contained in the connected component $\d$ of 
$\tL$ with $\d\sub D$. Assume that $S$ contains a unipotent $L$-conjugacy class 
$\boc_1$ (necessarily unique); then $S={}^\d\cz_L^0\boc_1$ and multiplication gives an
isomorphism ${}^\d\cz_L^0\T\boc_1@>\si>>S$. Let $\bY$ be the closure of $Y_{L,S}$ in 
$D$. Let $\ph:X@>>>\bY$ be as in 3.14. Let $A=A_{L,S}$. Then $\o\in A$ (see 7.1). We 
have $\bY^\o=\{g\in\bY; g \text{ unipotent}\}$, $X^\o=\ph\i(\bY^\o)\sub X$. Let 
$\ph^\o:X^\o@>>>\bY^\o$ be the restriction of $\ph:X@>>>\bY$. Let 
$d_0=2\nu-2\nu_L+\dim\boc_1$. We show: 

(a) {\it $X^\o$ and $\bY^\o$ are irreducible varieties of dimension $d_0$.}
\nl
The second projection makes $X^\o$ into a fibration over $G^0/P$ with all fibres 
isomorphic to $\bboc_1U_P$. To show that $X^\o$ is irreducible it is enough to notice
that $\bboc_1$ is irreducible. Since 
$\bY^\o=\ph^\o(X^\o)$, we see that $\bY^\o$ is irreducible. The statement about
dimension is a special case of Lemma 7.16(b). This proves (a).

Let $\cf_1$ be an irreducible $L$-equivariant local system on $\boc_1$ and let $\ce$ 
be the inverse image of $\cf_1$ under $S@>>>\boc_1,g\m g_u$. Then $\ce\in\cs(S)$. We 
assume that $(S,\ce)$ is a cuspidal pair of $\tL$. 

Let $X_S,\bce$ be as in 5.6. Let $K=IC(X,\bce)\in\cd(X)$, $K^*=IC(X,\bce^*)\in\cd(X)$.
Let $K^\o=K|_{X^\o},K^*{}^\o=K^*|_{X^\o}$. Define $\EE$ in terms of $L,S,\ce$ as in 
7.10.

Let $\cn_D$ be the set of all pairs $(\boc,\cf)$ where $\boc$ is a unipotent 
$G^0$-conjugacy class in $D$ and $\cf$ is an irreducible $G^0$-equivariant local 
system on $\boc$ (up to isomorphism). 

\proclaim{Proposition 8.2} (a) The restriction map
$\End_{\cd(\bY)}(\ph_!K)@>>>\End_{\cd(\bY^\o)}(\ph^\o_!K^\o)$ is an isomorphism.

(b) For any $\r\in\Irr\EE$, there is a unique $(\boc,\cf)\in\cn_D$ such that
$(\ph_!K)_\r|_{\bY^\o}[d_0]$ is $IC(\bboc,\cf)[\dim\boc]$ regarded as a simple
perverse sheaf on $\bY^\o$ (zero outside $\bboc$). Moreover, $\r\m(\boc,\cf)$ is 
an injective map $\g:\Irr\EE@>>>\cn_D$.

(c) The map $\g:\Irr\EE@>>>\cn_D$ in (b) depends only on $(L,S,\ce)$ and not on $P$.
\endproclaim
We prove (a). As in 7.14 we have a decomposition 
$\ph_!K=\op_{\r\in\Irr\EE}\r\ot(\ph_!K)_\r$ where $(\ph_!K)_\r[\dim\bY]$ are 
simple perverse sheaves on $\bY$. Restricting to $\bY^\o$ and using Lemma 7.16(c) we 
obtain a decomposition $\ph^\o_!K^\o=\op_{\r\in\Irr\EE}\r\ot(\ph_!K)_\r|_{\bY^\o}$. 
Hence the map in (a) factorizes as
$$\align&\End_{\cd(\bY)}(\ph_!K)=\op_{\r\in\Irr\EE}\End_{\cd(\bY)}(\r\ot(\ph_!K)_\r)\\&
@>b>>\op_{\r\in\Irr\EE}\End_{\cd(\bY^\o)}(\r\ot(\ph_!K)_\r|_{\bY^\o})
@>c>>\End_{\cd(\bY^\o)}(\ph^\o_!K^\o)\endalign$$   
where $c$ is injective. We have $b=\op_\r b_\r$ where
$$b_\r:\End(\r)\ot\End_{\cd(\bY)}((\ph_!K)_\r)@>>>
\End(\r)\ot\End_{\cd(\bY^\o)}((\ph_!K)_\r|_{\bY^\o})$$
is injective. (We use that 
$\End_{\cd(\bY)}((\ph_!K)_\r)=\bbq\sub\End_{\cd(\bY^\o)}((\ph_!K)_\r|_{\bY^\o})$,
since $(\ph_!K)_\r|_{\bY^\o}\ne 0$ by Lemma 7.16(e)). Hence $b$ is injective, so that 
the map in (a) is injective. It remains to show that 
$$\dim\End_{\cd(\bY^\o)}(\ph^\o_!K^\o)=\dim\End_{\cd(\bY)}(\ph_!K).$$
We have
$$\align&\dim\End_{\cd(\bY^\o)}(\ph^\o_!K^\o)
=\dim H^0_c(\bY^\o,\ph^\o_!K^\o[d_0]\ot\fD(\ph^\o_!K^\o[d_0]))\\&
=\dim H^0_c(\bY^\o,\ph^\o_!K^\o[d_0]\ot\ph^\o_!K^*{}^\o[d_0]))=
\dim H^{2d_0}_c(\bY^\o,\ph^\o_!K^\o\ot\ph^\o_!K^*{}^\o)\\&=
\dim H^{2d_0}_c(\bY^\o,\ph_!K\ot\ph_!K^*)=\dim\tct_\o
=\sum_{w\in\cw_S}\dim{}^{E_w}\ct_\o.\endalign$$
(The first equality follows by the argument in \cite{\CSII, 7.4} applied to the 
semisimple perverse sheaf $\ph^\o_!K^\o[d_0]$ on $\bY^\o$, see Lemma 7.16(d). The 
second equality holds since $\ph^\o$ is proper. The third equality is obvious. The 
fourth equality follows from Lemma 7.16(c). The fifth equality holds by definition. 
The sixth equality follows from Lemma 7.6 and Lemma 7.11.)

Let $w\in\cw_S$. As in the proof of Lemma 7.8 we have
$${}^{E_w}\ct_\o=H^{2\dim S-2\dim A}_c(\boc_1,\cf_1\ot\Ad(n_w\i)^*\cf_1^*)=
H^{2\dim\boc_1}_c(\boc_1,\cf_1\ot\Ad(n_w\i)^*\cf_1^*).$$
Since $\boc_1$ is irreducible and $\cf_1$ is an irreducible local system, this is $1$ 
dimensional if $\Ad(n_w\i)^*\cf_1\cong\cf_1$ (or equivalently, if
$\Ad(n_w\i)^*\ce\cong\ce$) and is $0$, otherwise. We see that 
$$\sum_{w\in\cw_S}\dim{}^{E_w}\ct_\o=|\cw_{\ce}|=\dim\EE$$
(see 7.10). It is then enough to show that $\dim\End_{\cd(\bY)}(\ph_!K)=\dim\EE$. 
This follows from the fact that $\ph_!K=IC(\bY,\p_!\tce)$ (see Proposition 5.7) and
$\End(IC(\bY,\p_!\tce))=\End(\p_!\tce)$, by the definition of an intersection
cohomology complex. This proves (a).

From the proof of (a) we see that both $b$ and $c$ are isomorphism. It follows that
the perverse sheaf $(\ph_!K)_\r|_{\bY^\o}[d_0]$ on $\bY^\o$ (see Lemma 7.16(d)) is 
simple and that, for $\r,\r'\in\Irr\EE$, we have
$(\ph_!K)_\r|_{\bY^\o}[d_0]\cong(\ph_!K)_{\r'}|_{\bY^\o}[d_0]$ if and only if $\r=\r'$.
Since the simple perverse sheaf $(\ph_!K)_\r|_{\bY^\o}[d_0]$ is $G^0$-equivariant and
$\bY^\o$ is a union of finitely many unipotent $G^0$-conjugacy classes, we see that 
$(\ph_!K)_\r|_{\bY^\o}[d_0]$ must be as in (b).

In the definition of $\g$ in (b), $P$ enters through the use of $\ph_!K$. However
$\ph_!K$ may be replaced by $IC(\bY,\p_!\tce)$ which does not depend on $P$. This 
proves (c). The proposition is proved.

\proclaim{Lemma 8.3} Let $\ph^S:X_S@>>>\bY$ be the restriction of $\ph:X@>>>\bY$. 
Let $(\boc,\cf)\in\cn_D$, $2d=d_0-\dim\boc$. The following four conditions are 
equivalent:

(i) $(\boc,\cf)=\g(\r)$ for some $\r\in\Irr\EE$ (notation of Proposition 8.2);

(ii) $\boc\sub\bY$ and $\cf$ is a constituent of the local system 
$(\ch^{2d}\ph^S_!\bce)|_\boc$;

(iii) $\boc\sub\bY$ and $\cf$ is a constituent of the local system 
$(\ch^{2d}\ph_!K)|_\boc$;

(iv) $\cf_1$ is a constituent of the local system $\ch^{\dim\boc-\dim\boc_1}f_!\cf$ 
on $\boc_1$ where $f:\boc\cap\boc_1U_P@>>>\boc_1$ is the restriction of 
$pr_1:\tL U_P@>>>\tL$.
\endproclaim
We first show:

(a) {\it the homomorphism $(\ch^{2d}\ph^S_!\bce)|_\boc@>>>(\ch^{2d}\ph_!K)|_\boc$ 
induced by the imbedding $X_S\sub X$ is an isomorphism}.
\nl
It is enough to show that for any $g\in\boc$, in the natural exact sequence
$$\align H^{2d-1}_c(\ph\i(g)\cap(X-X_S),K)&@>f>>H^{2d}_c(\ph\i(g)\cap X_S,\bce)@>>>
H^{2d}_c(\ph\i(g),K)\\&@>>>H^{2d}_c(\ph\i(g)\cap(X-X_S),K),\endalign$$
$$f:H^{2d-1}_c(\ph\i(g)\cap(X-X_S),K)@>>>H^{2d}_c(\ph\i(g)\cap X_S,\bce)\text{ is 
zero, and}\tag b$$
$$H^{2d}_c(\ph\i(g)\cap(X-X_S),K)=0.\tag c$$
The argument is similar to one in Lemma 7.6. For any stratum $S'$ of $\tL$ contained 
in $\bS$ let $X_{S'}$ be as in 5.6. We prove  (c). Using the partition 
$$\ph\i(g)\cap(X-X_S)=\cup_{S'\ne S}(\ph\i(g)\cap X_{S'})\tag d$$
we see that it is enough to show that $H^{2d}_c(\ph\i(g)\cap X_{S'},K)=0$ for any 
$S'\ne S$. Using the hypercohomology spectral sequence we see that it is enough to 
show that 

$H^i_c(\ph\i(g)\cap X_{S'},\ch^jK)\ne 0 \implies i+j<2d$ for $S'\ne S$.
\nl
If this last group is non-zero we would have

$i\le 2\dim(\ph\i(g)\cap X_{S'})\le 2d-(\dim S-\dim S')$
\nl
(see 4.2(b)) and $j<\dim X-\dim X_{S'}=\dim S-\dim S'$ (by the definition of $K$, 
since $S'\ne S$) hence $i+j<2d$, as desired. This proves (c).

To prove (b) we may assume that $\kk$ is an algebraic closure of a finite field 
$\FF_q$, that $G$ has a fixed $\FF_q$-structure with Frobenius map $F:G@>>>G$, that 
$P,L,S$ (hence $X$) are defined over $\FF_q$, that any stratum $S'$ as above is 
defined over $\FF_q$, that $F(g)=g$ and that we have an isomorphism 
$F^*\ce@>\si>>\ce$ which makes $\ce$ into a local system of pure weight $0$. Then we
have natural (Frobenius) endomorphisms of 
$H^{2d-1}_c(\ph\i(g)\cap(X-X_S),K),H^{2d}_c(\ph\i(g)\cap X_S,\bce)$ compatible with
$f$. To show that $f=0$, it is enough to show that

(e) $H^{2d}_c(\ph\i(g)\cap X_S,\bce)$ is pure of weight $2d$;

(f) $H^{2d-1}_c(\ph\i(g)\cap(X-X_S),K)$ is mixed of weight $\le 2d-1$.
\nl
Now (e) is clear since $\dim(\ph\i(g)\cap X_S)\le d$ (see 4.2(b)). We prove (f). 
Using the partition (d) we see that it is enough to prove:

$H^{2d-1}_c(\ph\i(g)\cap(X_{S'}),K)$ is mixed of weight $\le 2d-1$ for any 
$S'\ne S$.
\nl
Using the hypercohomology spectral sequence we see that it is enough to prove:

if $i,j$ are such that $2d-1=i+j$ then $H^i_c(\ph\i(g)\cap X_{S'},\ch^jK)$ is 
mixed of weight $\le 2d-1$ for any $S'$.
\nl
By Gabber's theorem \cite{\BBD, 5.3.2}, the local system $\ch^jK$ on $X_{S'}$ is 
mixed of weight $\le j$. Using Deligne's theorem \cite{\BBD, 5.1.14(i)}, we deduce 
that 

$H^i_c(\ph\i(g)\cap X_{S'},\ch^jK)$
\nl
is mixed of weight $\le i+j=2d-1$. This proves (a). In particular, conditions 
(ii),(iii) are equivalent.

Now $(\ch^{2d}\ph_!K)|_\boc=(\ch^{-\dim\boc}\ph_!K[d_0])|_\boc$ is isomorphic to

$\op_{\r'\in\Irr\EE}\r'\ot\ch^{-\dim\boc}(IC(\bboc',\cf')[\dim\boc'])|_\boc$
\nl
where $(\boc',\cf')=\g(\r')$. The sum may be restricted to those $\r'$ such that
$\boc\sub\bboc'$. If $\boc\ne\boc'$ then 
$\ch^{-\dim\boc}IC(\bboc',\cf')[\dim\boc']|_\boc=0$, by the definition of an 
intersection cohomology complex, while if $\boc=\boc'$, we have
$\ch^{-\dim\boc}IC(\bboc',\cf')[\dim\boc']|_\boc=\cf'$. We see that 
$$(\ch^{2d}\ph_!K)|_\boc\cong\op_{\r'\in\Irr\EE;\g(\r')=(\boc,\cf')}\r'\ot\cf'.\tag g
$$ 
Hence conditions (i) and (iii) are equivalent.

If (iv) holds then we have automatically $\boc\cap SU_P\ne\em$ hence 
$\boc\sub\bY$. The equivalence of (ii),(iv) is a special case of the equality (a) in
8.4 applied to the diagram $V@<f_2<<V'@>f_1>>\boc$ where 
$$V'=\{(g,xP)\in\boc\T G^0/P;x\i gx\in\boc_1U_P\},\qua V=P\bsl(\boc_1\T G^0)$$
with $P$ acting by $p:(y,x)\m(\p(p)y\p(p)\i,xp\i)$ where $\p$ is 
$pr_1:N_GP=\tL U_P@>>>\tL$ and

$f_2(g,xP)=\text{$P$-orbit of }(\p(x\i gx),x)$, $f_1(g,xP)=g$. 
\nl
The $G^0$-actions on $V$ by $g':(y,x)\m(y,g'x)$, on $V'$ by 
$g':(g,xP)\m(g'gg'{}\i,g'xP)$ and on $\boc$ by $g':g\m g'gg'{}\i$, are compatible 
with $f_1,f_2$ and are transitive on $V$ and $\boc$. We take 

$d_1=(\nu-\fra{1}{2}\dim\boc)-(\nu_L-\fra{1}{2}\dim\boc_1),
d_2=\fra{1}{2}(\dim\boc-\dim\boc_1)$. 
\nl
Then all fibres of $f_1$ have dimension $\le d_1$ and all fibres of $f_2$ have 
dimension $\le d_2$, see 4.2(a),(b). We set $N=d_1+\dim\boc=d_2+\dim V$. Let 
$\ce_1=\cf$. Let $\ce_2$ be the local system on $V$ whose inverse image under
$\boc_1\T G^0@>>>V$ is $\cf_1\ot\bbq$.

\subhead 8.4\endsubhead
Let $H$ be a connected algebraic group and let $H_1,H_2$ be two closed subgroups of
$H$. Let $V'$ be an algebraic variety with action of $H$. Let $f_1:V'@>>>H/H_1$,
$f_2:V'@>>>H/H_2$ be two $H$-equivariant maps. Let $e_1=\dim H/H_1,e_2=\dim H/H_2$.
Let $N$ be an integer such that all fibres of $f_1$ have dimension $\le N-e_1$ and 
all fibres of $f_2$ have dimension $\le N-e_2$.

Let $\ce_1$ be an $H$-equivariant irreducible local system on $H/H_1$. Let $\ce_2$ 
be an $H$-equivariant irreducible local system on $H/H_2$. Then 
$\ca_1=\ch^{2N-2e_1}f_1{}_!(f_2^*\ce_2)$ is an $H$-equivariant local system on 
$H/H_1$ and $\ca_2=\ch^{2N-2e_2}f_2{}_!(f_1^*\ce_1)$ is an $H$-equivariant local 
system on $H/H_2$. We have the following result.

(a) {\it The multiplicity of $\ce_1$ in $\ca_1$ is equal to the multiplicity of 
$\ce_2$ in $\ca_2$.}

These multiplicities are equal to 

$\dim H^{2e_1}_c(H/H_1,\ca_1\ot\ce_1^*)$, $\dim H^{2e_2}_c(H/H_2,\ca_2\ot\ce_2^*)$
\nl
respectively. It is enough to show

(b) $H^{2e_1}_c(H/H_1,\ca_1\ot\ce_1^*)\cong 
H^{2N}_c(V',f_1^*\ce_1^*\ot f_2^*\ce_2)$,

(c) $H^{2e_2}_c(H/H_2,\ca_2\ot\ce_2^*)\cong 
H^{2N}_c(V',f_1^*\ce_1\ot f_2^*\ce_2^*)$.
\nl
(Indeed, $\dim H^{2N}_c(V',f_1^*\ce_1^*\ot f_2^*\ce_2)=
\dim H^{2N}_c(V',f_1^*\ce_1\ot f_2^*\ce_2^*)$ since $\dim V'\le N$.) We have
$$\align&H^{2N}_c(V',f_1^*\ce_1^*\ot f_2^*\ce_2)=
H^{2N}_c(H/H_1,f_1{}_!(f_1^*\ce_1^*\ot f_2^*\ce_2))\\&=
H^{2N}_c(H/H_1,\ce_1^*\ot f_1{}_!f_2^*\ce_2)=
H^{2e_1}_c(H/H_1,\ce_1^*\ot\ch^{2N-2e_1}f_1{}_!f_2^*\ce_2)\endalign$$
where the last equality comes from a spectral sequence argument. This proves (b). 
The proof of (c) is entirely similar. This proves (a).

\subhead 8.5\endsubhead
In this subsection we preserve the setup of 7.4. 
Let $\ph'{}^{S'}:X'_{S'}@>>>\bY'$ be the restriction of $\ph':X'@>>>\bY'$. Let 
$\ph''{}^{S''}:X''_{S''}@>>>\bY''$ be the restriction of $\ph'':X''@>>>\bY''$. 
Assume that $S'={}^{\d'}\cz_{L'}^0\boc'_1,S''={}^{\d''}\cz_{L''}^0\boc''_1$ where 
$\boc'_1$ is a unipotent $L'$-conjugacy class and $\boc''_1$ is a unipotent 
$L''$-conjugacy class. Assume that $\ce'$ (resp. $\ce''$) is the inverse image under
$S'@>>>\boc'_1$ (resp. $S''@>>>\boc''_1$), $g\m g_u$, of an irreducible 
$L'$-equivariant (resp. $L''$-equivariant) local system on $\boc'_1$ (resp. 
$\boc''_1$). Assume that $(\boc,\cf)\in\cn_D$ satisfies

$\cf$ is a constituent of the local system $(\ch^{2d'}\ph'{}^{S'}_!\bce'{}^*)|_\boc$,

$\cf$ is a constituent of the local system $(\ch^{2d''}\ph''{}^{S''}_!\bce'')|_\boc$,
\nl 
where $d'=\nu-\nu_{L'}-\fra{1}{2}\dim\boc+\fra{1}{2}\dim\boc'_1$,
$d''=\nu-\nu_{L''}-\fra{1}{2}\dim\boc+\fra{1}{2}\dim\boc''_1$.

\proclaim{Lemma 8.6}In the setup of 8.5, the triples 
$(L',S',\ce'{}^*),(L'',S'',\ce'')$ are conjugate under an element of $G^0$.
\endproclaim
We have 
$(\ch^{2d'}\ph'{}^{S'}_!\bce'{}^*)|_\boc=((\ch^{2d'}\ph'{}^{S'}_!\bce')|_\boc)^*$ as 
local systems on $\boc$. (Indeed, for any $x\in\boc$ we have 
$H^{2d'}_c((\ph'{}^{S'})\i(x),\bce'{}^*)=H^{2d'}_c((\ph'{}^{S'})\i(x),\bce')^*$ 
since \lb $\dim(\ph'{}^{S'})\i(x)\le d'$.) Hence $\cf'{}^*$ is a direct summand of 
$(\ch^{2d'}\ph'{}^{S'}_!\bce')|_\boc$ and $\bbq$ is a direct summand of
$(\ch^{2d'}\ph'{}^{S'}_!\bce'\ot\ch^{2d''}\ph''{}^{S''}_!\bce'')|_\boc$. It follows 
that
$$H^{2\dim\boc}_c(\boc,
(\ch^{2d'}\ph'{}^{S'}_!\bce'\ot\ch^{2d''}\ph''{}^{S''}_!\bce'')|_\boc)\ne 0.\tag a$$
Let 
$$\align\fZ^\o_{S',S''}=&\{(g,x'P',x''P'')\in G\T G^0/P'\T G^0/P''; \\&g 
\text{ unipotent },x'{}\i gx'\in S'U_{P'},x''{}\i gx''\in S''U_{P''}\},\endalign$$
a special case of $\fZ^a_{S',S''}$ in 7.4. We have a partition
$$\fZ^\o_{S',S''}=\cup_{\boc'}\fZ^{\boc'}_{S',S''}\tag b$$
where $\boc'$ runs through the unipotent $G^0$-conjugacy classes in $D$ and

$\fZ^{\boc'}_{S',S''}=\{(g,x'P',x''P'')\in\fZ^\o_{S',S''};g\in\boc'\}$.
\nl
From (a) and the Leray spectral sequence of $pr_1:\fZ^\boc_{S',S''}@>>>\boc$ all of 
whose fibres have dimension $\le d'+d''$, we deduce
$$H^{2\ti d_0}_c(\fZ^\boc_{S',S''},\bce'\bxt\bce'')\ne 0$$
where
$\ti d_0=d'+d''+\dim\boc=2\nu-\nu_{L'}-\nu_{L''}+\fra{1}{2}(\dim\boc'_1+\dim\boc''_1)$.
Using this and the partition (b) where $\dim\fZ^\o_{S',S''}\le d_0$ (see 7.5(a)) we
deduce 
$$H^{2\ti d_0}_c(\fZ^\o_{S',S''},\bce'\bxt\bce'')\ne 0,$$
see 7.1(b). In other words, $\ct_\o$ (see 7.5) is $\ne 0$. By the argument in the 
proof of Lemma 7.9 we see that there exists a $G^0$-orbit $E$ on $G^0/P'\T G^0/P''$ 
such that ${}^E\ct_\o\ne 0$. Using Lemma 7.8 and its proof we see that there exists 
$n\in G^0$ such that $(P',nP'')\in E,n\i L'n=L'',n\i S'n=S''$ hence 
$n\i\boc'_1n=\boc''_1$ and 
${}^E\ct_\o=H^{2\dim\boc'_1}_c(\boc'_1,\ce'\ot\Ad(n\i)^*\ce'')$. Since this is 
$\ne 0$, the restriction of $\ce'\ot\Ad(n\i)^*\ce''$ to $\boc'_1$ contains $\bbq$ as
direct summand. Hence $\ce'{}^*|_{\boc'_1}\cong\Ad(n\i)^*\ce''|_{\boc'_1}$. Using 
our assumptions on $\ce',\ce''$ we deduce that $\ce'{}^*\cong\Ad(n\i)^*\ce''$. The 
lemma is proved.

\subhead 8.7\endsubhead
Let $\cn_D^0$ be the set of all pairs $(\boc,\cf)$ in $\cn_D$ such that 
$({}^D\cz_{G^0}^0\boc,\bbq\bxt\cf)$ is a cuspidal pair for $G$ (see 6.3); here we 
identify ${}^D\cz_{G^0}^0\T\boc={}^D\cz_{G^0}^0\boc$ using multiplication in $G$.

The following three conditions for $(\boc,\cf)\in\cn_D$ are equivalent:

{\it (i) $(\boc,\cf)\in\cn_D^0$;

(ii) for any proper parabolic $P$ of $G^0$ and any $g\in\boc\cap N_GP$ we have 

$H^{\dim\boc-\dim\dd}_c(\boc\cap gU_P,\cf)=0$
\nl
where $\dd$ is the $P/U_P$-conjugacy class of the image of $g$ in $N_GP/U_P$.

(iii) for any proper parabolic $P$ of $G^0$ and any unipotent $P/U_P$-conjugacy 
class $\dd$ in $N_GP/U_P$ we have $\ch^{\dim\boc-\dim\dd}h_!\cf=0$ where 
$h:\boc\cap\p\i(\dd)@>>>\dd$ is the restriction of the obvious map}
$\p:N_GP@>>>N_GP/U_P$.

Clearly, if (i) holds then (ii) holds. Assume now that (ii) holds. We show that (i)
holds. Let $g\in\boc\cap N_GP,z\in{}^D\cz_{G^0}^0$. Let $\boc_1$ be as in (ii). We 
must show that 
$H^{\dim\boc-\dim\boc_1}_c({}^D\cz_{G^0}^0\boc\cap zgU_P,\bbq\bxt\cf)=0$. This 
follows from (ii) since ${}^D\cz_{G^0}^0\boc\cap zgU_P=z(\boc\cap gU_P)$. (We use 
the fact that $gu$ is unipotent for any $u\in U_P$.) The equivalence of (ii),(iii) 
is clear.

\proclaim{Lemma 8.8} Let $(\boc,\cf)\in\cn_D$. There exists a quadruple 
$(P,L,\boc_1,\cf_1)$ such that (i)-(iii) below hold:

(i) $P$ is a parabolic of $G^0$ with Levi $L$, $\boc_1$ is a unipotent $L$-conjugacy
class in $\tL=N_GL\cap N_GP$ with $\boc_1\sub D$ and $\cf_1$ is an irreducible 
$L$-equivariant local system on $\cf_1$;

(ii) $(\boc_1,\cf_1)\in\cn_\d^0$ where $\d$ is the connected component of $\tL$ that
contains $\boc_1$;

(iii) $\cf_1$ is a constituent of the local system $\ch^{\dim\boc-\dim\boc_1}f_!\cf$ 
on $\boc_1$ where $f:\boc\cap\boc_1U_P@>>>\boc_1$ is the restriction of the 
projection $\tL U_P@>>>\tL$.
\endproclaim
We can clearly find a quadruple $(P,L,\boc_1,\cf_1)$ so that (i) and (iii) hold. 
For example, we can take $(P,L,\boc_1,\cf_1)=(G^0,G^0,\boc,\cf)$. We can further 
assume that $\dim P$ is minimum possible. We show that in this case, (ii) is 
automatically satisfied. Assume that $(\boc_1,\cf_1)\notin\cn_\d^0$. Then there 
exists a parabolic $P'$ of $G^0$ strictly contained in $P$ and a unipotent element 
$y\in N_GP'/U_{P'}$ such that 
$H^{\dim\boc_1-\dim\dd}_c(\boc_1\cap\p'{}\i(y),\cf_1)\ne 0$ where $\dd$ is the 
$P'/U_{P'}$-conjugacy class of $y$ in $N_GP'/U_{P'}$ and $\p':N_GP'@>>>N_GP'/U_{P'}$
is the obvious map. It follows that
$$H^{\dim\boc_1-\dim\dd}_c(\boc_1\cap\p'{}\i(y),\ch^{\dim\boc-\dim\boc_1}f_!\cf)
\ne 0,$$
($f$ as in (iii).) Using the Leray spectral sequence for the map
$f\i(\boc_1\cap\p'{}\i(y))=\boc\cap\boc_1U_P\cap\p'{}\i(y)@>>>\boc_1\cap\p'{}\i(y)$
(restriction of $f$) all of whose fibres have dimension 
$\le\fra{1}{2}(\dim\boc-\dim\boc_1)$ (see 4.2(a)) we deduce
$$H^{\dim\boc-\dim\dd}_c(\boc\cap\boc_1U_P\cap\p'{}\i(y),\cf)\ne 0.$$
We have a partition 
$\boc\cap\p'{}\i(y)=\cup_{\boc_2}(\boc\cap\boc_2U_P\cap\p'{}\i(y))$
where $\boc_2$ runs over the unipotent $L$-conjugacy class in $\d$. Since 
$\dim(\boc\cap\boc_2U_P\cap\p'{}\i(y))\le\fra{1}{2}(\dim\boc-\dim\dd)$ (see 4.2(a))
we deduce 
$$H^{\dim\boc-\dim\dd}_c(\boc\cap\p'{}\i(y)),\cf)\ne 0.$$
(We use repeatedly 7.1(b).) Hence there exists an irreducible constituent $\cf'_1$ 
of the local system $\ch^{\dim\boc-\dim\dd}f'_!\cf$ on $\dd$ where 
$f':\boc\cap\p'{}\i(\dd)@>>>\dd$ is the restriction of $\p'$. Using the minimality
of $\dim P$ we see that $\dim P\le\dim P'$. Since $P'\sub P$, we have $P'=P$. This 
contradiction proves that (ii) holds. The lemma is proved.

\subhead 8.9\endsubhead
Let $\tcm_D$ be the set of all quadruples $(P,L,\boc_1,\cf_1)$ that satisfy 
conditions (i),(ii) in 8.8. Let $\cm_D$ be the set of all triples $(L,\boc_1,\cf_1)$
such that $(P,L,\boc_1,\cf_1)$ for some $P$.

Now $G^0$ acts by conjugation on $\cm_D$; let $G^0\bsl\cm_D$ be the set of orbits. 
By associating to $(\boc,\cf)\in\cn_D$ the triple $(L,\boc_1,\cf_1)$ (part of the 
quadruple $(P,L,\boc_1,\cf_1)$ described in Lemma 8.8) we obtain a map 
$$\Ph_D:\cn_D@>>>G^0\bsl\cm_D;$$
the fact that $\Ph_D$ is well defined follows from the equivalence of (ii),(iv) in 
Lemma 8.3 and from Lemma 8.6. For any $(L,\boc_1,\cf_1)\in\cm_D$, we have a natural 
bijection
$$\Irr\EE\lra\Ph_D\i(L,\boc_1,\cf_1)$$
(the restriction of $\g$ in Proposition 8.2(b), where $S,\ce$ is related to 
$\boc_1,\cf_1$ as in 8.1).

\head 9. A restriction theorem\endhead
\subhead 9.1\endsubhead 
In this section we extend the results in \cite{\IC, \S8} to the disconnected case.  
Let $D$ be a connected component of $G$. Let $(P,L,\boc_1,\cf_1)\in\tcm_D$. Let $P'$
be a parabolic of $G^0$ such that $P\sub P'$ and $P'$ is normalized by some element
of $D$. Let $L'$ be the unique Levi of $P'$ such that $L\sub L'$. If $g\in\boc_1$ 
then $g$ normalizes $P$ and $gP'g\i$ is in the $G^0$-orbit of $P'$; also both 
$gP'g\i$ and $P'$ contain $P$ hence $gP'g\i=P'$. By the uniqueness of $L'$ we have 
also $gL'g\i=L'$. Thus, $\boc_1\sub\tL'$ where $\tL'=N_GL'\cap N_GP'$. Let $\d$ be
the unique connected component of $N_GL\cap N_GP$ that is contained in $D$. Then 
$\boc_1\sub\d$. Let $S={}^\d\cz_L^0\boc_1$. Define $\ce\in\cs(S)$ by 
$\ce=\bbq\bxt\cf_1$. Let $D'$ be the unique connected component of $\tL'$ that is 
contained in $D$. Define $\EE$ in terms of $L,S,\ce,G$ as in 7.10. Define $\EE'$ 
just like $\EE$, but in terms of $L,S,\ce,\tL'$. 

Let $Y=Y_{L,S}\sub D$ (relative to $G$). Let $\p:\tY=\tY_{L,S}@>>>Y$ (a principal 
$\cw_S$-bundle) be as in 3.13. Let $\cw'_S=\{n\in N_{L'}L;nSn\i=S\}/L$ (a subgroup 
of $\cw_S$) and let $\tY_1=\cw'_S\bsl\tY$. Then $\p$ is a composition 
$\tY@>\p'>>\tY_1@>\p''>>Y$ where $\p'$ is the obvious map and $\p''$ is induced by 
$\p$. Let $S^*$ be as in 3.11 and let ${}'S^*$ be the analogous set when $G$ is 
replaced by $\tL'$. Note that ${}'S^*$ is an open subset of $S$ containing $S^*$. 
We have a commutative diagram with cartesian squares
$$\CD
                                       \tY @<f<<G^0\T\tY'@>>>   \tY' @>>>  {}'\tY\\
                                      @V\p'VV      @V(1,q')VV    @Vq'VV   @V{}'qVV\\
                                         \tY_1@<f'<< G^0\T Y' @>>>     Y'  @>>>{}'Y
\endCD$$
where 

${}'\tY=\{(h,xL)\in\tL'\T L'/L;x\i hx\in{}'S^*\}$,

      ${}'Y=\cup_{x\in L'}x{}'S^*x\i$,

      ${}'q$ is $(h,xL)\m h$, the analogue of $\p:\tY@>>>Y$ when $G$ is replaced by
                 $\tL'$, (a principal $\cw'_S$-bundle),
      
     $\tY'=\{(h,xL)\in\tL'\T L'/L;x\i hx\in S^*\}$ (an open dense $\cw'_S$-stable 
     subset of ${}'\tY$),

     $Y'=\cw'_S\bsl\tY'$  (an open subset of ${}'Y$),

     $q'$ is the restriction of ${}'q$,

     $f(g,(h,xL))=(ghg\i,gxL)$ (a principal $L'$-bundle),

     $f'$ is defined by the commutativity of the left square (a principal 
                        $L'$-bundle),
\nl
and the unnamed maps are the obvious ones. 

Let $\tce$ be the local system on $\tY$ 
defined in 5.6; let $\tce'$ be the analogous local system on ${}'\tY$ (with $G$ 
replaced by $\tL'$). The inverse image of $\tce'$ under 
$G^0\T\tY'@>>>\tY'@>>>{}'\tY$ equals $f^*\tce$ hence the inverse image of 
${}'q_!\tce'$ under $G^0\T Y'@>>>Y'@>>>{}'Y$ equals $f'{}^*(\p'_!\tce)$. Since $f'$
is a principal $L$-bundle, $G^0\T Y'@>>>Y'$ is a principal $G^0$-bundle and 
$Y'@>>>{}'Y$ is an imbedding of an open dense subset, we have canonically
$$\End({}'q_!\tce')=\End(f'{}^*(\p'_!\tce))=\End(\p'_!\tce).$$
Since $\EE'=\End({}'q_!\tce')$, we see that 
$$\EE'=\End(\p'_!\tce).\tag a$$
It follows that
$$\p'_!\tce=\op_{\r'\in\Irr\EE'}\r'\ot(\p'_!\tce)_{\r'}$$
where $(\p'_!\tce)_{\r'}$ are irreducible local systems on $\tY_1$. (See 7.14.) The
obvious algebra homomorphism $\End(\p'_!\tce)@>>>\End(\p''_!(\p_!\ce))=\End(\p_!\ce)$,
that is, $\EE'@>>>\EE$, maps the summand $\EE'_w$ ($w\in\cw'_S)$ of $\EE'$ 
isomorphically onto the summand $\EE_w$ ($w\in\cw'_S\sub\cw_S$) of $\EE$; hence it 
is injective. We use this to identify $\EE'$ with a subalgebra of $\EE$.

Recall from 7.14 that $\p_!\tce=\op_{\r\in\Irr\EE}\r\ot(\p_!\tce)_\r$ where 
$(\p_!\tce)_\r$ is a local system on $Y$. We have
$$\p_!\tce=\p''_!\p'_!\tce=\op_{\r'\in\Irr\EE'}\r'\ot\p''_!(\p'_!\tce)_{\r'}$$
hence 
$$\p''_!(\p'_!\tce)_{\r'})=\Hom_{\EE'}(\r',\p_!\tce)=
\Hom_{\EE'}(\r',\op_{\r\in\Irr\EE}\r\ot(\p_!\tce)_\r).$$
We see that
$$\p''_!(\p'_!\tce)_{\r'})=\op_{\r\in\Irr\EE}\Hom_{\EE'}(\r',\r)\ot(\p_!\tce)_\r.
\tag b$$

\subhead 9.2\endsubhead
Let $\bY$ be the closure of $Y$ in $D$. Let 
${}'\bY$ be the closure of ${}'Y$ in $D'$. Let $\ph:X@>>>\bY$ be as in 3.14. Let 
$X_1=\{(g,xP')\in G\T G^0/P';x\i gx\in{}'\bY U_{P'}\}$. Then $\ph$ factors as
$X@>\ph'>>X_1@>\ph''>>\bY$ where $\ph'(g,xP)=(g,xP'),\ph''(g,xP')=g$. Clearly,
$\ph',\ph''$ are proper, surjective. We have a commutative diagram
$$\CD
\tY@>\p'>>\tY_1@>\p''>>Y\\
@Vj_0VV   @Vj_1VV    @Vj_2VV\\
X@>\ph'>>X_1@>\ph''>>\bY
\endCD$$
where $j_2$ is the obvious imbedding (as an open set), $j_0(g,xL)=(g,xP)$ (an 
isomorphism of $\tY$ with the open subset $\ph\i(Y)$ of $X$, see Lemma 5.5) and 
$j_1\p'(g,xP)=(g,xP')$. Now $j_1$ is an isomorphism onto the open subset
$\ph''{}\i(Y)$ of $X_1$. (We show only that $j_1$ is a bijection 
$\tY_1@>\si>>\ph''{}\i(Y)$. Since $\ph'$ is surjective and 
$Y=\ph(\tY),\tY=\ph\i(Y)$, 
we have $j_1(\tY_1)=\ph'(\tY)=\ph''{}\i(Y)$. To show that $j_1$ is injective, it is 
enough to show that, if two points $(g,xL),(g',x'L)$ of $\tY$ have the same image 
under $\ph'$ then they are in the same $\cw'_S$-orbit. We have $g=g',x'=xp'$ where 
$p'\in P'$. From Lemma 3.13 we see that $x'=xn\i$ for some $n\in N_{G^0}L,nSn\i=S$.
We have $n\i=p'$. Since $n$ normalizes $L$ and $P'$, it also normalizes $L'$ (by the
uniqueness of $L'$). Since $N_{G^0}L'\cap P'=L'$, we have $n\in L'$. Thus 
$(g,xL),(g',x'L)$ are in the same $\cw'_S$-orbit, as required.) Note also that 
$\tY_1$ is smooth (since $\tY$ is smooth). We identify $\tY,\tY_1,Y$ with open 
subsets of $X,X_1,\bY$ via $j_0,j_1,j_2$. We have a commutative diagram with 
cartesian squares
$$\CD
X   @<p_1<<X''       @>p_2>>X'   \\
@V\ph'VV @V\psi VV           @V{}'\ph VV    \\
X_1@<p_3<<  G^0\T U_{P'}\T{}'\bY@>p_4>>{}'\bY
\endCD$$
where

$X'=\{(h,x(P\cap L'))\in\tL'\T L'/(P\cap L');x\i hx\in\bS(U_P\cap L')\}$,

$X''=\{(g_1,q,xP)\in G^0\T N_GP'\T P'/P;x\i qx\in\bS U_P\}$,

$p_1$ is $(g_1,q,xP)\m(g_1qg_1\i,g_1xP)$, a principal $P'$-bundle,

$p_2$ is $(g_1,l'u',l'_1P)\m(l',l'_1(P\cap L'))$ with
$l'\in\tL',l'_1\in L',u'\in U_{P'}$, a principal $G^0\T U_{P'}$-bundle,

$p_3$ is $(g_1,u,h)\m(g_1hug_1\i,g_1P')$, a principal $P'$-bundle,

$p_4$ is $(g_1,u,h)\m h$, a principal $G^0\T U_{P'}$-bundle,

${}'\ph$ is $(h,x(P\cap L'))\m h$,

$\psi$ is $(g_1,l'u',l'_1P)@>>>(g_1,u',l')$, with
$l'\in\tL',l'_1\in L',u'\in U_{P'}$.
\nl
Since $p_3,p_4$ are principal bundles with connected group we have
$p_3^*IC(X_1,p'_!\tce)=p_4^*IC({}'\bY,{}'q_!\tce')$. (Both can be identified with 
$IC(G^0\T U_{P'}\T{}'bY,\bbq\bxt\bbq\bxt{}'q_!\tce')$.)

Let $K,K^*\in\cd(X)$ be as in 5.7. Let $K'\in\cd(X')$ be the analogous object (with
$G$ replaced by $\tL'$). From the definitions we have $p_1^*K=p_2^*K'$. From the 
commutative diagram above it then follows that 
$p_3^*\ph'_!K=p_4^*({}'\ph_!K')=p_4^*IC({}'\bY,{}'q_!\tce')$ (the last equality comes 
from Proposition 5.7 for $\tL'$ instead of $G$) hence 
$p_3^*\ph'_!K=p_3^*IC(X_1,p'_!\tce)$. Since $p_3$ is a principal $P'$-bundle we see 
that 
$$\ph'_!K=IC(X_1,\p'_!\tce).$$
From this and 9.1(a) we see that $\End(\ph'_!K)=\EE'$ and 
$\ph'_!K=\op_{\r'\in\Irr\EE'}\r'\ot(\ph'_!K)_{\r'}$ where 
$$(\ph'_!K)_{\r'}=IC(X_1,(\p'_!\tce)_{\r'}).\tag a$$
We now show
$$\ph''_!(\ph'_!K)_{\r'}=IC(\bY,\p''_!(\p'_!\tce)_{\r'})\tag b$$
for any $\r'\in\Irr\EE'$. From (a) we see that $\ph''_!(\ph'_!K)_{\r'}|_Y$ is the 
local system $\p''_!(\p'_!\tce)_{\r'}$. Since $\ph''$ is proper, (b) is a 
consequence of (a), of the assertion (c) below and the analogous assertion with $K$
replaced by $K^*$:

(c) For any $i>0$ we have $\dim\supp\ch^i(\ph''_!(\ph'_!K)_{\r'})<\dim\bY -i$.
\nl
This is checked as follows. We have 
$$\supp\ch^i(\ph''_!(\ph'_!K)_{\r'})\sub\supp\ch^i(\ph''_!(\ph'_!K))=
\supp\ch^i(\ph_!K)$$
hence (c) is a consequence of a statement in the proof of Proposition 5.7. Thus, (b)
is verified. Combining (b) with 9.1(b), we see that for any $\r'\in\Irr\EE'$ we
have
$$\ph''_!(\ph'_!K)_{\r'}=\op_{\r\in\Irr\EE}\Hom_{\EE'}(\r',\r)\ot(\ph_!K)_\r.\tag d
$$
(Recall from 7.14 that $\ph_!K=\op_{\r\in\Irr\EE}\r\ot(\ph_!K)_\r$ where 
$(\ph_!K)_\r=IC(\bY,(\p_!\tce)_\r)$.)

\subhead 9.3\endsubhead
Let $(\boc,\cf)\in\cn_D,(\boc',\cf')\in\cn_{D'}$. (In particular, $\boc$ is a 
unipotent $G^0$-conjugacy class in $D$ and $\boc'$ is a unipotent $L'$-conjugacy 
class in $D'$.) Let 
$$d'_0=2\nu_{L'}-2\nu_L+\dim\boc_1,d=\nu-\nu_L-\fra{1}{2}(\dim\boc-\dim\boc_1).$$ 
Let $\Ph_D:\cn_D@>>>G^0\bsl\cm_D$ be as in 8.9. Let 
$\Ph_{D'}:\cn_{D'}@>>>L'\bsl\cm_{D'}$ be the analogous map (with $G$ replaced by 
$\tL'$). We assume that $\Ph_D(\boc,\cf)$ is the $G^0$-orbit of $(L,\boc_1,\cf_1)$ 
and that $(\boc,\cf)$ corresponds to $\r\in\Irr\EE$ (see 8.9). We assume that 
$\Ph_{D'}(\boc',\cf')$ is the $L'$-orbit of $(L,\boc_1,\cf_1)$ and that 
$(\boc',\cf')$ corresponds to $\r'\in\Irr\EE'$ (see 8.9). Let 
$X_1^\o=\{(g,xP')\in X^1; g \text{ unipotent}\}$. Let 
$$R=\{(g,xP')\in G\T G^0/P';x\i gx\in\bboc'U_{P'}\}\sub X_1^\o.$$
We show:
$$\supp(\ph'_!K)_{\r'}\cap X_1^\o\sub R.\tag a$$
Let $(g,xP')\in\supp(\ph'_!K)_{\r'}\cap X_1^\o$. The isomorphism 
$p_3^*\ph'_!K=p_4^*({}'\ph_!K')$ in 9.2 is compatible with the action of $\EE'$. 
Hence $p_3^*(\ph'_!K)_{\r'}=p_4^*({}'\ph_!K')_{\r'}$ and

$p_3\i(\supp(\ph'_!K)_{\r'})=p_4\i(\supp({}'\ph_!K')_{\r'})$.
\nl
Hence there exists
$(g_1,u,h)\in X''$ such that $(g,xP')=(g_1hug_1\i,g_1P')$ and 
$h\in\supp({}'\ph_!K')_{\r'})$. Since $g$ is unipotent, $hu$ is unipotent hence $h$
is unipotent. Now a unipotent element in $\supp({}'\ph_!K')_{\r'})$ must be in 
$\bboc'$ since, by Proposition 8.2 (for $\tL'$ instead of $G$),

(b) $({}'\ph_!K')_{\r'}$ {\it restricted to the unipotent set in $D'$ is 
$IC(\bboc',\cf')[\dim\boc'-d'_0]$ (extended by zero outside $\bboc'$).}
\nl
Thus, $h\in\bboc'$. We have $x=g_1p'$ for some $p'\in P'$ and $g=g_1hug_1\i$. 
Hence $x\i gx=p'hup'{}\i\in\bboc'U_{P'}$ and $(g,xP')\in R$. This proves (a).

We have a partition $R=\cup_{\boc''}R_{\boc''}$ where $\boc''$ runs over the 
unipotent $L'$-conjugacy classes in $\bboc'$ and
$R_{\boc''}=\{(g,xP')\in G\T G^0/P';x\i gx\in\boc''U_{P'}\}$. Then 
$R':=R_{\boc'}$ is open in $R$. Clearly, 
$p_1\i(R)=p_2\i(\bboc')=G^0\T U_{P'}\T\bboc'$ and
$p_1\i(R_{\boc''})=p_2\i(\boc'')=G^0\T U_{P'}\T\boc''$. 

Let $\tcf'$ be the local system on $R'$ whose inverse image under 
$p_1:G^0\T U_{P'}\T\boc'@>>>R'$ equals the inverse image of $\cf'$ under 
$p_2:G^0\T U_{P'}\T\boc'@>>>\boc'$.

Since $p_1,p_2$ are principal bundles with connected group it follows that the 
inverse image of $IC(R,\tcf')$ under $p_1:G^0\T U_{P'}\T\bboc'@>>>R$ equals the
inverse image of $IC(\bboc',\cf')$ under $p_2:G^0\T U_{P'}\T\bboc'@>>>\bboc'$. 
It follows that

(c) $(\ph'_!K)_{\r'}|_{X_1^\o}=IC(R,\tcf')[\dim\boc'-d'_0]$ (extended by zero 
outside $R$).
\nl
(Using $p_1^*$ this is reduced to (b).)

For any subvariety $T$ of $X_1$ we denote by ${}_T\ph'':T@>>>\bY$ the restriction of
$\ph'':X^1@>>>\bY$ to $T$.

\proclaim{Proposition 9.4}In the setup of 9.3, let 
$d'=\fra{1}{2}(\dim\boc-\dim\boc'),n=\nu-\nu_{L'}-d'$. The following five numbers 
coincide:

(i) $\dim\Hom_{\EE'}(\r',\r)$;

(ii) the multiplicity of $\cf$ in the local system 
$\cl_1=\ch^{2d}(\ph''_!(\ph'_!K)_{\r'})|_\boc$;

(iii) the multiplicity of $\cf$ in the local system
$\cl_2=\ch^{2n}({}_R\ph''_!IC(R,\tcf'))|_\boc$;

(iv) the multiplicity of $\cf$ in the local system

$\cl_3=\ch^{2n}({}_{R'}\ph''_!IC(R,\tcf'))|_\boc=\ch^{2n}({}_{R'}\ph''_!\tcf')|_\boc;$

(v) the multiplicity of $\cf'$ in the local system $\ch^{2d'}(f_!\cf)$ where 
$f:\boc'U_{P'}\cap\boc@>>>\boc'$ is the restriction of $pr_1:\tL'U_{P'}@>>>\tL'$. 
\endproclaim
The proof will be given in 9.5-9.7.

\subhead 9.5\endsubhead
From 8.3(g) we see that, for $\ti\r\in\Irr\EE$, the multiplicity of $\cf$ in the 
local system $(\ch^{2d}(\ph_!K)_{\ti\r})|_\boc$ is $1$ if $\ti\r=\r$ and is $0$ 
otherwise. Hence from 9.2(d) it follows that the numbers (i),(ii) in 9.4 are equal.

We show that $\cl_1=\cl_2$. By 9.3(c) we may replace $IC(R,\tcf')$ in $\cl_2$ by 

$(\ph'_!K)_{\r'}|_{X_1^\o}[-\dim\boc'+d'_0]$
\nl
so that $\cl_2=\ch^{2d}({}_R\ph''_!((\ph'_!K)_{\r'}|_R)|_\boc$. (We have 
$2n=2d+\dim\boc'-d'_0$.) It is enough to show that
${}_{(X_1-R)}\ph''_!((\ph'_!K)_{\r'}|_{X_1-R})|_\boc=0$. Assume this is not true. Then
there exist $(g,xP')\in X_1-R$ such that $g\in\boc$ and
$(g,xP')\in\supp(\ph'_!K)_{\r'}$. Since $g$ is unipotent, this contradicts 9.3(a). 
We see that the numbers (ii),(iii) in 9.4 are equal.

\subhead 9.6\endsubhead
We show that $\cl_2=\cl_3$. For any $g\in\boc$ we consider the natural exact 
sequence
$$\align&H^{2d-1}_c(\ph''{}\i(g)\cap(R-R'),(\ph'_!K)_{\r'})@>\xi>>
H^{2d}_c(\ph''{}\i(g)\cap R',(\ph'_!K)_{\r'})@>>>\\&
H^{2d}_c(\ph''{}\i(g)\cap R,(\ph'_!K)_{\r'})@>>>
H^{2d}_c(\ph''{}\i(g)\cap(R-R'),(\ph'_!K)_{\r'}).\endalign$$
It is enough to show that the the middle map is an isomorphism. It is enough to 
show that $H^{2d}_c(\ph''{}\i(g)\cap(R-R'),(\ph'_!K)_{\r'})=0$ and that $\xi=0$. By
9.3(c) we may replace $(\ph'_!K)_{\r'}|_{X_1^\o}$ by $IC(R,\tcf')[\dim\boc'-d'_0]$.
We see that it is enough to show that

(a) $H^{2n}_c(\ph''{}\i(g)\cap(R-R'),IC(R,\tcf'))=0$ and that  

(b) $H^{2n-1}_c(\ph''{}\i(g)\cap(R-R'),IC(R,\tcf'))@>\xi>>
H^{2n}_c(\ph''{}\i(g)\cap R',IC(R,\tcf'))$ is zero.
\nl
From Proposition 4.2(b) we have

(c) $\dim(\ph''{}\i(g)\cap R_{\boc''})\le
\nu-\nu_{L'}-\fra{1}{2}(\dim\boc-\dim\boc'')$
\nl
for any $L'$-conjugacy class $\boc''$ in $\bboc'$. 

If the cohomology group in (a) is non-zero, then, using the partition

(d) $\ph''{}\i(g)\cap(R-R')=\cup_{\boc''\ne\boc'}(\ph''{}\i(g)\cap R_{\boc''})$ 
\nl
we see that $H^{2n}_c(\ph''{}\i(g)\cap R_{\boc''},IC(R,\tcf'))\ne 0$ for some
$\boc''\ne\boc'$. Hence there exist $i,j$ such that $2n=i+j$ and
$H^i_c(\ph''{}\i(g)\cap R_{\boc''},\ch^j(IC(R,\tcf')))\ne 0$. It follows that 
$i\le 2\dim R_{\boc''}$. The local system $\ch^j(IC(R,\tcf'))|_{R_{\boc''}}$ is 
$\ne 0$ so that $R_{\boc''}\sub\supp\ch^j(IC(R,\tcf'))$ and 
$\dim R_{\boc''}\le\dim R-j$. It follows that 
$i+j\le 2\dim R_{\boc''}+\dim R-\dim R_{\boc''}=\dim R+\dim R_{\boc''}<2n$ (we use
(c)) in contradiction with $i+j=2n$. This proves (a).

To prove (b) we may assume that $\kk$ is an algebraic closure of a finite field 
$\FF_q$, that $G$ has a fixed $\FF_q$-structure with Frobenius map $F:G@>>>G$, that 
$P,P',L,L',S$ (hence $X_1,\ph''$) are defined over $\FF_q$, that any $\boc''$ as 
above is defined over $\FF_q$, that $F(g)=g$ and that we have an isomorphism 
$F^*\cf'@>\si>>\cf'$ which makes $\cf'$ into a local system of pure weight $0$. Then
we have natural (Frobenius) endomorphisms of 
$$\align&H^{2n-1}_c(\ph''{}\i(g)\cap(R-R'),IC(R,\tcf')),\\&H^{2n}_c
(\ph''{}\i(g)\cap R',IC(R,\tcf'))=H^{2n}_c(\ph''{}\i(g)\cap R',\tcf')\endalign$$
compatible with $\xi$. To show that $\xi=0$, it is enough to show that

(e) $H^{2n}_c(\ph''{}\i(g)\cap R',\tcf')$ is pure of weight $2n$;

(f) $H^{2n-1}_c(\ph''{}\i(g)\cap(R-R'),IC(R,\tcf'))$ is mixed of weight $\le 2n-1$.
\nl
Now (e) is clear since $\dim(\ph''{}\i(g)\cap R')\le n$ (see (c)). We prove (f). 
Using the partition (d) we see that it is enough to prove:

$H^{2n-1}_c(\ph''{}\i(g)\cap R_{\boc''},IC(R,\tcf'))$ is mixed of weight $\le 2n-1$
for any $\boc''\ne\boc'$.
\nl
Using the hypercohomology spectral sequence we see that it is enough to prove:

if $i,j$ are such that $2n-1=i+j$ then 
$H^i_c(\ph''{}\i(g)\cap R_{\boc''},\ch^j(IC(R,\tcf')))$ is mixed of weight 
$\le 2n-1$ for any $\boc''$.
\nl
By Gabber's theorem \cite{\BBD, 5.3.2}, the local system $\ch^j(IC(R,\tcf'))$ is 
mixed of weight $\le j$. Using Deligne's theorem \cite{\BBD, 5.1.14(i)}, we deduce 
that 

$H^i_c(\ph''{}\i(g)\cap R_{\boc''},\ch^j(IC(R,\tcf')))$
\nl
is mixed of weight $\le i+j=2n-1$. This proves (b). We have shown that $\cl_2=\cl_3$. 
It follows that the numbers (iii),(iv) in 9.4 are equal.

\subhead 9.7\endsubhead
Consider the diagram $V@<f_2<<V'@>f_1>>\boc$ where 

$V'=\{(g,xP')\in\boc\T G^0/P';x\i gx\in\boc'U_{P'}\}=\ph''{}\i(\boc)\cap R'$,

$V=P'\bsl(\boc'\T G^0)$
\nl
with $P'$ acting by $p':(y,x)\m(\p(p')y\p(p')\i,xp'{}\i)$, $\p$ is 
$pr_1:\tL'U_{P'}@>>>\tL'$, $f_2(g,xP')=\text{$P'$-orbit of }(\p(x\i gx),x)$, 
$f_1(g,xP')=g$.

The $G^0$-actions on $V$ by $g':(y,x)\m(y,g'x)$, on $V'$ by 
$g':(g,xP')\m(g'gg'{}\i,g'xP')$ and on $\boc$ by $g':g\m g'gg'{}\i$, are compatible 
with $f_1,f_2$ and are transitive on $V$ and $\boc$. 

Then all fibres of $f_1$ have dimension $\le n$ and all fibres of $f_2$ have 
dimension $\le d'$, see 4.2(a),(b). We set $N=n+\dim\boc=d'+\dim V$. We apply 8.4(a)
with $\ce_1=\cf$ and with $\ce_2$ the local system on $V$ whose inverse image under 
$\boc'\T G^0@>>>V$ is $\cf'\ot\bbq$. We see that the numbers (iv),(v) in 9.4 are 
equal. This completes the proof of Proposition 9.4.

\subhead 9.8\endsubhead
In the setup of 9.3, we now drop the assumption that 

(a) $\Ph_D(\boc,\cf)$ is the $G^0$-orbit of $(L,\boc_1,\cf_1)$
\nl
and replace it by the assumption that the multiplicity in 9.4(v) is non-zero. We 
show that in this case (a) automatically holds.

Our assumption implies that $f_!\cf$ in 9.4(v) is non-zero. In particular,
$\boc'U_{P'}\cap\boc\ne\em$ hence ${}'\bY U_{P'}\cap\boc\ne\em$. We have 
${}'\bY U_{P'}\sub\bY$ hence $\bY\cap\boc\ne\em$ and $\boc\sub\bY$. Then the arguments
in 9.5-9.7 still show that the multiplicities in 9.4(ii),(iii),(iv),(v) are equal 
and in particular, they are all non-zero. It follows that the multiplicity of $\cf$ in
the local system $\ch^{2d}(\ph''_!(\ph'_!K))|_\boc=\ch^{2d}(\ph_!K)|_\boc$ is non-zero.
We now use the equivalence of (i) and (iii) in Lemma 8.3; our assertion follows.

\head 10. Preparatory results\endhead
The following result is a generalization of known results from the connected case to 
the disconnected case. Thus, (a)-(c) generalizes \cite{\IC, 2.9}, (d) generalizes
\cite{\HS, Prop.3} and (e) generalizes \cite{\IC, 9.3}.

\proclaim{Lemma 10.1} Let $P$ be a parabolic of $G^0$ with Levi $L$. Let 
$g\in\tL=N_GL\cap N_GP$. Let $\la g\ra$ be the $G^0$-conjugacy class of $g$. Let 
$\la g\ra_L$ be the $L$-conjugacy class of $g$. Let $F=\{vgv\i;v\in U_P\}$. Let 
$V=Z_G(g)^0/(Z_G(g)^0\cap P)$. Then:

(a) $F$ is an irreducible component of $\la g\ra\cap gU_P$ of dimension
$\fra{1}{2}(\dim\la g\ra-\dim\la g\ra_L)$;

(b) $\dim V=(\nu-\fra{1}{2}\dim\la g\ra)-(\nu_L-\fra{1}{2}\dim\la g\ra_L)$;

(c) $Z_G(g)^0\cap P=Z_P(g)^0$;

(d) $Z_{U_P}(g)$ is connected;

(e) if $P\ne G^0$ and $g$ is unipotent then
$(\nu-\fra{1}{2}\dim\la g\ra)-(\nu_L-\fra{1}{2}\dim\la g\ra_L)=\dim Z_{U_P}(g)>0$.
\endproclaim
We prove (a). From the semidirect product decomposition $P=LU_P$ we obtain a 
semidirect product decomposition $Z_P(g)=Z_L(g)Z_{U_P}(g)$. Hence

$\dim Z_P(g)=\dim Z_L(g)+\dim Z_{U_P}(g)$. 
\nl
Exactly the same argument shows that, if $P'$ is the unique parabolic of $G^0$ with 
Levi $L$ such that $P\cap P'=L$, then $\dim Z_{P'}(g)=\dim Z_L(g)+\dim Z_{U_{P'}}(g)$. 
(By the uniqueness of $P'$ we have $g\in N_GP'$.) Consider the map 
$$Z_P(g)^0\T Z_{P'}(g)^0@>>>Z_G(g)^0,(g_1,g_2)\m g_1g_2.$$
The pairs $(g_1,g_2),(g'_1,g'_2)$ are mapped to the same element in
$Z_G(g)^0$ if and only if $g'_1=g_1g_0,g'_2=g_0\i g'_1$ where 
$g_0\in Z_P(g)^0\cap Z_{P'}(g)^0$. Note also that 
$$Z_L(g)^0\sub Z_P(g)^0\cap Z_{P'}(g)^0\sub Z_L(g).$$
It follows that 
$$\dim Z_G(g)^0=\dim Z_P(g)^0+\dim Z_{P'}(g)^0-\dim(Z_P(g)^0\cap Z_{P'}(g)^0)+\d,
\qua (\d\ge 0)$$
hence
$$\align&\dim Z_G^0(g)=\dim Z_P(g)^0+\dim Z_{P'}(g)^0-\dim Z_L(g)^0+\d\\&=
\dim Z_{U_P}(g)+\dim Z_{U_{P'}}(g)+\dim Z_L(g)^0+\d.\tag f\endalign$$
Introducing here $\dim Z_{U_P}(g)=\dim U_P-\dim F$ and the analogous equality
\lb $\dim Z_{U_{P'}}(g)=\dim U_{P'}-\dim F'$ where $F'=\{v'gv'{}\i;v'\in U_{P'}\}$, 
we obtain
$$\dim Z_G^0(g)=\dim U_P+\dim U_{P'}+\dim Z_L(g)^0-\dim F-\dim F'+\d$$
hence
$$\dim F+\dim F'=\dim\la g\ra-\dim\la g\ra_L+\d.\tag g$$
Now $F$ is contained in $\la g\ra\cap gU_P$ and is closed in $gU_P$ (it is an orbit
of a unipotent group action on an affine variety). By 4.2(a), any irreducible 
component of $\la g\ra\cap gU_P$ has dimension 
$\le\fra{1}{2}(\dim\la g\ra-\dim\la g\ra_L)$. Hence
$\dim F\le\fra{1}{2}(\dim\la g\ra-\dim\la g\ra_L)$ and similarly
$\dim F'\le\fra{1}{2}(\dim\la g\ra-\dim\la g\ra_L)$. Comparing with (g) we see that
$\d=0$ and that $F$ (resp. $F'$) is an irreducible component of $\la g\ra\cap gU_P$ 
(resp. of $\la g\ra\cap gU_{P'}$) of dimension 
$\fra{1}{2}(\dim\la g\ra-\dim\la g\ra_L)$.

We prove (b). From the proof of (a) we see that
$\dim Z_G(g)=\dim Z_P(g)+\dim Z_{U_{P'}}(g)$ hence $\dim V=\dim Z_{U_{P'}}(g)$. On 
the other hand we have
$$\align&(\nu-\fra{1}{2}\dim\la g\ra)-(\nu_L-\fra{1}{2}\dim\la g\ra_L)\\&=
\fra{1}{2}(\dim Z_G(g)-\dim Z_L(g))=\fra{1}{2}(\dim Z_{U_P}(g)+\dim Z_{U_{P'}}(g)).
\tag g\endalign$$
(We have used (f) and the equality $\d=0$.) Hence to prove (b) it is enough to prove
the equality 

(h) $\dim Z_{U_P}(g)=\dim Z_{U_{P'}}(g)$.
\nl
This follows from the equality $\dim F=\dim F'$ in the proof of (a) and from the 
equality $\dim U_P=\dim U_{P'}$.

We prove (c). Let $H=Z_G(g)$ and let $T=(\cz_L\cap H)^0$. Then $T$ is a torus 
contained in $H^0$ hence $Z_{H^0}(T)$ is connected. By 1.10(a) we have 
$L=Z_{G^0}(T)$. It follows that $L\cap H^0=Z_{H^0}(T)$ is connected hence 
$L\cap H^0\sub Z_L(g)^0$. The group
$Z_G(g)^0\cap P$ contains $Z_P(g)^0$ as a subgroup of finite index. Let $g_0$ be a 
fixed element of $Z_G(g)^0\cap P$. From the proof of (a) we have $\d=0$ hence the set 
$A$ of all products $g_1g_2$ with $g_1\in Z_P(g)^0,g_2\in Z_{P'}(g)^0$ is constructible
dense in $Z_G(g)^0$. Hence $g_0A$ is again constructible dense in $Z_G(g)^0$ so that it
must meet $A$. Thus $g_0g_1g_2=g'_1g'_2$ for some 
$g_1,g'_1\in Z_P(g)^0,g_2,g_2'\in Z_{P'}(g)^0$. Set $\tg_0=g'_1{}\i g_0g_1$. Then
$\tg_0\in Z_P(g)\cap Z_{P'}(g)\cap H^0$ hence 
$\tg_0\in L\cap H^0\sub Z_L(g)^0\sub Z_P(g)^0$. Thus, $g'_1{}\i g_0g_1\sub Z_P(g)^0$. 
Since $g_1,g'_1\in Z_P(g)^0$ it follows that $g_0\in Z_P(g)^0$. Since 
$g_0\in Z_G(g)^0\cap P$ was arbitrary we see that $Z_G(g)^0\cap P\sub Z_P(g)^0$ hence
$Z_G(g)^0\cap P=Z_P(g)^0$. This proves (c).

We prove (d). Let $T$
be a maximal torus of $Z_L(g)^0$. By the first line in 1.4, the (projective) variety of
Borel subgroups of $G^0$ that are normalized by $g$ is non-empty. Now $T$ acts by 
conjugation on this variety. Since $T$ is a torus this action must have a fixed point. 
Thus there exists a Borel $B$ of $G^0$ such that $g\in N_GB,T\sub B$. Since 
$(\cz_L\cap Z_G(g))^0\sub T$, we have $Z_{G^0}(T)\sub Z_{G^0}((\cz_L\cap Z_G(g))^0)=L$ 
(the last equality follows from 1.10). From $Z_{G^0}(T)\sub L$ we see that $T$ is also 
a maximal torus of $Z_G(g)^0$. Hence $T$ is a maximal torus of $Z_B(g)^0$. Let $R$ be a
connected component of $Z_B(g)$. If $x\in R$ then $x$ normalizes $Z_B(g)^0$ hence 
$xTx\i$ is a maximal torus of $Z_B(g)^0$ hence $xTx\i=x_0Tx_0\i$ for some 
$x_0\in Z_B(g)^0$. Then $x_0\i x\in R$ and $x_0\i x$ normalizes $T$. Since $T$ is a 
torus in the connected solvable group $B$, any element of $B$ that normalizes $T$ must 
centralize $T$. Thus $x_0\i x\in Z_{G^0}(T)\sub L$. We see that any connected component
of $Z_B(g)$ meets $L$. From the semidirect product decomposition $B=(B\cap L)U_P$ we 
obtain a semidirect product decomposition $Z_B(g)=Z_{B\cap L}(g)Z_{U_P}(g)$. It follows
that $Z_B(g)^0=Z_{B\cap L}(g)^0Z_{U_P}(g)^0$.

Let $y\in Z_{U_P}(g)$. Since the connected component of $Z_B(g)$ containing $y$ meets 
$L$ we have $y'y=y''$ for some $y'\in Z_B(g)^0$ and $y''\in L$ hence
$y''\in Z_{B\cap L}(g)$. We can write $y'=y'_1y'_2$ where $y'_1\in Z_{B\cap L}(g)^0$,
$y'_2\in Z_{U_P}(g)^0$. Hence $y'_1y'_2y=y''$ and 
$y'_2y=y'_1{}\i y''\in U_P\cap L=\{1\}$. Thus $y=y'_2{}\i\in Z_{U_P}(g)^0$. We see that
$Z_{U_P}(g)\sub Z_{U_P}(g)^0$ hence $Z_{U_P}(g)=Z_{U_P}(g)^0$. This proves (d).

We prove (e). By assumption we have $\dim U_P>0$. If $\kk$ has characteristic $p>1$ 
then $\Ad(g):U_P@>>>U_P$ has order a power of $p$. We can find a finite subgroup $U'$ 
of $U_P$ of order $p^n$ where $n>0$ such that $\Ad(g)(U')=U'$ (for example we can take 
$U'$ to be the group of rational points of $U_P$ over a suitable finite field). Since 
the cardinal of any orbit of $\Ad(g):U'@>>>U'$ is $1$ or a multiple of $p$ it follows 
that the number of fixed points of $\Ad(g):U'@>>>U'$ is divisible by $p$. Being 
$\ne 0$, it is $>1$. Thus $Z_{U_P}(g)$ has at least $2$ elements. Being connected (see 
(d)) it must have dimension $>0$. Assume now that $\kk$ has characteristic $0$. Let 
$U'$ be the centre of $U_P$. Then $\dim U'>0$ and under an isomorphism $U'\cong\kk^n$, 
$\Ad(g):U'@>>>U'$ becomes a unipotent automorphism of $\kk^n,n>0$ hence it has a fixed 
point set of dimension $>0$. Thus the inequality in (e) holds. The equality in (e) 
follows from (g) and (h).

\medpagebreak

The following result is a generalization of \cite{\IC, 2.8} to the disconnected case.

\proclaim{Proposition 10.2}Let $(C,\ce)$ be a cuspidal pair for $G$ with $\ce\ne 0$. 
Let $g\in C$. Then $Z_G(g)^0/(\cz_{G^0}\cap Z_G(g))^0$ is a unipotent group.
\endproclaim
Let $T$ be a maximal torus of $Z_G(g)^0$. Let $L=Z_{G^0}(T)$. We can find 
$\c\in\Hom(\kk^*,G^0)$ such that $\c(\kk^*)\sub T$ and $Z_{G^0}(\c(\kk)^*)=L$. Let
$P=P_\c$ (a parabolic of $G^0$ with Levi $L$, see 1.16). Now $g$ normalizes $P$ since 
it centralizes $\c(\kk^*)$. Similarly, $g$ normalizes $L$. We shall use notation in 
Lemma 10.1. In particular, $F$ is defined. Clearly, $F\sub C$ and $\ce|_F$ is a 
non-zero, $U_P$-equivariant local system (for the conjugation action of $U_P$). Since 
this action is transitive and the isotropy group of $g$ is $Z_{U_P}(g)$ which is 
connected by Lemma 10.1(d), we see that $\ce|_F\cong\bbq^m$ for some $m>0$. Hence 
$H^{2e}_c(F,\ce)\ne 0$ where $e=\dim F=\fra{1}{2}(\dim\la g\ra-\dim\la g\ra_L)$. (See 
Lemma 10.1(a).) Now there are only finitely many $G^0$-conjugacy classes 
$\boc_1,\do,\boc_t$ that meet $gU_P$ and are contained in $C$ (the semisimple part of 
an element in $\boc_i$ must be in the $G^0$-conjugacy class of $g_s$). These conjugacy
classes have the same dimension, and one of them is $\la g\ra$. Hence, by Lemma 4.2(a),
we have 
$$\dim(\boc_i\cap gU_P)\le\fra{1}{2}(\dim\la g\ra-\dim\la g\ra_L)=e.$$
Hence $\dim(C\cap gU_P)\le e$. Since $F$ is a closed irreducible subvariety of 
$C\cap gU_P$ of dimension $e$ it follows that $H^{2e}_c(C\cap gU_P)$ contains a 
subspace isomorphic to $H^{2e}_c(F,\ce)\ne 0$ hence $H^{2e}_c(C\cap gU_P)\ne 0$. Since 
$(C,\ce)$ is a cuspidal pair, it follows that $P=G^0$ hence $L=G^0$. Since 
$L=Z_{G^0}(T)$, it follows that $T\sub\cz_{G^0}$ hence $T\sub(\cz_{G^0}\cap Z_G(g))^0$.
It follows that $Z_G(g)^0/(\cz_{G^0}\cap Z_G(g))^0$ is a unipotent group.

\medpagebreak

The following result extends to the disconnected case results in \cite{\LSS},
\cite{\IC, \S7}.

\proclaim{Lemma 10.3}Let $P$ be a parabolic of $G^0$ with Levi $L$ and let $\boc_1$ be 
a unipotent $L$-conjugacy class in $N_GL\cap N_GP$. Let $\ph:X^\o@>>>\bY^\o$ be as in 
8.1.

(a) Let $\boc$ be the unique unipotent $G^0$-conjugacy class in
$G$ such that $\boc_1U_P\cap\boc$ is dense in $\boc_1U_P$. Then $\boc$ is the unique
unipotent $G^0$-conjugacy class in $G$ which is open dense in $\bY^\o$.

(b) $G^0$ acts transitively on $\ph\i(\boc)$.

(c) $P$ acts transitively on $\bboc_1U_P\cap\boc$. 

(d) We have $\bboc_1U_P\cap\boc=\boc_1U_P\cap\boc$.

(e) Let $g\in\boc_1U_P\cap\boc$. Define $\bg\in N_GL\cap N_GP$ by $\bg U_P=gU_P$. The 
natural map $\g:Z_P(g)/Z_P^0(g)@>>>Z_{G^0}(g)/Z_G(g)^0$ is injective and the natural 
map $\g':Z_P(g)/Z_P^0(g)@>>>Z_L(\bg)/Z_L(\bg)^0$ is surjective.
\endproclaim
We prove (a). As $\bY^\o$ is the union of the $G^0$-conjugates of $\bboc_1U_P$ which 
is contained in the closure of $\boc$, we see that $\bY^\o$ is contained in the 
closure of $\boc$. Since $\boc_1U_P\sub\bY^\o$ and $\boc_1U_P\cap\boc\ne\em$ we see 
that $\bY^\o$ contains some point of $\boc$ hence it contains $\boc$. Thus the 
closure of $\boc$ is equal to $\bY^\o$. Hence $\boc$ is open dense in $\bY^\o$. The 
uniqueness of such $\boc$ follows from the irreducibility of $\bY^\o$ (see 8.1(a)). 

We prove (b). Clearly, $\ph:X^\o@>>>\bY^\o$ is surjective and by 8.1(a),
$X^\o,\bY^\o$  are irreducible of the same dimension. Hence all fibres of 
$\ph:\ph\i(\boc)@>>>\boc$ are finite. Now $\ph$ maps any $G^0$-orbit in $\ph\i(\boc)$ 
onto $\boc$ since $G^0$ is transitive on $\boc$; hence any $G^0$-orbit on 
$\ph\i(\boc)$ must have dimension equal to $\dim\boc=\dim\bY^\o=\dim X^\o$ hence it 
is dense in $X^\o$. It follows that any two $G^0$-orbits on $\ph\i(\boc)$ must have 
non-empty intersection, so that there is only one orbit on $\ph\i(\boc)$ and (b) is 
proved.

We prove (c). Let $g,g'$ be two elements of $\bboc_1U_P\cap\boc$. Then $g'=x\i gx$ for 
some $x\in G^0$. Since $g'\in\bboc_1U_P$ it follows that $(g,xP)\in\ph\i(\boc)$. By 
(b), $(g,xP)$ must be in the same $G^0$-orbit as $(g,P)\in\ph\i(\boc)$. Hence there 
exists $y\in G^0$ such that $y\i gy=g$, $yP=xP$. Then $y=xz,z\in P$ and 
$g=z\i x\i gxz=z\i g'z$. Thus $g,g'$ are conjugate under $z\in P$ and (c) is proved.

We prove (d). Let $g\in\bboc_1U_P\cap\boc, g'\in\boc_1U_P\cap\boc$. By (c), $g$ is 
$P$-conjugate to $g'$. Since $\boc_1U_P\cap\boc$ is stable under $P$-conjugacy, it
follows that $g\in\boc_1U_P\cap\boc$. Thus, $\bboc_1U_P\cap\boc\sub\boc_1U_P\cap\boc$. 
The reverse inclusion is obvious and (d) is proved.

We prove (e). Since $\dim\ph\i(\boc)=\dim\boc$, the isotropy group of $g$ in $G^0$ has 
the same dimension as the isotropy group of $(g,P)$ in $G^0$. Thus $Z_{G^0}(g)$ has the
same dimension as its subgroup $Z_P(g)$. It follows that $Z_G(g)^0=Z_P(g)^0$. Hence
$\g$ is injective. We show that $\g'$ is surjective. By (a), $gU_P\cap\boc$ is open in
$gU_P$. Being non-empty, it is dense in $gU_P$. Hence $gU_P\cap\boc$ is irreducible.
From (c) we see that $Z_L(\bg)U_P$ acts transitively by conjugation on $gU_P\cap\boc$. 
Since $gU_P\cap\boc$ is irreducible it follows that $Z_L(\bg)^0U_P$ must also act 
transitively on $gU_P\cap\boc$. Hence for any element $z\in Z_L(\bg)$ there exist
$z_1\in Z_L(\bg),v\in U_P$ such that $zgz\i=z_1vgv\i z_1\i$ so that 
$v\i z_1\i z\in Z_P(g)$. Under $\g'$, the coset of $v\i z_1\i z$ is mapped to the coset
of $z\i z$ which is the same as the coset of $z$. Thus, $\g'$ is surjective. The lemma 
is proved.

\head 11. The structure of the algebra $\EE$\endhead
\subhead 11.1\endsubhead   
In this section we generalize results in \cite{\IC, \S9} to the disconnected case.

Let $D$ be a connected component of $G$. Let $(P,L,\boc_1,\cf_1)\in\tcm_D$. Let $\d$ be
the unique connected component of $\tL=N_GL\cap N_GP$ that is contained in $D$. Then 
$\boc_1\sub\d$. Let $S={}^\d\cz_L^0\boc_1$. Define $\ce\in\cs(S)$ by 
$\ce=\bbq\bxt\cf_1$.
Let $\cw_S$ be as in 3.13. Define $\EE$ in terms of $L,S,\ce,G$ as in 7.10. 
Let $\bY$ be the closure of $Y_{L,S}$ in $D$. Let 
$\ph:X@>>>\bY$ be as in 3.14. Let $X_S,\bce$ be as in 5.6. Let $K\in\cd(X)$ be as in 
5.7. Let $\boc$ be the unipotent $G^0$-conjugacy class in $D$ such that 
$\boc\cap\boc_1U_P$ is open dense in $\boc_1U_P$ (see Lemma 10.3). Let $\hboc$ be 
the (unipotent) $G^0$-conjugacy class in $D$ such that $\boc_1\sub\hboc$.

\proclaim{Lemma 11.2} (a) We have $\ch^{2d}(\ph_!K)|_{\hboc}\ne 0$ where 
$d=\nu-\nu_L-\fra{1}{2}(\dim\hboc-\dim\boc_1)$.

(b) We have $\ch^0(\ph_!K)|_\boc\ne 0$ and
$0=\nu-\nu_L-\fra{1}{2}(\dim\boc-\dim\boc_1)$.

(c) If $P\ne G^0$ then $\boc\ne\hboc$.

(d) If $P\ne G^0$ then $\dim\EE\ge 2$.
\endproclaim
We prove (a). Let $e>0$ be the rank of $\ce$.
Using the equivalence of (ii),(iii) in Lemma 8.3 we see that it is enough to show that
$H^{2d}_c(\ph\i(g)\cap X_S,\bce)\ne 0$ where $g\in\boc_1$. 
Note that $\dim(\ph\i(g)\cap X_S)\le d$ by 4.2(b). The irreducible variety
$V=Z^0_G(g)/(Z_G(g)^0\cap P)$ is contained in $\ph\i(g)\cap X_S$ by $i:x\m(g,xP)$ and 
has dimension $d$ (see Lemma 10.1(b)). Hence it is enough to show that $\bce|_V\cong
\bbq^e$. Consider the commutative diagram
$$\CD
Z_G^0(g)@>\hat i>>\hX_S@>f>>\boc_1\\
@VjVV        @Va'VV     @.\\
V@>i>>   X_S   @.  {}
\endCD$$
where $\hX_S=\{(g',x)\in G\T G^0; x\i gx\in SU_P\}$, $\hat i(x)=(g,x)$,
$j(x)=x(Z_G^0(g)\cap P)$, $a'(g',x)=(g',xP)$, $f(g',x)=\boc_1$-component of
$x\i gx\in{}^\d\cz_L^0\boc_1U_P$. By definition we have $a'{}^*\bce=f^*\cf_1$. Since 
$f\hat i$ maps $Z^0_G(g)$ to a point, the local system $\hat i^*f^*\ce=j^*i^*\bce$ on
$Z^0_G(g)$ is isomorphic to $\bbq^e$. Since $j$ is a principal bundle with connected 
group $Z_G^0(g)\cap P$ (see Lemma 10.1(c)) it follows that $i^*\bce\cong\bbq^e$. This 
proves (a).

We prove (b). Let $g\in\boc\cap\boc_1U_P$. To prove the first assertion of (b) it is 
enough to show that $H^0_c(\ph\i(g)\cap X_S,\bce)\ne 0$. We have $(g,P)\in\ph\i(\boc)$.
Since $g\in\boc_1U_P$, we have $(g,P)\in X_S$. Since $G^0$ acts transitively on 
$\ph\i(\boc)$ (see Lemma 10.3(b)) we see that the $G^0$-orbit of $(g,P)$, that is 
$\ph\i(\boc)$, is contained in $X_S$. In particular, $\ph\i(g)\sub X_S$. To show that 
$H^0_c(\ph\i(g)\cap X_S,\bce)\ne 0$ it is enough to show that $\ph\i(g)$ is a finite 
set. Since $G^0$ acts transitively on $\ph\i(\boc)$, we have a bijection 
$\ph\i(g)\cong Z_{G^0}(g)/Z_P(g)$. By the proof of Lemma 10.3(e), the groups 
$Z_{G^0}(g),Z_P(g)$
have the same identity component. Hence $Z_{G^0}(g)/Z_P(g)$ is finite. This proves the 
first assertion of (b). 

From Lemma 10.3(a) we see that $\dim\boc=\dim\bY^\o$. By 8.1(a) we have 
$\dim\bY^\o=2\nu-2\nu_L+\dim\boc_1$. Hence the equality in (b) holds. 

We prove (c). By Lemma 10.1(e) we have $\nu-\nu_L-\fra{1}{2}(\dim\hboc-\dim\boc_1)>0$. 
Combining this with the equality in (b) we see that $\dim\hboc<\dim\boc$. Hence 
$\boc\ne\hboc$.

We prove (d). Using (a),(b) and Lemma 8.3 we see that there exist irreducible 
$G^0$-equivariant local systems $\cf^1$ on $\hboc$ and $\cf^2$ on $\boc$ and 
$\r^1,\r^2\in\Irr\EE$ such that $(\hboc,\cf^1)=\g(\r^1)$, $(\boc,\cf^2)=\g(\r^2)$ 
(where $\g$ as in Lemma 8.2). Since $\hboc\ne\boc$ (see (c)) we have 
$\g(\r^1)\ne\g(\r^2)$. Since $\g$ is injective, we have $\r^1\ne\r^2$. Thus $\Irr\EE$ 
has at least two elements. Hence $\dim\EE\ge 2$. 

\proclaim{Lemma 11.3} The conjugation action of $\cw_S$ on ${}^\d\cz_L^0/{}^D\cz_G^0$ 
is faithful.
\endproclaim
Let $n\in G^0$ be such that $nLn\i=L,nSn\i=S$ and $nxn\i\in{}^D\cz_G^0x$ for any
$x\in{}^\d\cz_L^0$. We must show that $n\in L$. Let $\cl=\Hom(\kk^*,{}^\d\cz_L^0)$, a 
free abelian group of finite rank and let $\bar n:\cl@>>>\cl$ be the endomorphism 
induced by $\Ad(n):{}^\d\cz_L^0@>>>{}^\d\cz_L^0$. This endomorphism has finite order
since $\cw_S$ is a finite group. By our assumption, the endomorphism 
$\t:{}^\d\cz_L^0@>>>{}^\d\cz_L^0,x\m nxn\i x\i$ satisfies $\t^2(x)=1$ for all $x$.
Hence $(\bar n-1)^2=0$. Since $\bar n$ has finite order it follows that $\bar n=1$ 
hence $\t(x)=x$ that is, $nxn\i=x$ for all $x\in{}^\d\cz_L^0$. We see that
$n\in Z_{G^0}({}^\d\cz_L^0)=L$ (see 1.10(a)). The lemma is proved.

\subhead 11.4\endsubhead
We can find a unipotent element $u\in\d$ such that $u$ is quasi-semisimple in $\tL$. 
Then we can find a Borel $B_1$ of $L$ and a maximal torus $T$ of $B_1$ such that
$uB_1u\i=B_1,uTu\i=T$. Let $B=B_1U_P$, a Borel of $G^0$ contained in $P$. According to 
1.8, $\Ad(u)$ preserves some \'epinglage of $G^0$ attached to $(B,T)$. Using 1.7(b), we
see that $\cz_{Z_L(u)^0}^0={}^\d\cz_L^0,\cz_{Z_G(u)^0}^0={}^D\cz_G^0$ hence

(a) $\cz_{Z_L(u)^0}^0/\cz_{Z_G(u)^0}^0={}^\d\cz_L^0/{}^D\cz_G^0$.
\nl
We will use the terminology "$D$-parabolic" instead of "parabolic normalized by some 
element of $D$". 

\proclaim{Lemma 11.5} Assume that $P\ne G^0$ and that there is no $D$-parabolic $P'$ 
of $G^0$ such that $P\sub P'$, $P'\ne G^0$, $P'\ne P$.

(a) We have $\cw_S=\cw_\ce$ and $|\cw_S|=\dim\EE=2$.

(b) Let $Q$ be a parabolic of $G^0$ with Levi $L$ such that $Q$ is normalized by some
(or equivalently, any) element $g\in\d$. Then either $Q=P$ or $Q=\ti Q$, the unique 
parabolic of $G^0$ with Levi $L$, opposed to $P$.

(c) Let $Q'=nPn\i$ where $n\in N_{G^0}L$ represents the non-trivial element of $\cw_S$.
Then $Q'=\ti Q$.
\endproclaim
We prove (a). Using 7.10 and Lemma 11.2(d) we have $2\le\dim\EE=|\cw_\ce|\le|\cw_S|$. 
It is enough to show that $|\cw_S|\le 2$. Let $u$ be as in 11.4. Then $Z_P(u)^0$ is a 
maximal parabolic of $Z_G(g)^0$. (From 1.7(c) we see that $Z_P(u)^0$ is a proper 
parabolic of $Z_G(u)^0$ with Levi $Z_L(u)^0$. If it is not maximal then from 1.7(c) 
we see that there exists a parabolic $P'$ of $G^0$ such that $P'\ne G^0$, $P\sub P'$, 
$P'\ne P$ and $uP'u\i=P'$, contradicting our assumption.) It follows that 
$\cz_{Z_L(u)^0}^0/\cz_{Z_G(u)^0}^0$ is $1$-dimensional. Using this and 11.4(a) we 
obtain 

(d) $\dim({}^\d\cz_L^0/{}^D\cz_G^0)=1$.
\nl
Since the automorphism group of a $1$-dimensional torus has order $2$ and $\cw_S$ acts 
faithfully on the torus in (c), by Lemma 11.3, it follows that $|\cw_S|\le 2$. This 
proves (a).

We prove (b). We have $\fg=\op_{\mu\in\cx}\fg_\mu$ where 
$\cx:=\Hom({}^\d\cz_L^0/{}^D\cz_G^0,\kk^*)$ and $\fg_\mu$ denotes a weight space of the
adjoint action of ${}^\d\cz_L^0/{}^D\cz_G^0$. If $\mu$ is the trivial character, we 
have $\fg_\mu=\Lie L$ (see 1.10(a)). We show: 

(e) {\it there exists $\c\in\Hom(\kk^*,G^0)$ such that $\c(\kk^*)\sub{}^\d\cz^0_L$ and 
$P_\c=Q$ (see 1.16).}
\nl
(A variant of 1.17(a).) First we can find $\c'\in\Hom(\kk^*,G^0)$ such that 
$\c'(\kk^*)\sub\cz^0_L$ and $P_{\c'}=Q$. We can find $n\ge 1$ such that $g^n$ is in 
the identity component of $N_GQ\cap N_GL$, that is $g^n\in L$. Define 
$f:\cz^0_L@>>>\cz^0_L$ by $f(z)=gzg\i$. We have $f^n=1$. Define 
$\c_j:\kk^*@>>>G^0$ by $\c_j(a)=f^j(\c(a))$ for $j\in[0,n-1]$. Define
$\c\in\Hom(\kk^*,G^0)$ by $\c(a)=\c_0(a)c_1(a)\do c_{n-1}(a)$. Then 
$f(\c(a))=\c(a)$ hence $\c(a)\in Z_G(g)$ for all $a\in\kk^*$. We see that
$\c(\kk^*)\sub{}^\d\cz^0_L$. Since $\c_j(a)=g^j\c'(a)g^{-j}$, we have
$P_{\c_j}=g^jP_{\c'}g^{-j}=g^jQg^{-j}=Q$. Hence the $\kk^*$-action 
$a\m\Ad(\c_j(a))$ has $\ge 0$ weights on $\Lie Q$ and $<0$ weights on 
$\fg/\Lie Q$. Since these actions (for $j=0,1,\do,n-1$) commute with each other,
it follows that the $\kk^*$-action 
$a\m\Ad(\c(a))=\Ad(\c_0(a))\Ad(c_1(a))\do\Ad(c_{n-1}(a))$ has $\ge 0$ weights on 
$\Lie Q$ and $<0$ weights on $\fg/\Lie Q$. Hence $P_\c=Q$ and (e) is proved.

Let $\c$ be as in (e). Then $\c$ induces a homomorphism of $1$-dimensional tori 
$\kk^*@>>>{}^\d\cz_L^0/{}^D\cz_G^0$ which must be either constant or surjective; it is 
not constant since $\c(\kk^*)$ acts non-trivially on $\fg$, hence it is surjective. It 
follows that any weight space of $\Ad(\c(\kk^*))$ on $\fg$ coincides with one of the 
weight spaces $\fg_\mu$. In particular, $\Lie Q=\op_{\mu\in\cx'}\fg_\mu$ where $\cx'$ 
is the subset of $\cx$ which under an isomorphism $\ZZ@>\si>>\cx'$ corresponds to 
$\{0,1,2,\do\}$ or to $\{0,-1,-2,\do\}$. We see that there are only two possibilities 
for $Q$. This proves (b).

Applying (b) to $Q=Q'$, we see that $Q'$ must be either $P$ or $\ti Q$. Since 
$Q'\ne P$ we must have $Q'=\ti Q$. The lemma is proved.

\subhead 11.6\endsubhead
The $D$-parabolics of $G^0$ that contain strictly $P$
and are minimal with this property, form a finite set
$\{P_r;r\in\ci\}$. For each $r\in\ci$ let $L_r$ be the unique Levi of $P_r$ such 
that $L\sub L_r$; let $\boc^r$ be the unipotent $L_r$-conjugacy class in 
$\tL_r=N_GL_r\cap N_GP_r$ such that $\boc^r\cap\boc_1U_{P\cap L_r}$ is open dense in
$\boc_1U_{P\cap L_r}$ (see Lemma 10.3); let $\hboc^r$ be the (unipotent) 
$L_r$-conjugacy class in $\tL_r$ such that $\boc_1\sub\hboc^r$. By Lemma 11.5(a) 
applied to $\tL_r$ instead of $G$, we see that $\{n\in N_{L_r}L;nSn\i=S\}/L$ has 
order $2$; let $s_r$ be the unique non-trivial element of this group.

\proclaim{Proposition 11.7} Let $N_{G^0}\d=\{x\in G^0;x\d x\i=\d\}$, a subgroup of
$N_{G^0}L$.

(a) The inclusions $\cw_{\ce}\sub\cw_S\sub N_{G^0}\d/L$ are equalities. 

(b) $N_{G^0}\d/L$ is a Coxeter group with simple reflections $\{s_r;r\in\ci\}$.
\endproclaim
Since $s_r\in\cw_{\ce}$ (by Lemma 11.5 for $\tL_r$ instead of $G$), the subgroup of
$N_{G^0}\d/L$ generated by $\{s_r;r\in\ci\}$ is contained in $\cw_{\ce}$. Thus, (a) 
is a consequence of (b). It remains to prove (b).

Let $u,B_1,B,T$ be as in 11.4. Let $\{P^i;i\in I\}$ be the parabolic subgroups of 
$G^0$ that contain $B$ and are minimal with this property. Let $W=N_{G^0}T/T$. For 
$i\in I$ let $s^i$ be the unique non-trivial element in the image of the inclusion 
$N_{P^i}T/T@>>>W$. Then $W$ is a Coxeter group with simple reflections 
$\{s^i;i\in I\}$. Let $F:W@>>>W$ be the automorphism induced by 
$\Ad(u):N_{G^0}T@>>>N_{G^0}T$. There is a unique bijection $F:I@>>>I$ such that
$F(s^i)=s^{F(i)}$ for $i\in I$ or equivalently, $uP^iu\i=P^{F(i)}$. For any subset 
$J$ of $I$ let $P^J$ be the subgroup of $G^0$ generated by $\{P^i;i\in J\}$ and let
$W_J$ be the subgroup of $W$ generated by $\{s^i;i\in J\}$. The condition that $P^J$
is a $D$-parabolic is equivalent to the condition that $F(J)=J$. 
Hence the parabolics $\{P_r;r\in\ci\}$ are precisely the parabolics $P^J$ where $J$ 
is an $F$-orbit on $I-K$ where $K\sub I$ is defined by $P=P^K$ (we have $F(K)=K$). 
Thus we may identify $\ci$ with the set of $F$-orbits on $I-K$ so that $P_r=P^r$ for
any $r\in\ci$. We may identify canonically $N_{G^0}L/L$ with $N_W(W_K)/W_K$ and 
$N_{G^0}\d/L$ with the fixed point set $(N_W(W_K)/W_K)^F$ of 
$F:N_W(W_K)/W_K@>>>N_W(W_K)/W_K$. Now $s_r$ is a non-trivial element in the fixed 
point set of $F$ on $N_{W_{K\cup r}}(W_K)/W_K$. Applying Lemma 11.5(b) to $\tL_r$ 
instead of $G$ we see that $s_r$ is the $W_K$-coset of the longest element of 
$W_{K\cap r}$ (which must in our case normalize $W_K$). We can now apply 
\cite{\COX, 5.9(i)}; we see that (b) holds. (Note that \cite{\COX, 5.9} has an 
additional assumption, namely that for any $F$-stable subset $J$ of $I$ that 
contains $K$, the longest element of $W_J$ normalizes $W_K$. However in the proof of
\cite{\COX, 5.9(i)} that assumption is only used for $J=K\cup r$, $r$ as above.) The
proposition is proved.

\proclaim{Lemma 11.8} (a) The local system $\ch^0(\ph_!K)|_\boc$ is  
irreducible. Hence there is a unique $\r\in\Irr\EE$ such that
$\ch^0(\ph_!K)|_\boc=\ch^0(\ph_!K)_\r|_\boc$. 

(b) We have $\dim\r=1$.
\endproclaim
For $r\in\ci$ let $\EE^r$ be the algebra defined like $\EE$ (for $\tL_r$ instead of
$G$). We have naturally $\EE^r\sub\EE$ and $\dim\EE^r=2$ (see Lemma 11.5(a)). We 
have $\EE^r=\EE_1\op\EE_{s_r}$ (notation of 7.10). Let $D_r$ be the unique connected
component of $\tL_r$ that is contained in $D$. Using Lemma 11.2 for $\tL_r$ instead
of $G$ we see that $\Ph_{D_r}\i(L,\boc_1,\cf_1)$ consists of two elements; one is of
the form $(\boc^r,?)$, the other of the form $(\hboc^r,?)$. They correspond (as in 
8.9 for $\tL_r$ instead of $G$) to $\r_r,\hat\r_r$ (respectively) in $\Irr\EE^r$. We 
have $\Irr\EE^r=\{\r_r,\hat\r_r\}$.

By Lemma 11.2(b), we have $\ch^0(\ph_!K)_\boc\ne 0$. Hence we can find 
$\r\in\Irr\EE$ such that $\ch^0(\ph_!K)_\r|_\boc\ne 0$. 

If $\Hom_{\EE^r}(\hat\r_r,\r)\ne 0$, then from the equivalence of (i),(iv) in 
Proposition 9.4 we see that $\ch^{2n}({}_{R'}\ph''_!\tcf')\ne 0$ hence $n\ge 0$ 
where $2n=2\nu-2\nu_{L_r}-\dim\boc+\dim\hboc_r$. However, from 
$$0=\nu-\nu_L-\fra{1}{2}(\dim\boc-\dim\boc_1),\qua
0=\nu_{L_r}-\nu_L-\fra{1}{2}(\dim\boc^r-\dim\boc_1)$$
(see Lemma 11.2(b) for $G$ and $\tL_r$) and $\dim\hboc^r<\dim\boc^r$ 
(as in the proof of Lemma 11.2(c)) we see that $n<0$. This contradiction shows that
$\Hom_{\EE^r}(\hat\r_r,\r)=0$. It follows that $\r|_{\EE^r}$ is a direct sum of 
copies of $\r_r$. Hence if $b_r$ is a basis element of $\EE_{s_r}$ then $b_r$ acts 
on $\r$ as a scalar times the identity. Since $\EE_w\EE_{w'}=\EE_{ww'}$ for 
$w,w'\in\cw_{\ce}=\cw_S$ and $\{s_r,r\in\ci\}$ generates $\cw_S$ we see that any 
element of $\EE$ acts on $\r$ as a scalar times the identity. Since $\r$ is 
irreducible, it must be one dimensional and $\r|_{\EE^r}=\r_r$ for any $r$. This 
last property determines uniquely $\r$. The lemma is proved.

\proclaim{Proposition 11.9} There is a unique isomorphism of algebras 
$\EE@>\si>>\bbq[\cw_{\ce}]$ which maps $\EE_w$ onto $\bbq w$ for any $w\in\cw_{\ce}$
and makes $\r$ in Lemma 11.8 correspond to the unit representation of $\cw_{\ce}$.
\endproclaim
In each $\EE_w$ we choose as basis element $b_w$ the unique element which acts as 
the identity on the $\EE$-module $\r$. It is clear that $b_wb_{w'}=b_{ww'}$ and 
$(b_w)$ provides the required isomorphism.

\subhead 11.10\endsubhead
Using Proposition 11.9 we can reformulate the results in 8.9 as a bijection
$$\cn_D@>\si>>\sqc_{(L,\boc_1,\cf_1)}\Irr N_{G^0}\d/L\tag a$$
where $(L,\boc_1,\cf_1)$ runs over a set of representatives for the $G^0$-orbits on
$\cm_D$ and $\Irr N_{G^0}\d/L$ is the set of (isomorphism classes of) irreducible
representations of the Coxeter group $N_{G^0}\d/L$. This is called the {\it 
generalized Springer correspondence}.

\head 12. Classification of objects in $\cn_D^0$\endhead
\subhead 12.1\endsubhead
Let $D$ be a connected component of $G$ that contains some unipotent elements. We 
would like to classify the objects in $\cn_D^0$. We will make a number of reductions.

Replacing $G$ by the subgroup generated by $D$ (with the same identity component as 
$G$) we see that $\cn_D^0$ does not change; hence we may assume that 

(a) $G/G^0$ is cyclic with generator defined by $D$.
\nl
Then $G/G^0$ is a unipotent group. 

\subhead 12.2\endsubhead
In the remainder of this section we assume that 12.1(a) holds. Let 
$\p:G@>>>G_{ss}=G/\cz_{G^0}^0$ be the obvious map. Let $D'=\p(D)$ (a connected 
component of $G_{ss}$). We have $D=\p\i(D')$. Let $\boc'$ be a unipotent 
$G_{ss}^0$-conjugacy class in $D'$. Let 
$\boc=\{g\in\p\i(\boc'),g\text{ unipotent}\}$. We show:

(a) $\boc$ {\it is a single unipotent $G^0$-conjugacy class in $D$;}

(b) {\it for $g\in\boc$, the homomorphism $Z_{G^0}(g)@>>>Z_{G_{ss}^0}(\p(g))$ 
induced by $\p$ is surjective and its kernel is connected; hence the resulting 
homomorphism 

$Z_{G^0}(g)/Z_{G^0}(g)^0@>>>Z_{G_{ss}^0}(\p(g))/Z_{G_{ss}^0}(\p(g))^0$
\nl
is an isomorphism.}
\nl
If $g\in\p\i(\boc')$ then $g_u\in\boc$. Thus $\boc$ is non-empty. Let $g,g'\in\boc$.
Since $\p(g),\p(g')$ are $G_{ss}^0$-conjugate, we see that there exists
$z\in\cz_{G^0}^0$ such that $g',zg$ are $G^0$-conjugate. Using 1.3(a) we can write 
$z=xtgx\i g\i$ with $t,x\in\cz_{G^0}^0,tg=gt$. Then $zg=xtgx\i$ is $G^0$-conjugate 
to $tg$ hence $g',tg$ are $G^0$-conjugate. In particular, $tg$ is unipotent. Since 
$tg=gt$ with $t$ semisimple, $g$ unipotent, we see that $t=1$ hence $g,g'$ are 
$G^0$-conjugate. This proves (a).

We prove (b). Let $y\in G^0$ be such that $yg=zgy$ for some $z\in\cz_{G^0}^0$. Using 
again 1.3(a) we can write $z=xtgx\i g\i$ with $t,x\in\cz_{G^0}^0,tg=gt$. Then 
$ygy\i=zg=xtgx\i$ is unipotent hence $tg$ is unipotent. As in the proof of (a) we 
deduce that $t=1$. Hence $ygy\i=xgx\i$. Setting $y'=x\i y$ we have $y'\in Z_{G^0}(g)$ 
and $y=xy'\in\cz_{G^0}^0Z_{G^0}(g)$. This shows surjectivity. The kernel of 
$Z_{G^0}(g)@>>>Z_{G_{ss}^0}(\p(g))$ is $\cz_{G^0}^0\cap Z_G(g)$. This group is 
connected by \cite{\BO, 9.6} applied to the unipotent automorphism $\Ad(g)$ of the 
torus $\cz_{G^0}^0$. This proves (b).

From (a),(b), we see that we have a bijection $\cn_{D'}@>\si>>\cn_D$ given by
$(\boc',\cf')\m(\boc,\cf)$ where $\boc$ is as above and $\cf$ is the inverse image
of $\cf'$ under $\boc@>>>\boc'$ (restriction of $\p$). 

A standard argument (similar to one in the proof of Lemma 6.4) shows that this 
bijection restricts to a bijection $\cn_{D'}^0@>\si>>\cn_D^0$.

\subhead 12.3\endsubhead
Assume that $G$ is such that $G^0$ is semisimple. We can find a reductive group 
$\tG$ with $\tG^0$ semisimple, simply connected, and a surjective homomorphism of 
algebraic groups $\p:\tG@>>>G$ such that $\Ker\p\sub\cz_{\tG^0}$. Then 
$\tG^0=\p\i(G^0)$ and $\tD=\p\i(D)$ is a connected component of $\tG$. Let 
$(\boc,\cf)\in\cn_D$. Let $\tbc=\{g\in\p\i(\boc);g\text{ unipotent}\}$. We show:

(a) {\it $\tbc$ is a single unipotent $\tG^0$-conjugacy class in $\tD$.}

(b) {\it if $g\in\tbc$, then the obvious homomorphism $Z_{\tG^0}(g)@>>>Z_{G^0}(\p(g))$ 
is surjective.}
\nl
If $g\in\p\i(\boc)$ then $g_u\in\p\i(\boc)$. Thus, $\tbc\ne\em$. Let $g,g'\in\tbc$. 
We can find $x\in\tG^0$ such that $g'=xgx\i z$ for some $z\in\Ker\p$. Since 
$\Ad(g\i)$ is an automorphism of $\Ker\p$ of order relatively prime to $|\Ker\p|$, any 
element of $\Ker\p$ is of the form $\Ad(g\i)(z_1)z_1\i z_2$ where $z_1,z_2\in\Ker\p$ 
and $\Ad(g\i)z_2=z_2$. In particular, $z=g\i z_1gz_1\i z_2$ with
$z_1,z_2$ as above. Then $gz=z_1gz_1\i z_2$. Since $gz=x\i g'x$ is unipotent, we see
that $z_1gz_1\i z_2$ is unipotent. Now $z_2$ is semisimple and it commutes with 
$z_1gz_1\i$ which is unipotent. By the uniqueness of Jordan decomposition we have 
$z_2=1$ and $gz=z_1gz_1\i=x\i g'x$. Since $xz_1\in\tG^0$ we see that $g,g'$ are 
$\tG^0$-conjugate. This proves (a).

We prove (b). Let $x\in\tG^0$ be such that $xg=gxz$ for some $z\in\ker\p$. It is 
enough to show that $xz_1$ commutes with $g$ for some $z_1\in\ker\p$. We write 
$z=g\i z_1gz_1\i z_2$ as above. As in the proof of (a) (with $g'=g$) we see that 
$gz=z_1gz_1\i=x\i gx$. Hence $xz_1$ commutes with $g$, as required.

From (a),(b) it follows that any object $(\boc,\cf)\in\cn_D$ gives rise to an object
$(\tbc,\tcf)\in\cn_{\tD}$ where $\tbc$ is as above and $\tcf$ is the inverse image 
of $\cf$ under the map $\tbc@>>>\boc$ induced by $\p$. (The local system $\tcf$ is 
irreducible by (b).) From (b) we see also that the resulting map
$\cn_D@>>>\cn_{\tD}$ is injective.

We show:

(c) {\it for $(\boc,\cf)\in\cn_D$ we have $(\boc,\cf)\in\cn_D^0$ if and only if
$(\tbc,\tcf)\in\cn_{\tD}^0$.}
\nl
Assume first that $(\boc,\cf)\in\cn_D^0$. Let $P'$ be a proper parabolic of $\tG^0$
and let $g\in\tbc\cap N_{\tG}P'$. Let $\dd'$ be the $P'/U_{P'}$-conjugacy class 
of the image of $g$ in $N_{\tG}P'/U_{P'}$. Then $P:=\p(P')$ is a proper parabolic of
$G^0$ and $\p(g)\in\boc\cap N_GP$. Let $\dd$ be the $P/U_P$-conjugacy class of the 
image of $\p(g)$ in $N_GP/U_P$. By assumption we have
$H^{\dim\boc-\dim\dd}_c(\boc\cap\p(g)U_P,\cf)=0$. Now the morphism 
$\p_0:\tbc\cap gU_{P'}@>>>\boc\cap\p(g)U_P$ induced by $\p$ is an isomorphism. (We 
show only that it is a bijection. Let $u\in U_P$ be such that $\p(g)u\in\boc$. We 
can find a unique $u'\in U_{P'}$ such that $\p(u')=u$. Then $gu'\in\p\i(\boc)$ and 
$gu'$ is unipotent, since $g$ normalizes $P'$, hence $gu'\in\tbc$. Since
$\p(gu')=\p(g)u$ we see that $\p_0$ is surjective. The injectivity follows from the
fact that $gU_{P'}@>>>\p(g)U_P$ is an isomorphism.) Also, $\dim\boc=\dim\tbc$ and
$\dim\dd'=\dim\dd$. It follows that 
$H^{\dim\tbc-\dim\dd'}_c(\tbc\cap gU_{P'},\tcf)=0$ so that
$(\tbc,\tcf)\in\cn_{\tD}^0$. The proof of the reverse implication follows 
essentially the same argument, in the reverse.

It is easy to see that the kernel of the obvious homomorphism

$Z_{\tG^0}(g)/Z_{\tG^0}(g)^0@>>>Z_{G^0}(\p(g))/Z_{G^0}(\p(g))^0$
\nl
is a homomorphic image of $\ker\p$. It follows that the image of $\cn_D@>>>\cn_{\tD}$ 
consists of all pairs $(\boc',\cf')\in\cn_{\tD}$ such that the natural action of 
$\ker\p$ on each fibre of $\cf'$ is trivial.

Thus the objects in $\cn_D^0$ are in natural bijection with the objects in
$\cn_{\tD}^0$ with trivial action of $\ker\p$.

\subhead 12.4\endsubhead
Next we assume that $G$ is such that $G^0$ is semisimple, simply connected. We can 
write uniquely $G^0$ as a product $G^0=G_1\T G_2\T\do\T G_k$ where each $G_i$ is a 
closed connected normal subgroup of $G$ different from $\{1\}$ and minimal with 
these properties. For $i\in[1,k]$, let 
$G'_i=G/(G_1\T\do\T G_{i-1}\T G_{i+1}\T\do\T G_k)$. Then $G'_i$ is a reductive group
with $G'_i{}^0=G_i$ and the image of $D$ in $G'_i$ is a connected component $D_i$ of
$G'_i$. Also we have an obvious homomorphism $G@>>>G'_1\T G'_2\T\do\T G'_k$ which is
an imbedding of algebraic groups by which we identify $G$ with a closed subgroup of
$G'_1\T G'_2\T\do\T G'_k$ with the same identity component; then $D$ becomes 
$D_1\T D_2\T\do\T D_k$. From the definitions it is clear that we have a natural 
bijection $\cn_{D_1}\T\cn_{D_2}\T\do\T\cn_{D_k}@>\si>>\cn_D:$
$$((\boc_1,\cf_1),(\boc_2,\cf_2),\do,(\boc_k,\cf_k))\m
(\boc_1\T\boc_2\T\do\T\boc_k,\cf_1\bxt\cf_2\bxt\do\bxt\cf_k)$$
and this restricts to a bijection
$$\cn_{D_1}^0\T\cn_{D_2}^0\T\do\T\cn_{D_k}^0@>\si>>\cn_D^0.$$

\subhead 12.5\endsubhead
Next we assume that $G$ is such that $G^0$ is semisimple, simply connected,
$\ne\{1\}$ and that $G$ has no closed connected normal subgroups other than $G^0$ 
and $\{1\}$. (If $G^0=\{1\}$ then $\cn_D=\cn_D^0$ consists of a single object 
$(D,\bbq)$.) We have $G^0=H_0\T H_1\T\do \T H_{m-1}$ where $H_i$ are connected, 
simply connected, almost simple, closed subgroups of $G^0$. Let $\g\in D$ be a 
unipotent quasi-semisimple element. We can assume that $H_i=\g^i H_0\g^{-i}$ for 
$i\in[1,m-1]$ and $\g^mH_0\g^{-m}=H_0$. Let $G'$ be the subgroup of $G$ generated by
$H_0$ and $\g^m$. Since $\g$ has finite order, $G'$ is closed, $G'{}^0=H_0$ and 
$D'=\g^mH_0$ is a connected component of $G'$. Consider the diagram
$$D'@<a<<G^0\T D'@>b>>D$$
where $a(g,\g^mh)=\g^mh,b(g,\g^mh)=g\g hg\i$ (with $h\in H_0$). Now $G^0\T H_0$ acts 
on $G^0\T D'$ by 
$$\align&(y,u_0):(h_0h_1\do h_{m-1},\g^mh)\m\\&
(yh_0u_0\i h_1\g^{-m+1}u_0\i\g^{m-1}h_2\g^{-m+2}u_0\i\g^{m-2}\do 
h_{m-1}\g\i u_0\i\g,u_0\g^mhu_0\i),\\&\quad (h_i\in H_i),\endalign$$
on $D'$ by $(y,u_0):\g^mh\m u_0\g^mhu_0\i$ and on $D$ by $(y,u_0):d\m ydy\i$; 
moreover, $a$ and $b$ are $G^0\T H_0$ equivariant. Note also that $a$ is a principal
$G^0$-bundle and $b$ is a principal $H_0$-bundle. Hence $a,b$ induce bijections
$$\align&\text{set of $H_0$-conjugacy classes in $D'$}@>a\i>>
\text{set of $G^0\T H^0$-orbits in $G^0\T D'$},\\&
\text{set of $G^0$-conjugacy classes in $D$}@>b\i>>
\text{set of $G^0\T H^0$-orbits in $G^0\T D'$}.\tag a\endalign$$
Moreover, if $h\in H_0$ then $a,b$ induce isomorphisms
$$Z_{H^0}(\g^mh)@<\si<<(\text{stabilizer of $(1,\g h)$ in } G^0\T H_0)@>\si>>
Z_{G^0}(\g h).\tag b$$
We show:

(c) {\it an $H_0$-conjugacy class in $D'$ is unipotent if and only if the 
$G^0$-conjugacy class in $D$ corresponding to it by (a) is unipotent.}
\nl
It is enough to show that for $h\in H_0$ we have $\g h$ unipotent if and only if
$\g^mh$ is unipotent. This is trivial if $m=1$. Assume now that $m\ge 2$. Then the 
characteristic exponent of $\kk$ is $>1$ and $m$ is a power of it. It is enough to 
show that the conditions $(\g h)^{m^k}=1$ for some $k>0$ and $(\g^mh)^{m^k}=1$ for 
some $k>0$ are equivalent. For $k>0$ we have
$$\align(\g h)^{m^k}&=(\g h\g\i)(\g^2h\g^{-2})\do(\g^{m^k}h\g^{-m^k})\\&=
x_k(\g x_k\g\i)(\g^2x_k\g^{-2})\do(\g^{m-1}x_k\g^{-m+1})\endalign$$
where $x_k=(\g h\g\i)(\g^{m+1}h\g^{-m-1})\do(\g^{m^k-m+1}h\g^{-m^k+m-1})$. Since 

$x_k\in H_1,\g x_k\g\i\in H_2,\g^2x_k\g^{-2}\in H_2,\do,\g^{m-1}x_k\g^{-m+1}\in H_0$,
\nl
the condition that $(\g h)^{m^k}=1$ is equivalent to the condition that 

$x_k=1,\g x_k\g\i=1,\do,\g^{m-1}x_k\g^{-m+1}=1$,
\nl
that is to the condition that $x_k=1$. On the other hand we have 
$(\g^mh)^{m^k}=\g^{m-1}x_{k+1}\g^{-m+1}$ and this is $1$ if and only if $x_{k+1}=1$. 
This proves (c).

From (a),(b),(c) we see that there is a natural bijection $\cn_D\lra\cn_{D'}$ in 
which $(\boc,\cf)\in\cn_D$ corresponds to $(\boc',\cf')\in\cn_{D'}$ when 
$b\i(\boc)=a\i(\boc')$ and the inverse image of $\cf$ under $b:b\i(\boc)@>>>\boc$ 
coincides with the inverse image of $\cf'$ under $a:a\i(\boc')@>>>\boc'$.

(d) {\it If $\boc,\boc'$ are as above, then the principal $H_0$-bundle 
$b:G^0\T\boc'@>>>\boc$ restricts to an isomorphism 
$b':H_1H_2\do H_{m-1}\T\boc'@>\sim>>\boc$.}
\nl
We show only that $b'$ is bijective. This follows from the fact that
$H_1H_2\do H_{m-1}\T D'$ meets each $H_0$-orbit on $G^0\T D'$ in exactly one point.

We show:

(e) {\it for $(\boc,\cf),(\boc',\cf')$ related as above, we have 
$(\boc,\cf)\in\cn_D^0$ if and only if $(\boc',\cf')\in\cn_{D'}^0$.}
\nl
Assume first that $(\boc',\cf')\in\cn_{D'}^0$. Let $P$ be a proper parabolic of 
$G^0$ and let $g\in\boc\cap N_GP$. We must show that $(\boc,\cf)$ satisfies the 
criterion 8.7(ii) with respect to $P,g$. Replacing $P,g$ by a $G^0$-conjugate we may
assume that $g=\g h$ where $h\in H_0,\g^mh\in\boc'$. We have $P=P_0P_1\do P_{m-1}$ 
where $P_i$ is a parabolic of $H_i$ for each $i$. Since $P$ is normalized by $g$ we
see that

$\g hP_0h\i\g\i=P_1,\g P_1\g\i=P_2,\do,\g P_{m-2}\g\i=P_{m-1},\g P_{m-1}\g\i=P_0$.
\nl
Thus,

$P=P_0(\g^{-m+1}P_0\g^{m-1})\do(\g\i P_0\g)$ and $\g^mh\in N_{G'}P_0$.
\nl
In particular, since $P\ne G^0$, we must have $P_0\ne H_0$.
Using (d) we can identify $\boc\cap\g hU_P$ with

$\{(x,\g^mh')\in H_1H_2\do H_{m-1}\T\boc';x\g h'x\i\in\g hU_P\}$
\nl
hence with
$$\align\{(x_1,x_2,\do x_{m-1},&\g^mh')\in H_1\T H_2\T\do\T H_{m-1}\T\boc';\\&
x_1x_2\do x_{m-1}\g h'x_1\i x_2\i\do x_{m-1}\i\in\g hU_P\}\endalign$$
and also with
$$\align&\{(x_1,x_2,\do x_{m-1},\g^mh',u_1,u_2,\do,u_{m-1})\\&
\in H_1\T H_2\T\do\T H_{m-1}\T\boc'\T U_{P_1}\T U_{P_2}\T\do\T U_{P_{m-1}};\\&
x_1x_2\do x_{m-1}\g h'x_1\i x_2\i\do x_{m-1}\i\in\g hU_{P_0}u_1u_2\do u_{m-1}\}.
\endalign$$
The last condition can be rewritten as 
$$\align&\g\i x_1\g h'\in hU_{P_0},\g\i x_2\g x_1\i=u_1,\g\i x_3\g x_2\i=u_2,\do,\\&
\g\i x_{m-1}\g x_{m-2}\i=u_{m-2},x_{m-1}\i=u_{m-1}\endalign$$
or as
$$\align&x_{m-1}=u_{m-1}\i,x_{m-2}=u_{m-2}\i\g\i u_{m-1}\i\g,\do,\\&
x_1=u_1\i\g\i u_2\i\g\do\g^{-m+2}u_{m-1}\i\g^{m-2},\\&
\g\i u_1\i\g   \g^{-2}u_2\i\g^2 \do\g^{-m+1}u_{m-1}\i\g^{m-1}h'\in hU_{P_0}.\tag f
\endalign$$
We have
$$\align&\g\i u_1\i\g\in\g^{-m}U_{P_0}\g^m,\g^{-2}u_2\g^2\in\g^{-m}U_{P_0}\g^m,
\do,\\&\g^{-m+1}u_{m-1}\i\g^{m-1}\in\g^{-m}U_{P_0}\g^m\endalign$$
and $\g^{-m}U_{P_0}\g^m h=hU_{P_0}$ hence the last condition in (f) is equivalent to
$h'\in hU_{P_0}$.
Hence in our variety we can drop the variables $x_1,x_2,\do,x_m$ and we obtain
$$\{(\g^mh',u_1,u_2,\do,u_{m-1})
\in \boc'\T U_{P_1}\T U_{P_2}\T\do\T U_{P_{m-1}};h'\in hU_{P_0}\}.$$
Clearly the first projection makes this last variety an affine space bundle over
$\boc'\cap\g^mh'U_{P_0}$. We see that $\boc\cap\g h'U_P$ is an affine space bundle 
over $\boc'\cap\g^mh'U_{P_0}$ with fibres of dimension $(m-1)\dim U_{P_0}$.

Let $\dd$ be the $P/U_P$-conjugacy class of the image of $\g h$ in $N_GP/U_P$. Let 
$\dd'$ be the $P_0/U_{P_0}$-conjugacy class of the image of $\g^m h$ in 
$N_{G'}P_0/U_{P_0}$.

From (d) we see that $\dim\boc=\dim\boc'+(m-1)\dim H_0$. Similarly, 
$\dim\dd=\dim\dd'+(m-1)\dim(P_0/U_{P_0})$. Hence 

$\dim\boc-\dim\dd=\dim\boc'-\dim\dd'+2(m-1)\dim U_{P_0}$. 
\nl
We now see that
$$\align H^{\dim\boc-\dim\dd}_c(\boc\cap\g hU_P,\cf)&=
H^{\dim\boc'-\dim\dd'+2(m-1)\dim U_{P_0}}_c(\boc\cap\g hU_P,\cf)\\&=
H^{\dim\boc'-\dim\dd'}_c(\boc'\cap\g^mhU_{P_0},\cf').\endalign$$
But the last cohomology space is $0$ since $(\boc',\cf')\in\cn_{D'}^0$. Hence the
first cohomology group is $0$. Thus, $(\boc,\cf)\in\cn_D^0$. 

Conversely, assume that $(\boc,\cf)\in\cn_D^0$. Let $P_0$ be a proper parabolic of 
$H_0$ and let $h\in H_0$ be such that $\g^mh\in\boc'\cap N_{G'}P_0$. We have 
$\g h\in\boc$. Now $P=P_0(\g^{-m+1}P_0\g^{m-1})\do(\g\i P_0\g)$ is a proper parabolic 
of $G^0$ such that $\g h\in N_GP$; the earlier argument can be applied in the reverse 
and gives \lb 
$H^{\dim\boc'-\dim\dd'}_c(\boc'\cap\g^mhU_{P_0},\cf')=0$ ($\dd'$ as above.)
Thus, $(\boc',\cf')\in\cn_{D'}^0$. 

\subhead 12.6\endsubhead
Next we assume that $G$ is such that $G^0$ is semisimple, simply connected, almost 
simple. Let $\D$ be the set of unipotent elements in $\cz_G$. Let $G'=G/\D$ and let
$\p:G@>>>G'$ be the obvious homomorphism. Since $\D\cap G^0=\{1\}$, the restriction
of $\p$ to $G^0$ is injective; it is in fact an isomorphism $G^0@>\si>>G'{}^0$. We 
show that 

(a) $\cz_{G'}\sub G'{}^0$.
\nl
Assume that $a\in G$ is such that $aba\i b\i\in\D$  for all $b\in G$. It is enough 
to show that the image $a'$ of $a$ in $G'$ is in $G'{}^0$. For $b\in G^0$ we have 
$aba\i b\i\in\D\cap G^0=\{1\}$ hence $a\in Z_G(G^0)$. Let $\g$ be a unipotent 
quasi-semisimple element in $D$. Then $a=\g^k a_0$ where $a_0\in G^0,k\in\NN$. It 
follows that $\Ad(\g^k):G^0@>>>G^0$ is an inner automorphism. Since 
$\Ad(\g^k):G^0@>>>G^0$ preserves an \'epinglage of $G^0$ it follows that 
$\g^k\in Z_G(G^0)$. Hence $a_0\in Z_G(G^0)$ and  $a_0\in\cz_{G^0}$. Since $\g^k$ 
commutes with all elements in $G^0$ and with $\g$, we have $\g^k\in\cz_G$. Since 
$\g^k$ is unipotent, we see that $\g^k\in\D$. Thus, $a\in\D G^0$. Hence 
$a'\in G'{}^0$. This proves (a).

Let $D'=\p(D)$, a connected component of $G'$. Then $\p$ restricts to an isomorphism
$D@>\si>>D'$. Let $\boc'$ be a unipotent $G'{}^0$-conjugacy class in $D'$. Let 
$\boc=D\cap\p\i(\boc')$. Then $\boc$ is a single unipotent $G^0$-conjugacy class in
$D$. For $g\in\boc$, the obvious homomorphism $Z_{G^0}(g)@>>>Z_{G'{}^0}(\p(g))$ is
an isomorphism. Hence there is a natural bijection $\cn_{D'}\lra\cn_D$ in which 
$(\boc',\cf')\in\cn_{D'}$ corresponds to $(\boc,\cf)\in\cn_D$ when 
$\boc=D\cap\p\i(\boc')$ and $\cf$ is the inverse image of $\cf'$ under
$\p:\boc@>>>\boc'$. It is clear that this bijection restricts to a bijection
$\cn_{D'}^0\lra\cn_D^0$.

\subhead 12.7\endsubhead
The results in 12.1-12.6 reduce the problem of classifying the objects in $\cn_D^0$
to the special case where $G$ is such that $G^0$ is semisimple, almost simple, that 
$\cz_G\sub G^0$ and that $D$ is a generator of $G/G^0$. In the
remainder of this section we assume that we are in this special case. We shall also
assume that $G\ne G^0$. (If $G=G^0$, the classification of objects in $\cn_D^0$ is 
known from \cite{\IC}.) Hence $\kk$ has characteristic $p>1$.
(We could also assume that $G^0$ is simply connected but we will not do so. In fact
we show below that in this case the classification of $\cn_D^0$ is independent of
isogeny.) By an argument in the proof of 12.6 we see that any unipotent 
quasi-semisimple element $u\in G$ such that $\Ad(u):G^0@>>>G^0$ is an inner 
automorphism must be in 
$\cz_G$ hence is $1$. Hence the cyclic group $G/G^0$ is isomorphic to a subgroup of
the group of automorphisms of the Dynkin graph of $G^0$. Hence it has order $p$.
There are four possibilities:

(i) $G^0$ is of type $A_{m-1},m\ge 3,p=2$;

(ii) $G^0$ is of type $D_m,m\ge 4,p=2$;

(iii) $G^0$ is of type $D_4,p=3$;

(iv) $G^0$ of type $E_6,p=2$.
\nl
In each case we have 

(a) ${}^D\cz_{G^0}=\{1\}$.
\nl
Let $G'=G/\cz_{G^0}$ and let $\p:G@>>>G'$ be the obvious map. Let $D'=\p(D)$, a
connected component of $G'$. We show:

(b) {\it if $\boc$ is a unipotent $G'{}^0$-conjugacy class in $D'$ then $\p\i(\boc)$ 
is a unipotent $G^0$-conjugacy class in $D$;}

(c) {\it if $x\in\p\i(\boc)$ then the obvious map $Z_{G^0}(x)@>>>Z_{G'{}^0}(\p(x)$ is 
an isomorphism.}
\nl
Let $x,y\in\p\i(\boc)$. There exists $g\in G^0$ such that $gxg\i=yz$ where 
$z\in\cz_{G^0}$. To prove (b), it is enough to show that $yz$ is $G^0$-conjugate to 
$y$. Now $z'\m y\i z'yz'{}\i$ is an endomorphism of the finite abelian group 
$\cz_{G^0}$ with kernel ${}^D\cz_{G^0}$ which is $\{1\}$. Hence 

(d) {\it this endomorphism is an isomorphism }
\nl
and we can write $z=y\i z'yz'{}\i$ for some $z'\in\cz_{G^0}$. Then $yz=z'yz'{}\i$ and 
(b) is proved.

The surjectivity (resp. injectivity) of the map in (c) follows from (d) (resp. (a)).

From (b),(c) we see that $(\boc,\cf)\m(\p\i(\boc),\p^*\cf)$ is a bijection between 
$\cn_{D'}$ (defined in terms of $G'$) and $\cn_D$ (defined in terms of $G$). It is 
clear that this bijection restricts to a bijection between $\cn_{D'}^0$ and 
$\cn_D^0$.

\subhead 12.8\endsubhead
In this subsection we assume that $\kk$ is an algebraic closure of a finite field 
$\FF_q$ of characteristic $p$, that $G$ has a fixed $\FF_q$-structure with 
Frobenius map $F:G@>>>G$ and that $FD=D$. Let $N'$ be the number of unipotent 
$G^F$-conjugacy classes in $D^F$. In case 12.7(i) and (ii), the author has shown, 
see \cite{\SP, I, 4.5, 4.6}, that $N'$ can be explicitly computed. More precisely, 
in case 12.7(i) we have

$N'=p(m)$, the number of partitions of $m$.
\nl
Using the classification of unipotent representations of a unitary group over a 
finite field in terms of cuspidal ones (as in the proof of \cite{\CLA, 9.2}) we 
deduce
$$N'=\sum_{k\ge 0,s\ge 0;\fra{1}{2}(s^2+s)+2k=m}p_2(k)\tag a$$
where $p_2(k)$ is the number of irreducible representations of a Weyl group of type
$B_k$ (we set $p_0=1$).

In case 12.7(ii), taking $G=O_{2m}(\kk)$, we see that 

(b) {\it the number of unipotent $G^F$-conjugacy classes in $G^F$ is equal to the 
number of irreducible representations of $G^F$ whose restriction to $(G^0)^F$ is a 
sum of unipotent representations},

(c) {\it the number of unipotent $(G^0)^F$-conjugacy classes in $(G^0)^F$ is equal to 
the number of unipotent representations of $(G^0)^F$.}
\nl
Using (b) for an $F$ such that $G^0$ is split over $\FF_q$ we obtain
$$N'+M_1/2+M_2=R_1/2+2R_2.$$
where

$N'$ is the number of unipotent $O_{2m}(\FF_q)$-conjugacy classes in 
$O_{2m}(\FF_q)-SO_{2m}(\FF_q)$ (with $O_{2m}(\FF_q)$ is split),

$M_1$ (resp. $M_2$) is the number of unipotent $SO_{2m}(\FF_q)$-conjugacy classes in 
the split group $SO_{2m}(\FF_q)$ which are not fixed (resp. fixed) by conjugation 
with some/any $g\in O_{2m}(\FF_q)-SO_{2m}(\FF_q)$,

$R_1$ (resp. $R_2$ is the number of unipotent representations of the split group 
$SO_{2m}(\FF_q)$ which do not extend (resp. do extend) to $O_{2m}(\FF_q)$.

Using (c) for an $F$ such that $G^0$ is split over $\FF_q$ we obtain
$$M_1+M_2=R_1+R_2.$$
Using (b) for an $F$ such that $G^0$ is non-split over $\FF_q$ we obtain
$$M=R$$
where

$M$ is the number of unipotent $SO_{2m}(\FF_q)$-conjugacy classes in the non-split 
group $SO_{2m}(\FF_q)$,

$R$ is the number of unipotent representations of the non-split group $SO_{2m}(\FF_q)$.

Using $M_1=R_1,M_2=M$ we see that $N'=R$. Now $R$ can be computed from the 
classification of unipotent representations of the non-split $SO_{2m}(\FF_q)$ in terms 
of cuspidal ones, we deduce (as in (a)); this gives the following formula for $N'$:
$$N'=R=\sum_{k\ge 0,s\ge 0,s\text{ odd };s^2+k=m}p_2(k).\tag d$$
A similar method can be used to compute $N'$ in the cases 12.7(iii),(iv) with $F$ 
such that $G^0$ is split over $\FF_q$. Alternatively, these numbers can be obtained
from \cite{\MA},\cite{\MAA}; they are

(e) $N'=7$ if $G$ is as in 12.7(iii) (with $G^0$ adjoint, split),

(f) $N'=28$ if $G$ is as in 12.7(iv) (with $G^0$ adjoint, split).
\nl
Next we note that $N'$ is equal to the number of unipotent $(G^0)^F$-conjugacy 
classes in $D^F$, since any element $g\in D^F$ is centralized by some element in
$D^F$ (for example by $g$ itself). It follows that $|\cn_D|=N'$. Thus, the method 
above yields $|\cn_D|$ (recall from 12.7 that $|\cn_D|$ is the same for $G$ as for 
$G/\cz_{G^0}$). 

\proclaim{Theorem 12.9} (a) In case 12.7(i), $|\cn_D^0|$ is $1$ if 
$m\in\{3,6,10,\do\}$ and is $0$ otherwise. 

(b) In case 12.7(ii), $|\cn_D^0|$ is $1$ if $m\in\{3^2,5^2,7^2,\do\}$ and is $0$ 
otherwise. 

(c) In cases 12.7(iii) and 12.7(iv), $|\cn_D^0|$ is $1$.
\endproclaim
Assume that we have a (possibly incomplete) list of triples 
$(L^i,\boc_1^i,\cf_1^i)_{i\in\cx}$ in $G^0\bsl\cm_D$ with $L^i\ne L$. For each 
$i\in\cx$, the fibre of the map $\Ph_D:\cn_D@>>>G^0\bsl\cm_D$ (see 8.9) at
$(L^i,\boc_1^i,\cf_1^i)$ is indexed by the irreducible representations of an 
explicit Coxeter group (see Proposition 11.9) hence its cardinal $f_i$ is known.
Then $|\cn_D^0|\le|\cn_D|-\sum_{i\in X}f_i$. Here the right hand side is explicitly
known since $|\cn_D|$ is known from 12.8. If $|\cn_D|-\sum_{i\in X}f_i=0$, then it 
follows that $|\cn_D^0|=0$. If $|\cn_D|-\sum_{i\in X}f_i=1$, then it follows that 
our list is complete (any additional member of that list would contribute at least
$2$ to $|\cn_D|$ since the corresponding Coxeter group (see 11.9) is non-trivial); 
it follows that in this case $|\cn_D|=1$. This method works in each case. We give 
in each case the list $(L^i,\boc_1^i,\cf_1^i)$.

In case 12.7(i) with $G^0=PGL_m(\kk)$, we take $L^i$ to be the image of 
$(\kk^*)^{2k}\T GL_{(s^2+s)/2}$ under $GL_m(\kk)@>>>PGL_m(\kk)$ (here 
$m=\fra{1}{2}(s^2+s)+2k,k>0$) and $(\boc_1^i,\cf_1^i)$ is uniquely determined by 
$L^i$ (we use an inductive hypothesis for $L^i/\cz^0_{L^i}$). The required formula 
for $|\cn_D^0|$ is equivalent to
$$|\cn_D|=\sum_{k\ge 0,s\ge 0;\fra{1}{2}(s^2+s)+2k=m}p_2(k)\tag d$$
which is known from 12.8.

In case 12.7(ii) with $G^0=SO_{2m}$, we take $L^i$ to be $(\kk^*)^{2k}\T SO_{2s^2}$
(here $m=s^2+k$, $s$ odd, $k>0$) and $(\boc_1^i,\cf_1^i)$ is uniquely determined by
$L^i$ (we use an inductive hypothesis for $L^i/\cz^0_{L^i}$). The required formula 
for $|\cn_D^0|$ is equivalent to 
$$|\cn_D|=\sum_{k\ge 0,s\ge 0,s\text{ odd };s^2+k=m}p_2(k)\tag e$$
which is known from 12.8.

In case 12.7(iii), we take $L^i$ to be a maximal torus of $G^0$ and 
$(\boc_1^i,\cf_1^i)$ is uniquely determined by $L^i$. The required formula for 
$|\cn_D^0|$ is equivalent to $|\cn_D|=1+6$ which is known from 12.8. Here $6$ is the 
number of irreducible representations of a Weyl group of type $G_2$.

In case 12.7(iv), we take $L^i$ to be such that $L^i/\cz^0_{L^i}$ has type $A_5$ or 
$L^i$ is a maximal torus of $G^0$; $(\boc_1^i,\cf_1^i)$ is uniquely determined by 
$L^i$ (we use (a) which is already proved). The required formula for $|cn_D^0|$ is 
equivalent to $|\cn_D|=1+2+25$ which is known from 12.8. Here $2$ (resp. $25$) is the 
number of irreducible representations of a Weyl group of type $A_1$ (resp. $F_4$). The 
theorem is proved.

\head 13. Symbols\endhead
\subhead 13.1\endsubhead
Symbols are combinatorial objects used in \cite{\CLA} to parametrize unipotent
representations of classical groups over a finite field and in \cite{\IC} to
describe the generalized Springer correspondence for classical groups in
characteristic $\ne 2$. In \cite{\LS} it has been observed that symbols can also be
used to describe the generalized Springer correspondence for (connected) classical 
groups in characteristic $2$. Since the combinatorics of unipotent classes in
disconnected classical groups in characteristic $2$ is very similar to that in the
connected case, it can be expected that in this case, again the generalized 
Springer correspondence can be described in terms of symbols.

\subhead 13.2\endsubhead
Let $\r,s\in\NN,n\in\ZZ$. For any $d\in\ZZ$ let ${}^\r\tX^s_{n,d}$ be the set of 
all ordered pairs $(A;B)$ of finite sequences of natural numbers 
$A:a_1,a_2,\do,a_m$ and $B:b_1,b_2,\do,b_{m'}$ (for some $m,m'$) that are subject 
to the following conditions:

$a_i-a_{i-1}\ge\r$ for $1<i\le m$;

$b_i-b_{i-1}\ge\r$ for $1<i\le m'$;

$b_i\ge s$ for all $1\le i\le m'$;

$m-m'=d$;

$\sum a_i+\sum b_i=n+\r(m+m')(m+m'-2)/4+s(m+m')/2$ if $d$ is even;

$\sum a_i+\sum b_i=n+\r(m+m'-1)^2/4+s(m+m'-1)/2$ if $d$ is odd.
\nl
(In \cite{\LS,\S1} the notation $\tX^{r,s}_{n,d}$ is used for the present
${}^\r\tX_{n,d}^s$ where $\r=r+s$.)

Let ${}^\r X^s_{n,d}$ be the set of equivalence classes on ${}^\r\tX^s_{n,d}$ for 
the equivalence relation generated by 
$$\align&(a_1,a_2,\do,a_m;b_1,b_2,\do,b_{m'})\si
\\&(0,a_1+\r,a_2+\r,\do,a_m+\r;s,b_1+\r,b_2+\r,\do,b_{m'}+\r).\endalign$$
For example, ${}^0X^0_{n,d}$ is in obvious bijection with the set of pairs of 
partitions $\a,\b$ with $\sum\a_i+\sum\b_i=n$; hence we have a bijection
$${}^0X^0_{n,d}\lra\Irr W_n\tag a$$
where $W_n$ is a Coxeter group of type $B_n$. (By convention, $W_0=\{1\}$ and 
$\Irr W_n=\em$ for $n<0$.) Let 

$n_{\r,s,d}=\r d^2/4-sd/2$ for $d$ even, $n_{\r,s,d}=\r(d-1)(d+1)/4-s(d-1)/2$ for $d$ 
odd.
\nl
We have a bijection 
${}^0X^0_{n-n_{\r,s,d},d}@>>>{}^\r X^s_{n,d}$ given by
$$\align&(c_1,c_2,\do,c_m;c'_1,c'_2,\do,c'_{m'})\m \\&
(c_1,c_2+\r,\do,c_m+(m-1)\r;c'_1+s,c'_2+s+\r,\do,c'_{m'}+s+(m'-1)\r).\endalign$$
Combined with the bijection (a) we obtain a bijection
$${}^\r X^s_{n,d}\lra\Irr W_{n-n_{\r,s,d}}.\tag b$$  
Now let $E$ be $2\ZZ$ or $2\ZZ+1$ (if $s>0$) and let $E=2\NN+1$ (if $s=0$). Let 
${}^\r X^s_{n,E}=\sqc_{d\in E}{}^\r X^s_{n,d}$. From (b) we obtain a bijection 
$${}^\r X^s_{n,E}\lra\sqc_{d\in E}\Irr W_{n-n_{\r,s,d}}.\tag c$$
Assume that $\r\ge 1$. If $(A;B)\in{}^\r\tX^s_{n,E}$ then $A,B$ may be considered as
subsets of $\NN$. Consider two elements of ${}^\r X^s_{n,E}$; we can represent them
in the form $(A;B),(A';B')$ where $|A|+|B|=|A'|+|B'|$ (by our assumption on $E$). We
say that these two elements are {\it similar} if $A\cup B=A'\cup B'$,
$A\cap B=A'\cap B'$. This defines an equivalence relation (similarity) on 
${}^\r X^s_{n,E}$. 

Assume that $\r>s$. Let $(A;B)\in{}^\r\tX^s_{n,d}$. Let $S=(A\cup B)-(A\cap B)$. A 
non-empty subset $I$ of $S$ is said to be an interval if 
$I=\{c_1<c_2<\do<c_k\}$ with $c_2-c_1<\r,c_3-c_2<\r,\do,c_k-c_{k-1}<\r$ and $I$ is
maximal with this property. Clearly, the intervals form a partition of $S$. An 
interval is said to be proper if it does not contain any number in $[0,s-1]$.

Let $\cc$ be a similarity class in ${}^\r X^s_{n,E}$ where $\r>s$. Let $V'_\cc$ be 
the $\FF_2$-vector space with basis indexed by the proper intervals in 
$S=(A\cup B)-(A\cap B)$ where $(A;B)$ represents an element in $\cc$. If $s>0$ let 
$V_\cc=V'_\cc$; if $s=0$ let $V_\cc$ be the quotient of $V'_\cc$ by the subspace 
spanned by the sum of all basis elements of $V'_\cc$. Note that $V_\cc$ is 
independent of the choice of $(A;B)$. As in \cite{\IC, 11.5}, \cite{\LS, 1.4} we see
that there is a canonical bijection $\cc\lra V_\cc$.
(In particular, $\cc$ has a natural structure of $\FF_2$-vector space.) Putting 
together the bijections above we obtain a bijection
$${}^\r X^s_{n,E}\lra\sqc_\cc V_\cc\tag d$$
where $\cc$ runs over the set of similarity classes in ${}^\r X^s_{n,E}$.

Note that the canonical basis of $V'_\cc$ is totally ordered by the requirement that
the basis element corresponding to a proper interval $I$ is $<$ than the basis 
element corresponding to a proper interval $I'$ if any number in $I$ is $<$ than 
any number in $I'$.

\medpagebreak

We describe the partition of ${}^\r X^s_{n,E}$ in similarity classes in a number of
cases. (In each case, each horizontal line contains the elements in a similarity 
class.)
 
${}^4X^1_{1,2\ZZ+1}$:

$(1;\em),(\em;1)$

$(0,4;2)$

\medpagebreak

${}^4X^1_{2,2\ZZ}$:

$(0;3)$

$(1;2),(2;1)$

$(0,4;2,6)$

$(1,5;1,5)$.

\medpagebreak

${}^4X^1_{2,2\ZZ+1}$:

$(2;\em),(\em,2)$,

$(0,5;2)$

$(1,5;1),(1;1,5)$

$(0,4;3)$

$(0,4,8;2,6)$

\medpagebreak

${}^4 X^0_{3,2\NN+1}$:

$(3,\em)$

$(0,6;1),(1,6;0)$

$(0,5;2)$

$(1,5;1)$

$(0,4;3)$

$(0,4,9;1,5),(1,5,9;0,4)$

$(0,4,8;1,6)$

$(0,4,8,12;1,5,9)$.

\subhead 13.3\endsubhead   
Let $\cv_{2n}$ be the set of all pairs $(\l,\da)$ where $\l$ is a sequence 
$\l_1\le\l_2\le\do\le\l_{2k+1}$ in $\NN$ with $\sum_i\l_i=2n$ and $\da$ is a 
partition of $[1,2k+1]$ into one and two element subsets called blocks (where each 
two element block consists of consecutive integers) such that 

if $\{i\}$ is a block then $\l_i$ is even;

if $\{i,i+1\}$ is a block then $\l_i=\l_{i+1}$;

if $\{i\}$ and $\{j,j+1\}$ are blocks then $\l_i\ne \l_j$;
\nl
also, at most one of the $\l_i$ is $0$. (This last condition determines uniquely $k$.)

Given $(\l,\da)\in\cv_{2n}$ as above we define a sequence 
$c_1\le c_2\le\do\le c_{2k+1}$ by

$c_i=\l_i/2+2(i-1)$ if $\{i\}$ is a block;               

$c_i=(\l_i+1)/2+2(i-1)$, $c_{i+1}=(\l_{i+1}-1)/2+2i=c_i+1$ if $\{i,i+1\}$ is a block
with $\l_i=\l_{i+1}=$ odd;

$c_i=(\l_i+2)/2+2(i-1)$, $c_{i+1}=(\l_{i+1}-2)/2+2i=c_i$  if $\{i,i+1\}$ is a block 
with $\l_i=\l_{i+1}=$ even.

Let $A_{\l,\da}$ be the $\FF_2$-vector space generated by the set $\{s_i\}$ where 
$i\in[1,2k+1]$ is such that either $\{i\}$ is a block or $\l_i$ is odd; the 
relations are:

$s_i=s_j$ if $\l_i=\l_j$;

$s_i=s_j$ if $\l_j=\l_i+1$;

$s_i=s_j$ if $\l_i,\l_j$ are even and $\l_j=\l_i+2$;

$s_i=0$  if $\l_i\le 2$.
                      
\subhead 13.4\endsubhead
Let $\cv'_{2n}$ be the subset of $\cv_{2n}$ defined by the condition that $\l_i>0$ 
for all $i$.

Given $(\l,\da)\in\cv'_{2n}$ as above we define a sequence 
$c'_1\le c'_2\le\do\le c'_{2k+1}$ by $c'_i=c_i-1$ (where $c_i$ is as in 13.3.)

Let $\tA'_{\l,\da}$ be the $\FF_2$-vector space generated by the set $\{s_i\}$ where
$i\in[1,2k+1]$ is such that either $\{i\}$ is a block or $\l_i$ is odd; the 
relations are:

$s_i=s_j$ if $\l_i=\l_j$;

$s_i=s_j$ if $\l_j=\l_i+1$;

$s_i=s_j$ if $\l_i,\l_j$ are even and $\l_j=\l_i+2$.

Let $A'_{\l,\da}$ be the subspace of $\tA'_{\l,*}$ generated by the sums $s_i+s_j$ 
where $i\ne j$ and both $s_i,s_j$ are defined.

\subhead 13.5\endsubhead
Consider the set $\cv''_N$ of all pairs $(\l,\da)$ where $\l$ is a sequence 
$\l_1\le\l_2\le\do\le\l_k$ in $\NN+1$ with $\sum_i\l_i=N$ and $\da$ is a partition of
$[1,k]$ into one and two element subsets called blocks (where each two element block 
consists of consecutive integers) such that 

if $\{i\}$ is a block then $\l_i$ is odd;

if $\{i,i+1\}$ is a block then $\l_i=\l_{i+1}$;

if $\{i\}$ and $\{j,j+1\}$ are blocks then $\l_i\ne \l_j$.
\nl
Given $(\l,\da)\in\cv''_N$ as above we define a sequence 
$c''_1\le c''_2\le\do\le c''_k$ by

$c''_i=(\l_i-1)/2+2(i-1)$ if $\{i\}$ is a block;               

$c''_i=\l_i/2+2(i-1)$, $c''_{i+1}=(\l_{i+1}-2)/2+2i=c_i+1$ if $\{i,i+1\}$ is a block
with $\l_i=\l_{i+1}=$ even;

$c''_i=(\l_i+1)/2+2(i-1)$, $c_{i+1}=(\l_{i+1}-3)/2+2i=c_i$  if $\{i,i+1\}$ is a block
with $\l_i=\l_{i+1}=$ odd.
\nl
In other words, $c''_i$ are defined by the same formulas as the $c_i$ in 13.3 but 
with $\l_i$ replaced by $(\l_i-1)$.

Let $A''_{\l,\da}$ be the $\FF_2$-vector space generated by the set $\{s_i\}$ where 
$i\in[1,k]$ is such that either $\{i\}$ is a block or $\l_i$ is even; the relations
are:

$s_i=s_j$ if $\l_i=\l_j$;

$s_i=s_j$ if $\l_j=\l_i+1$;

$s_i=s_j$ if $\l_i,\l_j$ are odd and $\l_j=\l_i+2$;

$s_i=0$  if $\l_i=1$.

\subhead 13.6\endsubhead
Note that $A_{\l,\da},\tA'_{\l,\da},A''_{\l,\da}$ have canonical bases consisting 
of the images of those $s_i$ that are non-zero. These bases are totally ordered by 
the requirement that the basis element represented by some $s_i$ is $\le$ than the 
basis element represented by some $s_j$ if $i<j$. We define maps
$$\sqc_{(\l,\da)\in\cv_{2n}}A_{\l,\da}^*@>>>\sqc_{\cc\sub{}^4X^2_{n,2\ZZ+1}}V_\cc
\tag a$$
$$\sqc_{(\l,\da)\in\cv'_{2n}}A'_{\l,\da}{}^*@>>>
                                     \sqc_{\cc\sub{}^4X^0_{n-1,2\NN+1}}V_\cc\tag b$$
$$\sqc_{(\l,\da)\in\cv''_{2n+1}}A''_{\l,\da}{}^*@>>>
                                      \sqc_{\cc\sub{}^4X^1_{n,2\ZZ+1}}V_\cc\tag c$$
$$\sqc_{(\l,\da)\in\cv''_{2n}}A''_{\l,\da}{}^*@>>>
\sqc_{\cc\sub{}^4X^1_{n,2\ZZ}}V_\cc\tag d$$
where $\cc$ runs over the set of similarity classes in the appropriate 
${}^4X^s_{N,E}$ and ${}^*$ denotes the dual $\FF_2$-vector space, as follows.

In (a) any $(\l,\da)$ determines as in 13.3 a sequence
$c_1\le c_2\le\do\le c_{2k+1}$. Then $\cc$ is the similarity class of
$(c_1,c_3,c_5,\do,c_{2k+1};c_2,c_4,\do,c_{2k})$ in ${}^4X^2_{n,2\ZZ+1}$. There is a
unique $\FF_2$-vector space isomorphism $A_{\l,\da}^*@>\si>>V_\cc$ which preserves 
the canonical bases and their natural total order; this defines the map (a).

In (b) any $(\l,\da)$ determines as in 13.4 a sequence
$c'_1\le c'_2\le\do\le c'_{2k+1}$. Then $\cc$ is the similarity class of
$(c'_1,c'_3,c'_5,\do,c'_{2k+1};c'_2,c'_4,\do,c'_{2k})$ in ${}^4X^0_{n,2\NN+1}$.
There is a unique $\FF_2$-vector space isomorphism $\tA'_{\l,\da}{}^*@>\si>>V'_\cc$ 
which preserves the canonical bases and their natural total order; by passage to a 
quotient, this induces an isomorphism $A'_{\l,\da}{}^*@>\si>>V_\cc$ and defines the
map (b).

In (c) and (d) any $(\l,\da)$ determines as in 13.5 a sequence
$c''_1\le c''_2\le\do\le c''_k$ where $k$ is odd in (c) and is even in (d). Then 
$\cc$ is the similarity class of $(c'_1,c'_3,\do;c'_2,c'_4,\do)$ in 
${}^4X^1_{n,E}$ where $E$ is $2\ZZ+1$ in (c) and is $2\ZZ$ in (d). There is a 
unique $\FF_2$-vector space isomorphism $A''_{\l,\da}{}^*@>\si>>V_\cc$ which preserves
the canonical bases and their natural total order; this defines the maps (c),(d).

The maps (a),(b),(c),(d) are well defined bijections. For (a) this is shown in
\cite{\LS, 2.2}; exactly the same proof applies for (b),(c),(d). (A partial result
in this direction can be found in \cite{\MS}.)

\subhead 13.7\endsubhead
Combining the bijection 13.2(d) with the bijection (a),(b),(c) or (d) in 13.6, we 
obtain bijections

$\sqc_{(\l,\da)\in\cv_{2n}}A_{\l,\da}^*@>\si>>{}^4X^2_{n,2\ZZ+1}$

$\sqc_{(\l,\da)\in\cv'_{2n}}A'_{\l,\da}{}^*@>\si>>{}^4X^0_{n-1,2\NN+1}$
                                          
$\sqc_{(\l,\da)\in\cv''_{2n+1}}A''_{\l,\da}{}^*@>\si>>{}^4X^1_{n,2\ZZ+1}$

$\sqc_{(\l,\da)\in\cv''_{2n}}A''_{\l,\da}{}^*@>\si>>{}^4X^1_{n,2\ZZ}$.

\subhead 13.8\endsubhead
Assume now that $\kk$ has characteristic $2$. Let $G=G_m$ be as in 12.7(i). If
$G,D$ is one of 
$$\align&(Sp_{2n}(\kk),Sp_{2n}(\kk)),(O_{2n}(\kk),O_{2n}(\kk)-SO_{2n}(\kk)),\\&
(G_{2n+1},G_{2n+1}-G_{2n+1}^0),(G_{2n},G_{2n}-G_{2n}^0),\tag a\endalign$$
the generalized Springer correspondence can be viewed as a bijection

$\cn_D\lra\sqc_{d\in 2\ZZ+1}\Irr W_{n-d(d-1)}$

$\cn_D\lra\sqc_{d\in 2\NN+1}\Irr W_{n-d^2}$

$\cn_D\lra\sqc_{s\in\NN;s(s+1)/2\text{ odd}}\Irr W_{((2n+1)-s(s+1)/2)/2}=
                 \sqc_{d\in 2\ZZ+1}\Irr W_{n-(d^2-1-\fra{1}{2}(d-1))}$

$\cn_D\lra\sqc_{s\in\NN;s(s+1)/2\text{ even}}\Irr W_{(2n-s(s+1)/2)/2}=
                                 \sqc_{d\in 2\ZZ}\Irr W_{n-(d^2-\fra{1}{2}d)}$
\nl
respectively.

\subhead 13.9\endsubhead
According to \cite{\SP, I, 2.6, 2.9}, for $(G,D)$ as in 13.8(a), the set $\cn_D$ is 
naturally in bijection with 
$$\sqc_{(\l,\da)\in\cv_{2n}}A_{\l,\da}^*,\qua
\sqc_{(\l,\da)\in\cv'_{2n}}A'_{\l,\da}{}^*,\qua
\sqc_{(\l,\da)\in\cv''_{2n+1}}A''_{\l,\da}{}^*,\qua
\sqc_{(\l,\da)\in\cv''_{2n}}A''_{\l,\da}{}^*,$$
respectively. (In \cite{\SP, I, 2.11} the notation for unipotent classes is in terms 
of a partition in which to certain parts of a fixed size one attaches an index $0$.
The partition is the first coordinate in our $(\l,\da)$; when the index $0$ is 
attached to the parts of size $a$, then the parts of size $a$ form blocks of size 
$2$ according to $\da$. For example $1^2\op 2\op 4^2_0$ (resp. $1^2\op 2\op 4^2$) in
\cite{\SP} becomes our $(11)(2)(44)$ (resp. $(11)(2)(4)(4)$) where the brackets 
represent the partition $\da$.)

Composing the bijections above with the bijection in 13.7 we see that $\cn_D$ is 
naturally in bijection with 
$${}^4X^2_{n,2\ZZ+1},{}^4X^0_{n-1,2\NN+1},{}^4X^1_{n,2\ZZ+1},{}^4X^1_{n,2\ZZ}$$
respectively. Combining this with the bijections 13.2(c) we obtain bijections

$\cn_D\lra\sqc_{d\in 2\ZZ+1}\Irr W_{n-(d^2-d)}$,

$\cn_D\lra\sqc_{d\in 2\NN+1}\Irr W_{n-1-(d^2-1)}$,

$\cn_D\lra\sqc_{d\in 2\ZZ+1}\Irr W_{n-(d^2-1-\fra{1}{2}(d-1))}$,

$\cn_D\lra\sqc_{d\in 2\ZZ}\Irr W_{n-(d^2-\fra{1}{2}d)}$,
\nl
respectively.

\proclaim{Proposition 13.10} The previous four bijections coincide with the 
respective bijections given by the generalized Springer correspondence (see 13.8).
\endproclaim
(This shows that the generalized Springer correspondence is explicitly computable.)
When $(G,D)=(Sp_{2n}(\kk),Sp_{2n}(\kk))$, this is shown in \cite{\LS, 2.4} using the
restriction formula \cite{\IC, 8.3}. The proof in the other three cases is entirely
similar, once one has the analogue of the restriction formula for disconnected 
groups. That analogue is given by Proposition 9.4. (In these last three cases, a 
special case of the proposition, namely the part pertaining to pairs 
$(\boc,\cf)\in\cn_D$ with $\cf=\bbq$, is proved in \cite{\MS}.)

\proclaim{Corollary 13.11}(a) If $G,D$ are as in 12.7 with $G^0$ of type $A_{m-1}$ 
and $(\boc,\cf)\in\cn_D^0$ then the partition corresponding to $\boc$ is 
$m=3+7+11+\do$ or $m=1+5+9+\do$.

(b) If $G,D$ are as in 12.7 with $G=O_{2n}(\kk)$ and $(\boc,\cf)\in\cn_D^0$ then the
partition corresponding to $\boc$ is $2n=2+6+10+\do$ (an odd number of parts).
\endproclaim

\head 14. Spin groups\endhead
\subhead 14.1\endsubhead
In this section we assume that $\kk$ has characteristic $\ne 2$ and that $G=G^0$ is
$\text{\rm Spin}_n(\kk),n\ge 3$. We have a partition $\cn_G=\sqc_\c\cn_G^\c$ where 
$\c$ runs over the characters $\cz_G@>>>\kk^*$ and $\cn_G^\c$ consists of all 
$(\boc,\cf)\in\cn_G$ such that the $\cz_G$-action on $\cf$ (coming from the 
$G$-equivariance of $\cf$) is via $\c$ on each fibre of $\cf$. Assume now that $\c$
is such that its restriction to the kernel of the obvious isogeny 
$\text{\rm Spin}_n(\kk)@>>>SO_n(\kk)$ is non-trivial. According to \cite{\IC, \S14}, 
the generalized Springer correspondence for $G$ restricts to a bijection
$$\cn_G^\c\lra\sqc_{t\in 4\ZZ+n}\Irr W_{\fra{1}{4}(n-(2t^2-t))}.$$
Moreover, the map which to each $(\boc,\cf)\in\cn_G^\c$ associates the partition of 
$n$ whose parts are the sizes of the Jordan blocks of the image in $SO_n(\kk)$ of an
element in $\boc$ is a bijection
$$\cn_G^\c\lra X_n$$
where $X_n$ is the set of all partitions $\l=(\l_1\le\l_2\le\do\le\l_m)$ of $n$ 
(with $\l_i\in\NN+1$) such that

(a) for any even $p$, $|\{i|\l_i=p\}|$ is even;

(b) for any odd $p$ we have $|\{i|\l_i=p\}|\le 1$.
\nl
Thus the generalized Springer correspondence restricted to $\cn_G^\c$ may be regarded
as a bijection
$$X_n\lra\sqc_{t\in 4\ZZ+n}\Irr W_{\fra{1}{4}(n-(2t^2-t))}.\tag c$$
We would like to describe this bijection in a combinatorial way. In principle, the 
results in \cite{\LS, \S4} provide such a description, which involves an inductive
procedure instead of a closed formula. In this section we provide a closed formula 
for (c). 

\subhead 14.2\endsubhead
Let $\l=(\l_1\le\l_2\le\l_3\le\dots\le\l_m)\in X_n$. For $i\in[1,m]$ we set 
$$t_i=\sum_{j=1}^{i-1}d(\l_j)$$ 
where for any integer $s$ we set 
$$\text{$d(s)=0$ if $s$ is even, $d(s)=(-1)^{(s-1)/2}$ if $s$ is odd}.$$
We modify the entries of $\l$ as follows and we mark the modified entries by an 
indeterminate $a$ or $b$ as follows:

(1) if an entry $e=\l_i$ satisfies $e\in 4\ZZ+1$ then $e$ is replaced by 
$\fra{1}{4}(e-1)-t_i$ with mark $a$;

(2) if an entry $e=\l_i$ satisfies $e\in 4\ZZ+3$, then $e$ is replaced by 
$\fra{1}{4}(e-3)+t_i$ with mark $b$;

(3) if an entry $e\in 4\ZZ$ appears exactly $2p>0$ times with 
$\l_i=\l_{i+1}=\do=\l_{i+2p-1}=e$, we replace $e,e,\do,e$ ($2p$ terms) by
$$\fra{1}{4}e-t_i,\fra{1}{4}e+t_i,\fra{1}{4}e-t_i,\fra{1}{4}e+t_i,\do,
\fra{1}{4}e-t_i,\fra{1}{4}e+t_i$$
($2p$ terms, marked by $a,b,a,b,\do,a,b$);

(4) if an entry $e\in 4\ZZ+2$ appears exactly $2p>0$ times with 
$\l_i=\l_{i+1}=\do=\l_{i+2p-1}=e$, we replace $e,e,\do,e$ ($2p$ terms) by
$$\fra{1}{4}(e+2)-t_i,\fra{1}{4}(e-2)+t_i,\fra{1}{4}(e+2)-t_i,
\fra{1}{4}(e-2)+t_i,\do,\fra{1}{4}(e+2)-t_i,\fra{1}{4}(e-2)+t_i$$
($2p$ terms, marked by $a,b,a,b,\do,a,b$).

\proclaim{Lemma 14.3} (a) The modified entries (marked with $a$) form an increasing
sequence $\a$ in $\NN$.

(b) The modified entries (marked with $b$) form an increasing sequence $\b$ in $\NN$.
\endproclaim
We prove (a). Consider two consecutive (modified) entries $\l'_i,\l'_j$ marked with
$a$, coming from $\l_i,\l_j$. We show that $\l'_i\le\l'_j$.

If they both come from a group as in (3) or (4), then there is nothing to prove. 

If $\l'_i$ comes from a group as in (3) (resp. (4)) and $\l'_j$ comes from another 
group as in (3) (resp. (4)), the result is clear. 

If $\l'_i$ comes from a group as in (3) and $\l'_j$ comes from a group as in (4), we
have $\l_i/4<(\l_j+2)/4$ since $\l_i<\l_j$.

If $\l'_i$ comes from a group as in (4) and $\l'_j$ comes from a group as in (3), we
have $(\l_i+2)/4\le\l_j/4$ since $\l_i\le\l_j-2$.

If $\l'_i$ comes from an entry as in (1) and $\l'_j$ comes from a group as in (3),
we have $(\l_i-1)/4\le\l_j/4-1$ since $\l_i\le\l_j-3$.

If $\l'_i$ comes from an entry as in (1) and $\l'_j$ comes from a group as in (4),
we have $(\l_i-1)/4\le(\l_j+2)/4-1$ since $\l_i\le\l_j-1$.

If $\l'_i$ comes from a group as in (3) and $\l'_j$ comes from an entry as in (1),
we have $\l_i/4\le(\l_j-1)/4$ since  $\l_i\le\l_j-1$.

If $\l'_i$ comes from a group as in (4) and $\l'_j$ comes from an entry as in (1),
we have $(\l_i+2)/4\le(\l_j-1)/4$ since  $\l_i\le\l_j-3$.

If $\l'_i,\l'_j$ come from entries as in (1) we have $(\l_i-1)/4\le(\l_j-1)/4-1$ 
since $\l_i\le\l_j-4$.

We now consider the first (modified) entry $\l'_i$ of type $a$, coming from $\l_i$.
We show that $\l'_i\ge 0$.

If it is of type (1), it is $(\l_i-1)/4-(-1-1-\do-1)\ge(\l_i-1)/4\ge 0$.

If it is of type (3), it is $\l_i/4-(-1-1-\do-1)\ge\l_i/4\ge 0$.

If it is of type (4), it is $(\l_i+2)/4-(-1-1-\do-1)\ge\l_i/4\ge 0$.

We prove (b). Consider two consecutive (modified) entries $\l'_i,\l'_j$ marked with
$b$, coming from $\l_i,\l_j$. We show that $\l'_i\le\l'_j$.

If they both come from a group as in (3) or (4), then there is nothing to prove. 

If $\l'_i$ comes from a group as in (3) (resp. (4)) and $\l'_j$ comes from another 
group as in (3) (resp. (4)), the result is clear. 

If $\l'_i$ comes from a group as in (3) and $\l'_j$ comes from a group as in (4), 
we have $\l_i/4\le(\l_j-2)/4$ since  $\l_i\le\l_j-2$.

If $\l'_i$ comes from a group as in (4) and $\l'_j$ comes from a group as in (3), we
have $(\l_i-2)/4\le\l_j/4$ since  $\l_i<\l_j$.

If $\l'_i$ comes from an entry as in (2) and $\l'_j$ comes from a group as in (3),
we have $(\l_i-3)/4\le\l_j/4-1$ since $\l_i\le\l_j-1$.

If $\l'_i$ comes from an entry as in (2) and $\l'_j$ comes from a group as in (4),
we have $(\l_i-3)/4\le(\l_j-2)/4-1$ since $\l_i\le\l_j-3$.

If $\l'_i$ comes from a group as in (3) and $\l'_j$ comes from an entry as in (2),
we have $\l_i/4\le(\l_j-3)/4$ since  $\l_i\le\l_j-3$.

If $\l'_i$ comes from a group as in (4) and $\l'_j$ comes from an entry as in (2),
we have $(\l_i-2)/4\le(\l_j-3)/4$ since  $\l_i\le\l_j-1$.

If $\l'_i,\l'_j$ come from entries as in (2) we have $(\l_i-3)/4\le(\l_j-3)/4-1$ 
since $\l_i\le\l_j-4$.

We now consider the first (modified) entry $\l'_i$ of type $b$, coming from $\l_i$.
We show that $\l'_i\ge 0$.

If it is of type (2), it is $(\l_i-3)/4+(1+1+\do+1)\ge(\l_i-3)/4\ge 0$.

If it is of type (3), it is $\l_i/4+(1+1+\do+1)\ge\l_i/4\ge 0$.

If it is of type (4), it is $(\l_i-2)/4+(1+1+\do+1)\ge(\l_i-2)/4\ge 0$.

The lemma is proved.

\proclaim{Lemma 14.4} The sum $S$ of the modified entries is $(n-2t^2+t)/4$ where 
$t=\sum_id(\l_i)$.
\endproclaim
We have
$$\align S&=\fra{1}{4}\sum_{i;\l_i\in 2\ZZ}\l_i+\fra{1}{4}\sum_{i;\l_i\in 2\ZZ+1}\l_i
-\fra{1}{4}|\{i;\l_i\in 4\ZZ+1\}|\\&-\fra{3}{4}|\{i;\l_i\in 4\ZZ+3\}|
+\sum_{i;\l_i\in 4\ZZ+3}t_i-\sum_{i;\l_i\in 4\ZZ+1}t_i\endalign$$
hence $S=n/4+k$ where
$$\align k&=-\sum_it_id(\l_i)+\sum_i(d(\l_i)-2d(\l_i)^2)/4\\&
=-\sum_{j<i}d(\l_i)d(\l_j)-\sum_id(\l_i)^2/2+\sum_id(\l_i)/4\\&
=-(\sum_id(\l_i))^2/2+\sum_id(\l_i)/4=-t^2/2+t/4=(-2t^2+t)/4.\endalign$$
The lemma is proved.

\subhead 14.5\endsubhead
To $\l$ we associate the ordered pair of partitions $(\a,\b)$ (if $t\ge 1$) or
$(\b,\a)$ (if $t\ge 0$). This pair of partition may be regarded as an element of 
$\Irr W_{(n-2t^2+t)/4}$. Since $(n-2t^2+t)/4\in\ZZ$ we have $t\in 4\ZZ+n$. We have 
thus defined in a combinatorial way a map 

$X_n@>>>\sqc_{t\in 4\ZZ+n}\Irr W_{\fra{1}{4}(n-(2t^2-t))}$.
\nl
From the definitions one can check that this coincides with the map defined in 
\cite{\LS,\S4} by an inductive procedure. It therefore coincides with the 
generalized Springer correspondence 14.1(c).

\enddocument